\documentclass{amsart}

\usepackage{graphicx} 
\usepackage{amsmath,amsfonts,amssymb,amsthm}
\usepackage{algorithmic}
\usepackage{latexsym, graphicx, epsfig}
\usepackage{verbatim}
\usepackage[linesnumbered,vlined,commentsnumbered,ruled]{algorithm2e}
\usepackage{pgfplots}
\usepackage{url}
\usepackage{mathrsfs}
\usepackage[a4paper]{geometry}
\usepackage{color}
\usepackage{array}
\usepackage{subfig}
\newtheorem{remark}{Remark}[section]

\usepackage{bm}
\usepackage{float}
\usepackage[colorlinks,linkcolor=blue,anchorcolor=blue,citecolor=red]{hyperref}
\usepackage{booktabs}

\SetKwInOut{Offline}{Offline}
\SetKwInOut{Online}{Online}

\newcommand{\z}{\mathbf{z}}
\newcommand{\ud}{\, \mathrm{d}}

\title[Multi-fidelity for kinetic models with uncertain contact dynamics]{Multi-fidelity methods for kinetic models of epidemic dynamics with uncertain contact structure}

\date{}

\author{Liu Liu}
\address{The Chinese University of Hong Kong, Hong Kong}
\email{liuliu@cuhk.edu.hk}

\author{Andrea Medaglia}
\address{Department of Mathematics and Computer Science, University of Ferrara, Italy}
\email{andrea.medaglia@unife.it}

\author{Hao Xie}
\address{The Chinese University of Hong Kong, Hong Kong}
\email{haoxie@link.cuhk.edu.hk}

\author{Mattia Zanella}
\address{Department of Mathematics "F. Casorati", University of Pavia, Italy}
\email{mattia.zanella@unipv.it}

\begin{document}

\begin{abstract}
In this work, we develop a multi-fidelity strategy for kinetic models in epidemiology with uncertain contact dynamics. Assessing and controlling the population-level effects of contact dynamics requires the development of models for understanding observable effects of heterogeneous contact structures, whose formation depends on complex social phenomena. These can be captured taking into account high-dimensional uncertain quantities. The proposed approach combines high-fidelity kinetic solvers with a hierarchy of low-fidelity surrogates, including reduced macroscopic models and coarse kinetic descriptions, remaining applicable even in regimes where a macroscopic closure is unavailable. This hierarchical framework identifies representative parameter samples and reconstructs full solutions via projection-based techniques, enabling efficient uncertainty propagation while drastically reducing computational cost. Numerical experiments in high-dimensional stochastic settings demonstrate that accurate statistical estimates of epidemic observables can be obtained with significantly reduced computational costs compared to standard approaches.

\noindent\textbf{MSC:} 35Q84, 92D30, 65C20, 35Q92 \\
\textbf{Keywords:} Contact heterogeneity; Kinetic modelling; mathematical epidemiology; multi-fidelity methods; uncertainty quantification 
\end{abstract}

\maketitle

\section{Introduction}\label{Introduction}

Recent epidemics have highlighted that the spread of infectious diseases is not only determined by classical mechanisms of compartmental epidemiology, but also by the social contact structure of the population, which is often heterogeneous, poorly measured and variable in time~\cite{Albi2021,Albi2021.b,DPTZ,Gatto2020}. In addition, changes in the contact structure can have a significant impact on disease dynamics \cite{Sun21}. Variations in the average number of contacts, as well as in their distribution across individuals, may alter transmission pathways \cite{beraud15,britton20}. The computation of human mixing patterns has been approached through various methodologies, including surveys, contact diaries, sensors, and the observation of synthetic populations, as discussed in \cite{fumanelli,mossong08,visi226}. Accurately representing these features is therefore essential for reliable prediction and for the design of effective non-pharmaceutical interventions (NPIs)~\cite{block,BARTHELEMY2005275}.

In this context, kinetic modelling provides a natural framework to describe epidemic dynamics as the result of interactions among a large number of heterogeneous agents. These approaches enable a systematic link between microscopic agent-based descriptions and macroscopic observables, offering deeper insight into the emergence of collective behaviours~\cite{Dimarco2020,DTZ,GF,martalo26}; see also~\cite{Albi2022,zanella_rev} for a review. A key feature of kinetic epidemic models is the possibility to incorporate contact heterogeneity at the mesoscopic scale. In particular, they couple standard compartmental dynamics with an evolution equation for the distribution of contacts, where the number of daily interactions is treated as an internal variable. This distribution evolves under a thermalization operator describing contact formation mechanisms, and its asymptotic profile may exhibit either thin or fat tails depending on the underlying social structure~\cite{FMZ}.  

The derivation of effective macroscopic descriptions consistent with the underlying microscopic interaction dynamics is therefore of paramount importance to enhance the explanatory power of the developed models  and to understand how infections propagate in large heterogeneous interacting systems. However, in the presence of uncertainties affecting the microscopic dynamics or the interaction structure, the characterization of a closed macroscopic limit is not always guaranteed, and the resulting limit may not be explicitly computable. Therefore, uncertainties in the contact-formation dynamics or in the initial distribution may hinder the derivation of effective macroscopic models capable of describing the observable dynamics.  

Furthermore, the kinetic formulation entails a substantial computational cost. For each realization of the random input, one must solve a system of kinetic equations with operators depending on both the uncertainties and the evolving epidemiological state. In moderately or highly dimensional stochastic settings, standard uncertainty quantification approaches based on repeated full-resolution simulations rapidly become intractable. This is a well-known bottleneck in multiscale and highly oscillatory PDEs with random parameters, where simultaneous resolution of physical and stochastic scales is computationally prohibitive~\cite{BLPZ2022,DLPZ2021,JLZC2026,LL25,LPZ2022}.

A widely used way around this bottleneck is to resort to multi-fidelity constructions~\cite{ZNX14}. The key observation, emphasized in recent work on kinetic epidemic models, is that cheaper surrogates are naturally available: (i) macroscopic closures obtained from the steady Fokker--Planck equilibria, and (ii) coarsened kinetic solvers that keep the right parametric dependence but on a cheaper mesh. On the other hand, the fully resolved kinetic solver is still needed to capture the fast contact relaxation and the coupled epidemic dynamics accurately. Multi-fidelity then proceeds in two stages: first, one explores the random space with the cheapest model and selects, via a greedy algorithm, a small set of “important’’ parameter samples; next, one runs the expensive solver only on this set and reconstructs high-fidelity solutions at generic parameters by reusing the projection rule learned from the low-fidelity ensemble~\cite{Liu2020}. This is the classical bi-fidelity pattern, and it already yields large savings whenever the low-fidelity model captures well the parametric geometry. 

The goal of this paper is to adapt this multi-fidelity paradigm to a kinetic epidemic model with uncertain contact structure and to demonstrate, on epidemic observables, that only a handful of high-fidelity simulations are enough to recover the statistics of the full kinetic model. Our contributions are threefold. First, building on the Fokker--Planck derivation for uncertain social contacts, we recall a simple kinetic system for epidemic dynamics and its macroscopic limits, emphasizing how the tail behavior of the equilibrium depends on the random parameter and why this is critical for NPIs. Second, we design bi-fidelity algorithms that use kinetic and macroscopic solvers in a hierarchical way, closely mirroring the multi-fidelity constructions developed for oscillatory PDEs. Third, through numerical tests with up to ten random dimensions, we show that the proposed methods reproduce the high-fidelity kinetic solution with errors of practical size while reducing the computational cost by magnitudes.

\paragraph{Organization.}
Section \ref{Analysis} recalls the kinetic epidemic model with uncertain contact distributions and the associated macroscopic closures. Section \ref{Algorithms} presents the bi- and tri-fidelity algorithms specialized to this model. Section \ref{Numerical} reports numerical results that quantify accuracy and efficiency in multiple scenarios.

\section{Kinetic compartmental models with uncertain contact structure}\label{Analysis}

We consider a prototypical system of agents subdivided in the following epidemiological relevant states: susceptible ($S$) agents are the ones that can contract the disease, infected and infectious ($I$) agents are responsible for the spread of the disease, exposed ($E$) agents have been in contact with infectious agents, and removed ($R$) agents cannot spread the disease. In the following, we will indicate with $\mathcal C = \{S,E,I,R\}$ the compartmentalization of the population. We highlight that the present approach may be extended to other type of compartmentalization of the system of agents. 

We denote with $f_J = f_J(x,t,\z)$, $J \in \mathcal C$, the distribution of the number of contacts $x\in\mathbb R^+$ at time $t\ge0$ of agents in the compartment $J$. The random vector $\z \in \mathbb R^{d_\z}$, $d_\z\in \mathbb N$, with known distribution $p(\z)$, collects all the uncertainties affecting the formation of a large time equilibrium density.  Hence, the total contact distribution of the society is obtained as 
\[
\sum_{J \in \mathcal C} f_J(x,t,\z) = f(x,t,\z), \qquad \int_{\mathbb R^+}f(x,t,\z) dx = 1,
\]
and the mass fractions of the population in each compartment are defined as 
\[
\rho_J(t,\z) = \int_{\mathbb R^+} f_J(x,t,\z)dx, 
\]
while their moment of order $r>0$ are given by
\[
\rho_J(t,\z)m_{r,J}(t,\z) = \int_{\mathbb R^+}x^r f_J(x,t,\z)dx.
\]
In the following, to simplify notations we will indicate with $m_J(t,\z) = m_{1,J}(t,\z) $ the first order moment, corresponding then to the cases $r=1$. 

Following the approach presented in \cite{DPTZ} we are interested in the evolution of the kinetic densities $(f_J)_{J \in \mathcal C}$ solution to 
\begin{equation}
\label{eq:kinetic}
    \begin{cases}
\partial_t f_S(x,t,\z) = -K(f_S,f_I)(x,t,\z) + \dfrac{1}{\tau}Q_S(f_S)(x,t,\z), \\
\partial_t f_E(x,t,\z) =  K(f_S,f_I)(x,t,\z) - \gamma_E f_E(x,t,\z) +  \dfrac{1}{\tau}Q_E(f_E)(x,t,\z) \\
\partial_t f_I(x,t,\z) = \gamma_E f_E(x,t,\z) -\gamma_I f_I(x,t,\z) + \dfrac{1}{\tau}Q_I(f_I)(x,t,\z) \\
\partial_t f_R(x,t,\z) = \gamma_I f_I(x,t,\z) +  \dfrac{1}{\tau}Q_R(f_R)(x,t,\z). 
    \end{cases}
\end{equation}

In \eqref{eq:kinetic} the transmission of the infection is governed by the local incidence rate
\begin{equation}
\label{eq:incidence}
K(f_S,f_I)(x,t,\z) = f_S(x,t,\z) \int_{\mathbb R^+} \kappa(x,x_*) f_I(x_*,t,\z)dx_*,
\end{equation}
being $\kappa(x,x_*)$ the contact function weighting the frequency of contacts between susceptible and infected agents. In the following, we will assume 
\begin{equation}
\label{eq:kappa}
\kappa(x,x_*) = \beta x x_*,\qquad \beta>0,
\end{equation}
where the parameter  $\beta$ scales the overall intensity of the contact–infection process (i.e., the baseline transmission rate). We remark that the contact function in \eqref{eq:kappa} can be generalized to take into account different impacts of the contact dynamics to take into account superlinear amplifications of contact effect, accentuating the influence of superspreaders, as further discussed in \cite{proc_MZ,zanella_rev}. 
Hence, plugging \eqref{eq:kappa} into \eqref{eq:incidence} we get
\[
K(f_S,f_I)(x,t,\z) = \beta x\, f_S(x,t,\z) \rho_I(t,\z) m_{I}(t,\z). 
\]
Within the choice in \eqref{eq:kappa}, the incidence rate is proportional to the product of the number of contact of susceptible and infected agents. In \eqref{eq:kinetic} we also introduced the transition rate between exposed to infected compartments, $\gamma_E>0$, and the recovery rate $\gamma_I>0$. 

The operators $Q_J(x,t,\z)$ determine the emergence of the contact dynamics of the agents' system. In \cite{FMZ,zanella_rev} it has been derived from an agent-based perspective the Fokker-Planck operator
\begin{equation}
\label{eq:QJ}
\begin{split}
Q_J(f_J)(x,t,\z) =& \dfrac{\mu}{2\theta(\z)} \partial_x \left[ x^{1-\alpha(\theta(\z))}\left(\left(\dfrac{x}{m_J(t)} \right)^{\theta(\z)} -1 \right)f_J(x,t,\z) \right] \\
&+ \dfrac{\sigma^2}{2} \partial_x^2 (x^{2-\alpha(\theta(\z))}f_J(x,t,\z))
\end{split}
\end{equation}
coupled with no-flux boundary condition at $x=0$. In \eqref{eq:QJ} we introduced the coefficients $\mu,\sigma^2>0$, $\theta = \theta(\z) \in [-1,1]$ is a random parameter and 
\[
\alpha(\theta) = \dfrac{1+\theta(\z)}{2} \in [0,1]. 
\]
This introduced operator is always mass preserving since
\[
\int_{\mathbb R^+} Q_J(f_J)(x,t,\z) dx = 0,
\]
and it is momentum preserving if $\theta \equiv \pm 1$ since from \eqref{eq:QJ} we get
\[
\int_{\mathbb R^+} x Q_J(f_J)(x,t,\z) dx = -\dfrac{\mu}{2\theta} \int_{\mathbb R^+} x^{1-\alpha} \left( \left( \dfrac{x}{m_J}\right)^{\theta}-1 \right)f_J(x,t,\z)dx. 
\]
More generally, for any $\theta \in (-1,1)$ the momentum is not a conserved quantity for the introduced collision operator.  In particular, this non-conservative feature hinders the derivation of closed macroscopic laws for the first moment. Indeed, its evolution is intrinsically coupled with higher-order statistical moments of the distribution, preventing a self-contained macroscopic description.

\subsection{Equilibrium distribution of the collision operator}

The Fokker-Planck-type operator defined in \eqref{eq:QJ} is such that its equilibrium density $f^q_J(x,t,\z)$ parametrised by $m_J(t,\z)$ is obtained as the unique solution to the following differential equation 
\[
\dfrac{\mu}{\theta(\z)} x^{1-\alpha(\theta(\z))} \left(\left( \dfrac{x}{m_J(t,\z)} \right)^{\theta(\z)}-1 \right)f_J^q(x,t,\z) + \sigma^2 \partial_x (x^{2-\alpha(\theta(\z))}f_J^q(x,t,\z)) = 0, 
\]
which is given by
\begin{equation}
\label{eq:fqJ}
f_J^q(x,t,\z) = C^\theta_{\mu,\sigma^2} x^{\frac{\mu}{\sigma^2\theta(\z)}-2+\alpha(\theta(\z))}\exp\left\{-\dfrac{\mu}{\sigma^2\theta^2(\z)}\left( \dfrac{x}{m_J(t,\z)} \right)^{\theta(\z)} \right\},
\end{equation}
where $C^\theta_{\mu,\sigma^2}>0$ is a normalization constant.
The equilibrium solution in \eqref{eq:fqJ} inherits a direct dependence on the uncertain parameter $\theta(\z)$, which governs its asymptotic structure. In particular, different realizations of $\theta$ modify the tail decay of the equilibrium distribution, leading to a variability in its behaviour

Indeed, we may observe how, if the distribution of the uncertain parameter is such that $\theta(\z)\equiv \z$ 
and $p(\z) = \delta(\z-1)$ we get a Gamma distribution
\[
\mathbb E[f^q_J(x,t,\z)\big | \z \equiv 1] = \dfrac{\lambda^\lambda}{m_J(t)^\lambda \Gamma(\lambda)}x^{\lambda-1}\exp\left\{-\dfrac{\lambda x}{m_J(t)}\right\}, \qquad \lambda = \dfrac{\mu}{\sigma^2},
\]
having slim tails and such that 
\[
\int_{\mathbb R^+} x^2 f^q_J(x,t,\z) dx = \frac{\lambda+1}{\lambda}\, m_J(t)^{2}, \qquad \lambda=\frac{\mu}{\sigma^{2}}.
\] 
On the other hand, if $p(\z) = \delta(\z+1)$, we obtain the inverse Gamma distribution 
\[
\mathbb E[f^q_J(x,t,\z)\big | \z \equiv -1] = \dfrac{(\lambda m_J(t))^{\lambda+1}}{\Gamma(\lambda+1)}x^{-2-\lambda}\exp\left\{ -\dfrac{\lambda m_J(t)}{x}\right\}, \qquad \lambda = \dfrac{\mu}{\sigma^2}, 
\]
which is a fat-tailed distribution. More generally, the equilibrium distribution is influenced by the uncertain parameter affecting $\theta = \theta(\z)$. In particular, for any $\theta\ge0$ the equilibrium density exhibits an exponential decay for $x \gg0$.  On the other hand, for $\theta<0$ the equilibrium density possesses polynomial decay for $x \gg0$. Uncertainties in the parameter $\theta$ have a strong impact in terms of the behaviour of the multi-agent system. 

\subsection{Macroscopic equations}
From the kinetic model \eqref{eq:kinetic} we may recover classical compartmental models by computing the moments of the kinetic densities $(f_J)_{J \in \mathcal C}$, see \cite{DPTZ}. Since the operators $(Q_J)_{J \in \mathcal C}$ are mass preserving, if we integrate \eqref{eq:kinetic} in $x \in \mathbb R_+$ we get
\begin{equation}
\label{eq:macro1}
\begin{cases}
\dfrac{d}{dt}\rho_S(t,\z) = -\beta m_S(t,\z)\rho_S(t,\z) m_I(t,\z)\rho_I(t,\z), \\
\dfrac{d}{dt} \rho_E(t,\z)= \beta m_S(t,\z)\rho_S(t,\z) m_I(t,\z)\rho_I(t,\z) - \gamma_E\rho_E(t,\z), \\
\dfrac{d}{dt}\rho_I(t,\z) = \gamma_E \rho_E(t,\z) - \gamma_I \rho_I(t,\z), \\
\dfrac{d}{dt}\rho_R(t,\z) = \gamma_I \rho_I(t,\z).
\end{cases}
\end{equation}
We may observe that the system for mass fractions is not closed since $\rho_J(t,\z)$ depends on the first order moment of the kinetic density $m_J(t,\z)$ whose evolution can be computed from \eqref{eq:kinetic}
 and is given by
 \begin{equation}
\label{eq:macro2}
\begin{cases}
\dfrac{d}{dt} (\rho_S m_S)(t,\z) =& -\beta m_{2,S}(t,\z)\rho_S(t,\z) m_I(t,\z)\rho_I(t,\z) + \dfrac{1}{\tau} \int_{\mathbb R^+} x Q_S(f_S)(x,t,\z)dx \\
\dfrac{d}{dt} (\rho_E m_E)(t,\z) =& \beta m_{2,S}(t,\z)\rho_S(t,\z) m_I(t,\z)\rho_I(t,\z) - \gamma_E m_E(t,\z) \rho_E(t,\z) \\
&  + \dfrac{1}{\tau} \int_{\mathbb R^+} x Q_E(f_E)(x,t,\z)dx \\
\dfrac{d}{dt} (\rho_I m_I)(t,\z) =&\gamma_E m_E(t,\z)\rho_E(t,\z) - \gamma_I m_I(t,\z)\rho_I(t,\z) + \dfrac{1}{\tau} \int_{\mathbb R^+} x Q_I(f_I)(x,t,\z)dx  \\
\dfrac{d}{dt} (\rho_R m_R)(t,\z) =& \gamma_I m_I(t,\z)\rho_I(t,\z) + \dfrac{1}{\tau} \int_{\mathbb R^+} x Q_R(f_R)(x,t,\z)dx. 
\end{cases}
\end{equation}
Hence, for any $\z$ such that $\theta \in \{-1,1\}$, e.g. $\z \sim \textrm{Bernoulli}(q)$ we obtain that 
\[
\int_{\mathbb R^+}x Q_J(f_J)(x,t,\z)dx = 0, \qquad J \in \mathcal C, 
\]
and in the limit $\tau \ll 1$, 
\[
m_{2,J}(t,\z) = \int_{\mathbb R^+} x^2 f_{J,\rho_J,m_J}^q(x,\z)dx = \Lambda_\theta(\z) m_J^2(t,\z),\qquad \Lambda_\theta(\z) = \left( \dfrac{\lambda+\theta(\z)}{\lambda}\right)^{\theta(\z)}. 
\]
Therefore, in the limit $\tau\ll1$ we can approximate \eqref{eq:macro2} as follows
\begin{equation}
\label{eq:macro2_closed}
\begin{cases}
\dfrac{d}{dt} (\rho_S m_S)(t,\z) =& -\beta \Lambda_\theta(\z)\rho_S(t,\z)m_S^2(t,\z) m_I(t,\z)\rho_I(t,\z)  \\[4pt]
\dfrac{d}{dt} (\rho_E m_E)(t,\z) =& \beta \Lambda_\theta(\z)\rho_S(t,\z)m_S^2(t,\z) m_I(t,\z)\rho_I(t,\z) - \gamma_E m_E(t,\z) \rho_E(t,\z) \\[4pt]
\dfrac{d}{dt} (\rho_I m_I)(t,\z) =&\gamma_E m_E(t,\z)\rho_E(t,\z) - \gamma_I m_I(t,\z)\rho_I(t,\z) \\[4pt]
\dfrac{d}{dt} (\rho_R m_R)(t,\z) =& \gamma_I m_I(t,\z)\rho_I(t,\z).  
\end{cases}
\end{equation}

We conclude that the passage from the kinetic description \eqref{eq:kinetic} to macroscopic compartmental dynamics naturally generates a hierarchy of moment equations, where lower-order quantities are coupled with higher-order moments and with residual kinetic contributions through the operators \(Q_J\). In particular, closure at the level of mass fractions alone is not available in general, as the dynamics of $m_{2,J}(t,\z)$ depends on the specific closure distribution that shaped by uncertainties and their distribution. 

In special cases, fast relaxation towards local equilibrium allows one to express higher-order moments in terms of the first-order ones, leading to an approximate closure of the system. In this limit, uncertainty in the microscopic interactions is transferred to the macroscopic level through the parameter-dependent factor \(\Lambda_\theta(\mathbf z)\), which modulates the effective nonlinear incidence terms.

This shows that, while a closed macroscopic description can be recovered under suitable asymptotic regimes, the resulting effective dynamics retains a nontrivial dependence on the underlying kinetic uncertainty, which persists in the form of modified interaction coefficients.
\begin{remark}
If $p(\z) = \delta(\z-1)$, $\theta\equiv \z$, the resulting macroscopic equations are the ones derived from a Gamma-type equilibrium distribution for which we get
\[
\mathbb E[m_{2,J}(t,\z) | \theta\equiv 1] =  \dfrac{\lambda+1}{\lambda} m_J^2(t), \qquad \lambda = \dfrac{\mu}{\sigma^2} 
\]
which gives the macroscopic dynamics 
 \begin{equation}
\label{eq:macro3}
\begin{cases}
\dfrac{d}{dt} (\rho_S m_S)(t) =& -\beta \dfrac{\lambda+1}{\lambda}m_S^2\rho_S m_I\rho_I  \\[6pt]
\dfrac{d}{dt} (\rho_E m_E)(t) =& \beta \dfrac{\lambda+1}{\lambda}m_S^2\rho_S m_I\rho_I- \gamma_E m_E \rho_E \\[6pt]
\dfrac{d}{dt} (\rho_I m_I)(t) =&\gamma_E m_E\rho_E - \gamma_I m_I\rho_I   \\[6pt]
\dfrac{d}{dt} (\rho_R m_R)(t) =& \gamma_I m_I\rho_I 
\end{cases}
\end{equation}
While, if $p(\z) = \delta(\z+1)$, $\theta\equiv \z$, the equilibrium distribution is of inverse Gamma-type for which we have
\[
\mathbb E[m_{2,J}(t,\z) | \theta\equiv -1] =  \dfrac{\lambda}{\lambda-1} m_J^2(t), \qquad \qquad \lambda = \dfrac{\mu}{\sigma^2} 
\]
which gives the macroscopic dynamics 
 \begin{equation}
\label{eq:macro4}
\begin{cases}
\dfrac{d}{dt} (\rho_S m_S)(t) =& -\beta \dfrac{\lambda}{\lambda-1}m_S^2\rho_S m_I\rho_I  \\[6pt]
\dfrac{d}{dt} (\rho_E m_E)(t) =& \beta \dfrac{\lambda}{\lambda-1}m_S^2\rho_S m_I\rho_I- \gamma_E m_E \rho_E \\[6pt]
\dfrac{d}{dt} (\rho_I m_I)(t) =&\gamma_E m_E\rho_E - \gamma_I m_I\rho_I   \\[6pt]
\dfrac{d}{dt} (\rho_R m_R)(t) =& \gamma_I m_I\rho_I 
\end{cases}
\end{equation}
The obtained macroscopic systems are not equivalent  as they strongly depend on
the uncertainties in contact formation dynamics of the population of interest, see \cite{zanella_rev,zanella_medaglia}.
\end{remark}

\begin{remark}
For a general $p(\z)$ we cannot close the system of macroscopic equations since the first-order moment is not a conserved quantity.  
\end{remark}
\section{Multi-fidelity algorithms}\label{Algorithms}

In the previous section, we introduced the kinetic epidemic model and its macroscopic approximations. These models naturally have different levels of accuracy and computational cost. The microscopic kinetic model \eqref{eq:kinetic} describes the full contact distribution $f_J(x,t,\z)$ and is the most accurate model in this work. However, it is also expensive, since it requires solving the kinetic equations in the contact variable $x$ for each random sample $\z$. On the other hand, the macroscopic models \eqref{eq:macro1}-\eqref{eq:macro2} only evolve averaged quantities such as $\rho_J(t,\z)$ and $m_J(t,\z)$. They are cheaper, but less accurate.

This motivates the use of multi-fidelity methods. Instead of solving the microscopic model for all random samples, we use the cheaper models to explore the random space and select important parameter points. The expensive kinetic solver is then evaluated only at these selected points. In this section, we explain how the bi-fidelity and tri-fidelity methods are applied to the epidemic models introduced above. The numerical solvers adopted for the different models will be specified in the next Section.

For each random parameter $\z \in I_\z$, we denote by $U^H(\z)$ the high-fidelity solution, which in this work is always obtained using the microscopic kinetic solver. The quantity $U^H(\z)$ may represent either the full time histories of the compartmental densities, such as $\rho_S(t,\z)$, $\rho_E(t,\z)$, and the other compartment variables, or their values at the final time.
For the bi-fidelity method, we denote by $U^L(\z)$ the corresponding quantities obtained from the low-fidelity solver, which in our numerical tests is given by the macroscopic model. For the tri-fidelity method, we additionally introduce a medium-fidelity solution $U^M(\z)$.
The specific definitions of $U^L$, $U^M$, and $U^H$ will be detailed in the numerical experiments.

\begin{algorithm}
    \caption{Bi-fidelity construction for kinetic epidemic observables}
    \label{Bi-fid}
    \begin{algorithmic}[1]
        \STATE Choose a candidate set $\Gamma=\{z_1,z_2,\ldots,z_M\}\subset I_\z$ and compute the low-fidelity snapshots $\mathbf{U}^L(z_j)$ for all $z_j\in\Gamma$.
        \STATE Apply the greedy procedure of~\cite{DeVore}, which is in appendix~\ref{app:chol-selection}, to the set $\{\mathbf{U}^L(z_j)\}_{j=1}^M$ and select $r\ll M$ representative parameters $\gamma=\{z_{i_1},z_{i_2},\ldots,z_{i_r}\}\subset\Gamma$.
        \STATE Evaluate the high-fidelity solver at the selected samples and define the reduced snapshot spaces
        \[
        \mathscr{U}^L(\gamma)=\mathrm{span}\{\mathbf{U}^L(z_{i_1}),\ldots,\mathbf{U}^L(z_{i_r})\},
        \qquad
        \mathscr{U}^H(\gamma)=\mathrm{span}\{\mathbf{U}^H(z_{i_1}),\ldots,\mathbf{U}^H(z_{i_r})\}.
        \]
        \STATE For a new parameter value $\z\in I_\z$, compute $\mathbf{U}^L(z)$ and determine the coefficients $\mathbf{c}^L(z)=(c_1^L(z),\ldots,c_r^L(z))^T$ from the Galerkin system
        \begin{equation}\label{C-N}
            \mathbf{G}^L \mathbf{c}^L(\z)=\mathbf{f}^L(\z),
            \qquad
            G^L_{mn}=\langle \mathbf{U}^L(z_{i_m}),\mathbf{U}^L(z_{i_n})\rangle^L,
            \qquad
            f^L_m(\z)=\langle \mathbf{U}^L(\z),\mathbf{U}^L(z_{i_m})\rangle^L,
        \end{equation}
        where $\langle\cdot,\cdot\rangle^L$ denotes the low-fidelity inner product induced by the chosen observable discretization.
        \STATE Define the bi-fidelity approximation by transferring the same coefficients to the high-fidelity snapshot space:
        \begin{equation}\label{UB}
            \mathbf{U}^B(\z)=\sum_{m=1}^r c_m^L(\z)\,\mathbf{U}^H(z_{i_m}).
        \end{equation}
    \end{algorithmic}
\end{algorithm}

We next introduce the tri-fidelity approximation. In the numerical tests of this paper, $U^H$ is always the microscopic kinetic solver. The medium-fidelity model is a macroscopic solver that is more accurate than the cheapest model. The low-fidelity model is a simplified macroscopic solver used mainly for fast sample selection.
The main distinction between the bi-fidelity and tri-fidelity constructions lies in the computation of the projection coefficients. In the bi-fidelity framework, these coefficients are determined using the low-fidelity model. In contrast, the tri-fidelity method uses the low-fidelity model only to identify representative samples, whereas the medium-fidelity model is employed to compute the projection coefficients. The final approximation is then reconstructed using the high-fidelity solutions at the selected samples.

\begin{algorithm}
\caption{Tri-fidelity approximation for kinetic epidemic observables}
\begin{algorithmic}[1]
\STATE Select a candidate set
\[
    \Gamma_N=\{z_1,z_2,\ldots,z_N\}\subset I_\z .
\]
\STATE Run the low-fidelity solver for each $\z_j\in\Gamma_N$ and obtain
\[
    \mathbf{U}^L(\Gamma_N)=\{\mathbf{U}^L(z_1),\mathbf{U}^L(z_2),\ldots,\mathbf{U}^L(z_N)\}.
\]
\STATE Apply the Pivoted Cholesky sample selection procedure (Algorithm~\ref{alg:chol-greedy-epidemic}) to $\mathbf{U}^L(\Gamma_N)$ and select $K$ important samples
\[
    \gamma_K=\{z_{i_1},z_{i_2},\ldots,z_{i_K}\}\subset\Gamma_N,
    \qquad K\ll N .
\]
\STATE Run the medium-fidelity solver and the high-fidelity microscopic solver at each selected sample $z_{i_k}\in\gamma_K$.
\STATE Construct the medium-fidelity approximation space
\[
    \mathscr U^M(\gamma_K)
    =
    \operatorname{span}
    \{\mathbf{U}^M(z_{i_1}),\mathbf{U}^M(z_{i_2}),\ldots,\mathbf{U}^M(z_{i_K})\}.
\]
\STATE Construct the high-fidelity approximation space
\[
    \mathscr U^H(\gamma_K)
    =
    \operatorname{span}
    \{\mathbf{U}^H(z_{i_1}),\mathbf{U}^H(z_{i_2}),\ldots,\mathbf{U}^H(z_{i_K})\}.
\]

\STATE For a new parameter value $\z\in I_\z$, compute the medium-fidelity solution $U^M(\z)$.
\STATE Find the projection coefficients
\[
    \mathbf c^M(\z)
    =
    \big(c^M_1(\z),c^M_2(\z),\ldots,c^M_K(\z)\big)^T
\]
by solving
\begin{equation}\label{eq:TF_Galerkin}
    G^M\mathbf c^M(\z)=\mathbf f^M(\z),
\end{equation}
where
\[
    G^M_{mn}
    =
    \left\langle
    \mathbf{U}^M(z_{i_m}),\mathbf{U}^M(z_{i_n})
    \right\rangle^M,
    \qquad
    f^M_m(\z)
    =
    \left\langle
    \mathbf{U}^M(\z),\mathbf{U}^M(z_{i_m})
    \right\rangle^M .
\]
\STATE Construct the tri-fidelity approximation by
\begin{equation}\label{eq:TF_approx}
    \mathbf{U}^T(\z)
    =
    \sum_{m=1}^{K} c^M_m(\z)\mathbf{U}^H(z_{i_m}) .
\end{equation}
\end{algorithmic}
\end{algorithm}

Here $\langle\cdot,\cdot\rangle^M$ denotes the inner product for the medium-fidelity observables. The tri-fidelity method separates the tasks of point selection and coefficient computation. The cheapest model is used on the large candidate set $\Gamma_N$, while the more accurate medium-fidelity model is used to compute the projection rule. The high-fidelity kinetic model is still only used at the selected samples.

The above construction is particularly well suited to the present epidemic setting. In the present epidemic setting, the three fidelity levels are directly built from the model hierarchy introduced in Section~\ref{Analysis}. The high-fidelity model is the full kinetic system \eqref{eq:kinetic} with the Fokker--Planck contact operator \eqref{eq:QJ}; hence it evolves the complete contact distributions $f_J(x,t,\z)$ for all compartments $J\in\{S,E,I,R\}$. For the bi-fidelity tests, the low-fidelity model is the macroscopic closure obtained from the moment equations, namely the equation \eqref{eq:macro1} together with the system \eqref{eq:macro2_closed} when $\theta=\pm1$. For the tri-fidelity test, the medium-fidelity model is taken as the macroscopic moment model \eqref{eq:macro1}--\eqref{eq:macro2}, which keeps the evolution of both $\rho_J$ and $m_J$, while the low-fidelity model is a simplified version of this macroscopic system in which the integral terms of $x Q_J(f_J)(x,t,\z)$ are frozen at their initial values. Thus, the cheapest model is used only to select representative samples, the medium-fidelity model provides more reliable projection coefficients, and the high-fidelity kinetic model supplies the accurate snapshots used in the final reconstruction. This structure is precisely what makes the procedure effective in the numerical experiments reported in Section~\ref{Numerical}.

\section{Numerical Results}\label{Numerical}

In this section, we test the accuracy and efficiency of the proposed multi-fidelity methods for the kinetic epidemic model with uncertain contact dynamics. The main purpose is to check whether a small number of high-fidelity microscopic simulations can reproduce the main statistical quantities of the fully resolved kinetic solver. In all tests, the high-fidelity model is the microscopic kinetic system \eqref{eq:kinetic}, where the contact relaxation operator $Q_J(f_J)(x,t,\z)$ is defined by \eqref{eq:QJ}. The lower-fidelity models are chosen from the macroscopic systems derived in Section~\ref{Analysis}, or from simplified versions of these systems when the tri-fidelity method is used.

To measure the approximation error, we choose a fixed testing set $\{\hat{z}_i\}_{i=1}^{N_z}\subset I_\z$, which is independent of the candidate set used in the greedy selection step. At the final time $t=T$, we compute the average error between the high-fidelity solution and the multi-fidelity approximation as
\begin{equation} \label{error}
    \mathcal{E} = \frac{1}{N_z} \sqrt{\sum_{i=1}^{N_z} \left| \mathbf{U}^H(T,\hat{z}_i) - \mathbf{U}^{MF}(T,\hat{z}_i) \right|^2},
\end{equation}
where $\mathbf{U}^H$ denotes the quantity of interest obtained from the high-fidelity microscopic solver, and $\mathbf{U}^{MF}$ denotes the corresponding bi-fidelity or tri-fidelity approximation. In the following numerical examples, the quantities of interest are mainly the compartmental densities, since they directly describe the evolution of the susceptible and exposed populations.

For time integration, we use the fourth-order Runge--Kutta method. The role of the high-fidelity solver is to provide accurate reference solutions of the kinetic model. However, because this solver requires the evolution of the full contact distribution $f_J(x,t,\z)$ for each compartment and each random sample, it is computationally expensive. The low-fidelity solver, on the other hand, evolves only macroscopic quantities such as $\rho_J$ and $m_J$. It is therefore much cheaper, but it may not fully capture the kinetic relaxation effect. The multi-fidelity construction aims to combine these two features: the low-fidelity model is used to explore the random space and select representative samples, while the high-fidelity model is evaluated only at a small number of selected points.

The numerical results are organized as follows. In Section \ref{subsect:test1} we consider the special cases $\theta=\pm1$, for which the momentum-preserving property allows us to use closed macroscopic equations as a reliable low-fidelity model. In Section \ref{subsect:test2} we allow $\theta$ to vary in $[-1,1]$ and also introduces uncertainty in the epidemiological parameters. This case is more challenging because the macroscopic closure becomes less direct when $\theta$ is not restricted to $\pm1$. Finally, in Section \ref{subsect:test3} we study a tri-fidelity strategy, where a very cheap simplified model is used for point selection, while a more accurate macroscopic model is used to compute the projection coefficients.

\subsection{Test 1: Bi-fidelity when $\theta = \pm 1$} \label{subsect:test1}
The initial condition of first order moment is specified as

\begin{equation}
\label{eq:mJ_test1}
    m_J(t=0,\mathbf{z}) = 
    \begin{cases}\displaystyle
    3\left(1+\sum_{i=1}^{d} \dfrac{z_i\sin z_i}{i}\right),     & \text{if }J=S  \\
    \displaystyle 3\left(1+\sum_{i=1}^{d} \dfrac{z_i\sin z_i}{i}\right),     & \text{if }J=E  \\
    \displaystyle 3\left(1+\sum_{i=1}^{d} \dfrac{z_i\sin z_i}{2i}\right),     & \text{if }J=I  \\
    3,     & \text{if }J=R ,  \\
    \end{cases}
\end{equation}
where $d$ is the dimensionality of $\z$, and $\mathbf{z}=(z_1,\cdots,z_{d})$. This choice introduces uncertainty through the initial contact moments. In particular, we insert the same perturbation into the susceptible (S), exposed (E), and infected (I) classes. The removed class (R) is initialized with a constant first moment, since it does not directly contribute to the transmission process at the beginning of the simulation. The initial condition for the mass fractions is taken as 
$$
\rho_S(t=0,\z)=0.97,\qquad \rho_E(t=0,\z)=\rho_I(t=0,\z)=\rho_R(t=0,\z)=0.01, 
$$
while the initial mean values are defined as in \eqref{eq:mJ_test1}. 
In this test, the parameter $\theta$ is fixed as $\theta=\pm1$. These two values correspond to the cases where the collision operator preserves the first-order moment. This property is important because it allows us to use the closed macroscopic equations \eqref{eq:macro1} and \eqref{eq:macro2_closed} as the low-fidelity model. Hence, the macroscopic solver can directly use the initial values of $\rho_J$ and $m_J$.

The microscopic solver requires the initial kinetic density $f_J(x,0,\z)$ rather than only its moments. Since the first-order moment is defined by
\[
\rho_J(0,\z)m_J(0,\z)=\int_{\mathbb{R}^+}x f_J(x,0,\z)\,\ud x,
\]
we prescribe the initial contact density by the Gaussian-type distribution

\begin{equation*}
    f_J(x,0,\mathbf{z}) = C_{J}\frac{10}{\sqrt{2\pi}} \exp\Big\{-50\big(x-m_J(0,\z)\big)^2\Big\}, 
\end{equation*}
where $C_J$ is a normalization constant. This distribution is simple, positive, normalized on $\mathbb{R}^+$, and approximately has mean $m_J(0,\z)$. Therefore, it is consistent with the prescribed initial first-order moment and can be used as the initial condition for the microscopic kinetic equation.

We set $d=10$ and sample each component independently as $z_i\sim\mathcal U[-1,1]$. The candidate set contains $N_z=500$ random samples. For the high-fidelity solver, the time interval $[0,20]$ is divided into $200$ uniform steps, and the contact domain $[0,500]$ is discretized by $5000$ uniform cells. The model parameters are
\[
\beta=0.025,\qquad
\gamma_E=0.33,\qquad
\gamma_I=0.1,\qquad
\mu=0.5,\qquad
\sigma^2=0.1.
\]
We consider two relaxation regimes, $\tau=10^{-2}$ and $\tau=10^{-4}$. A smaller value of $\tau$ means faster relaxation of the contact distribution toward its local equilibrium. Therefore, comparing these two values allows us to see the stability of this method when the kinetic relaxation becomes stronger.

\begin{figure*} 
    \centering
    \subfloat[Average $L^2$ error of bi-fidelity approximations for $\rho_S$]{\includegraphics[width=.45\textwidth]{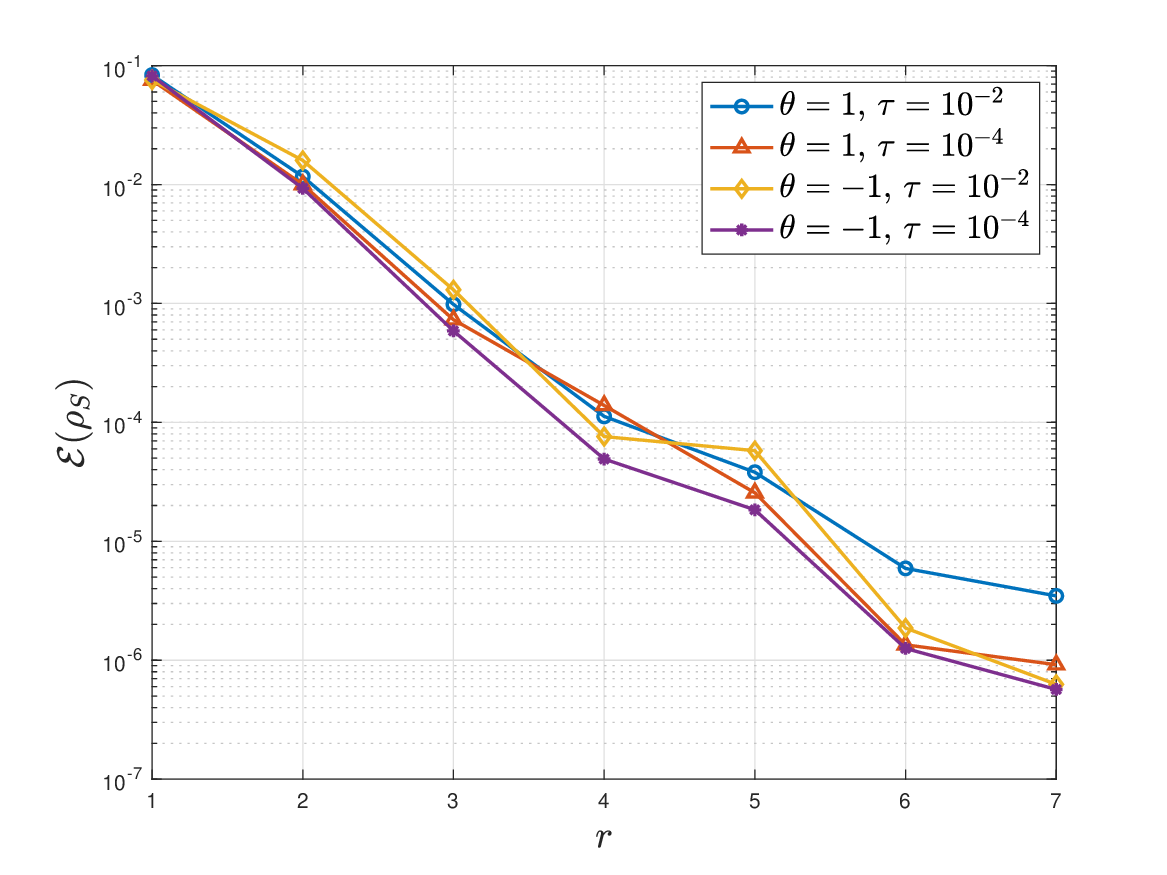}}
    \hspace{5mm}
    \subfloat[Average $L^2$ error of bi-fidelity approximations for $\rho_E$]{\includegraphics[width=.45\textwidth]{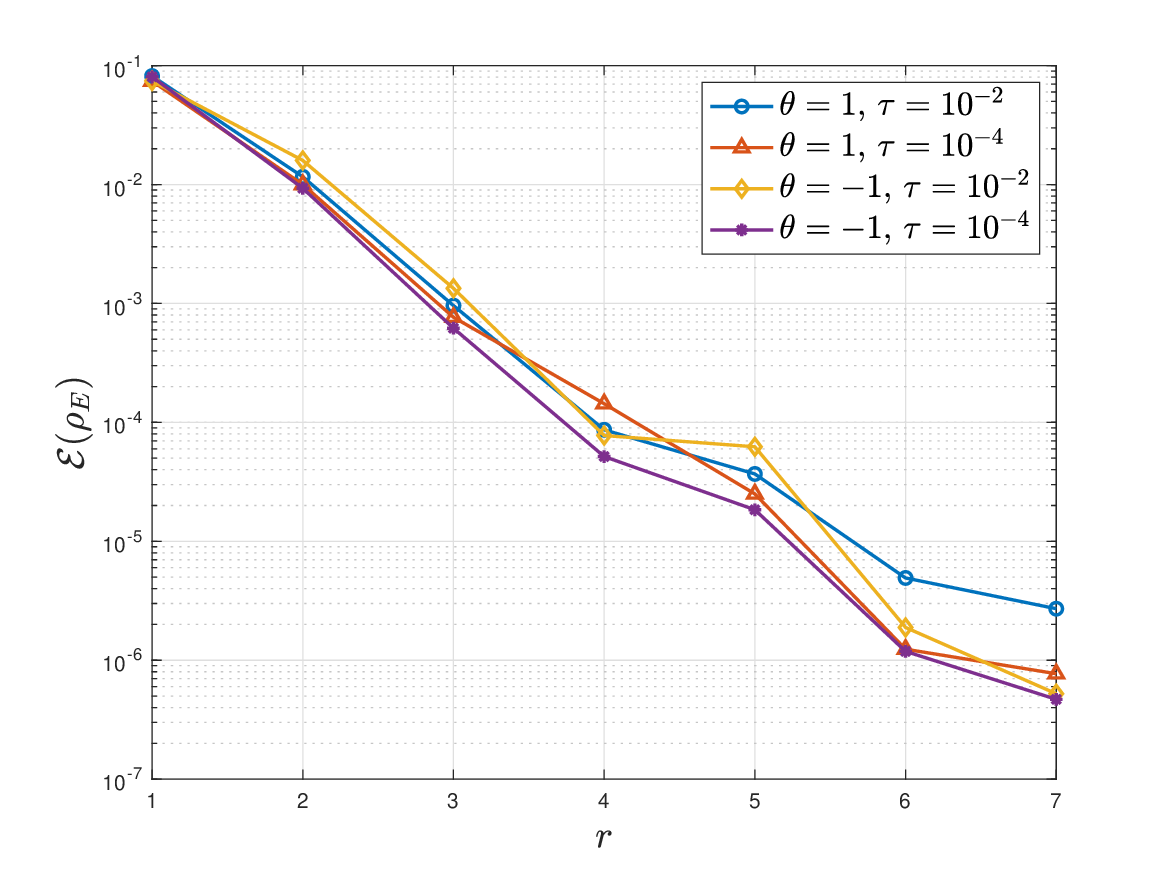}}

    \subfloat[Average $L^2$ error of bi-fidelity approximations for $\rho_I$]{\includegraphics[width=.45\textwidth]{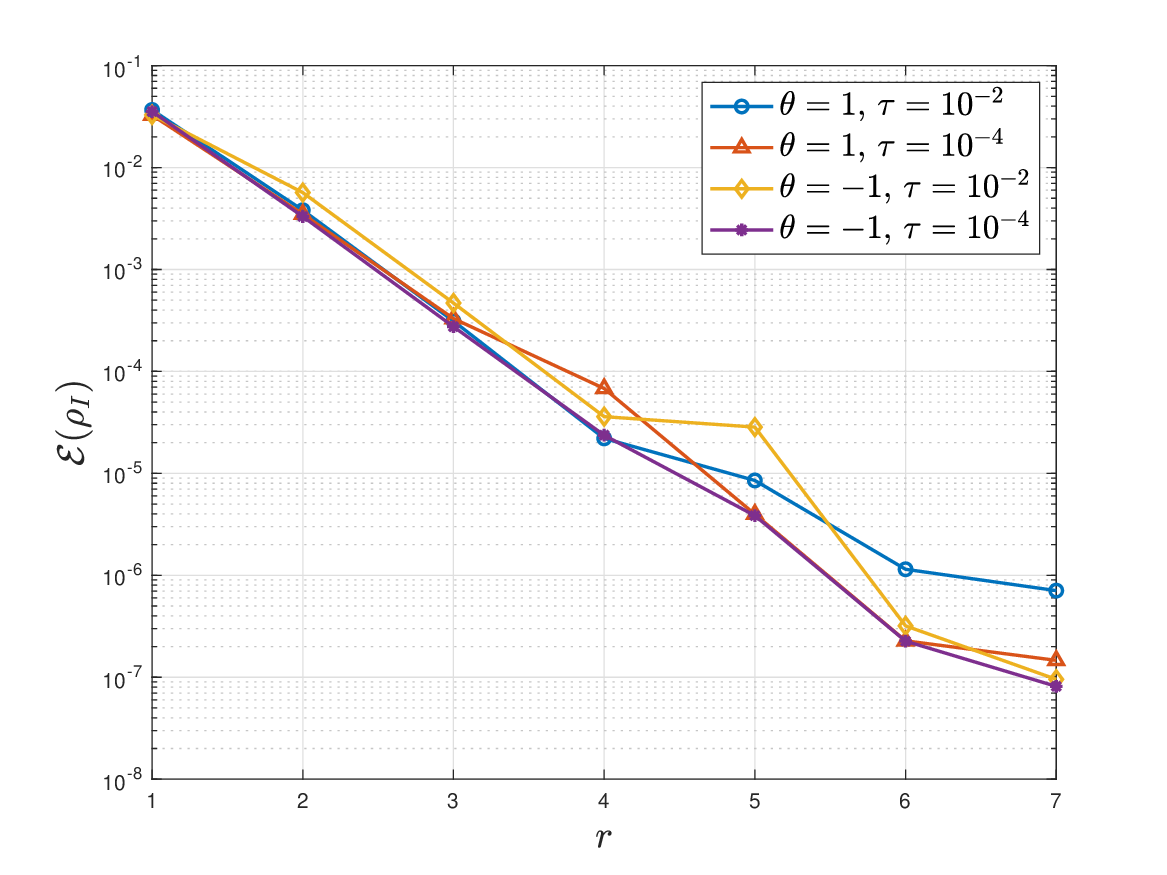}}
    \hspace{5mm}
    \subfloat[Average $L^2$ error of bi-fidelity approximations for $\rho_R$]{\includegraphics[width=.45\textwidth]{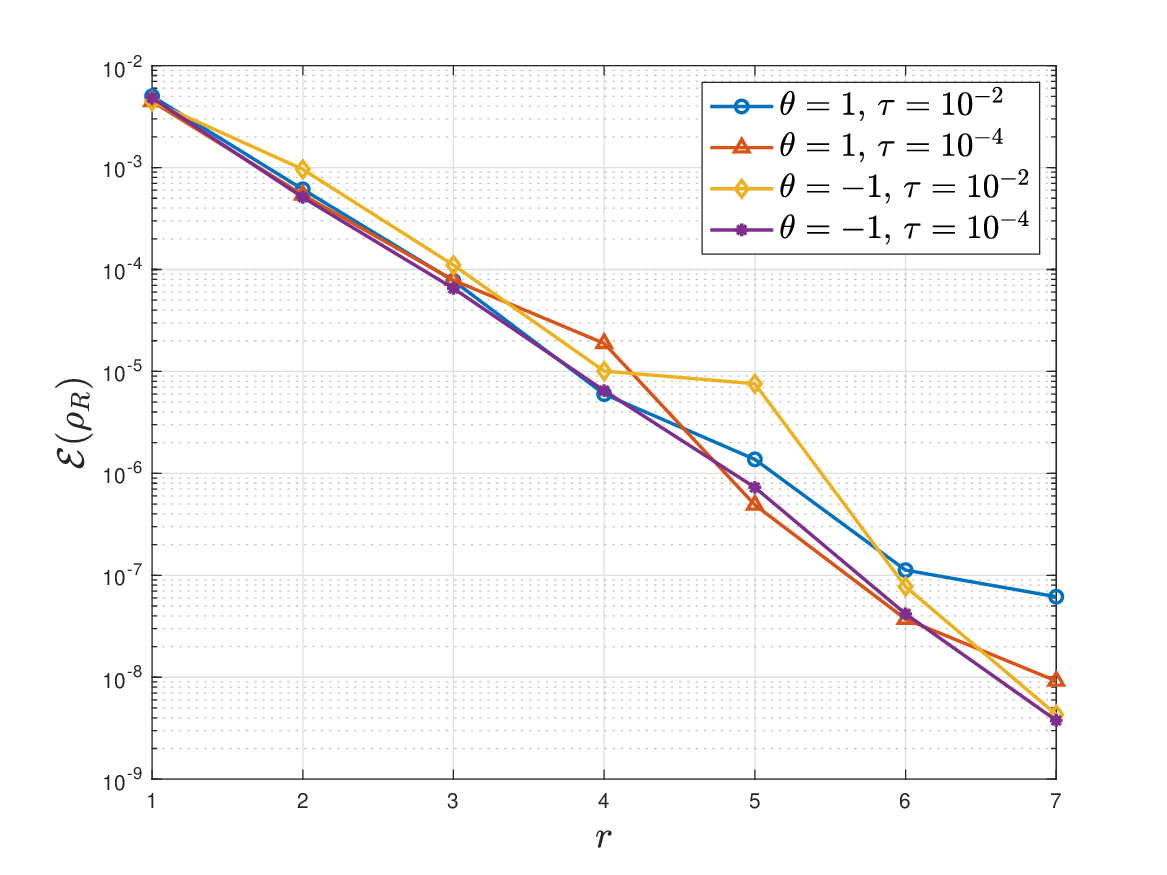}}
    \caption{Test 1: Average $L^2$ error of bi-fidelity approximations for $\rho_J$ with respect to the number of high-fidelity simulation runs at $\theta=\pm1$ and different $\tau$.}
    \label{Fig.Test1.1.m}
\end{figure*}

\begin{figure*}
    \centering
    \subfloat[Average $L^2$ error of bi-fidelity approximations for $m_S$]{\includegraphics[width=.45\textwidth]{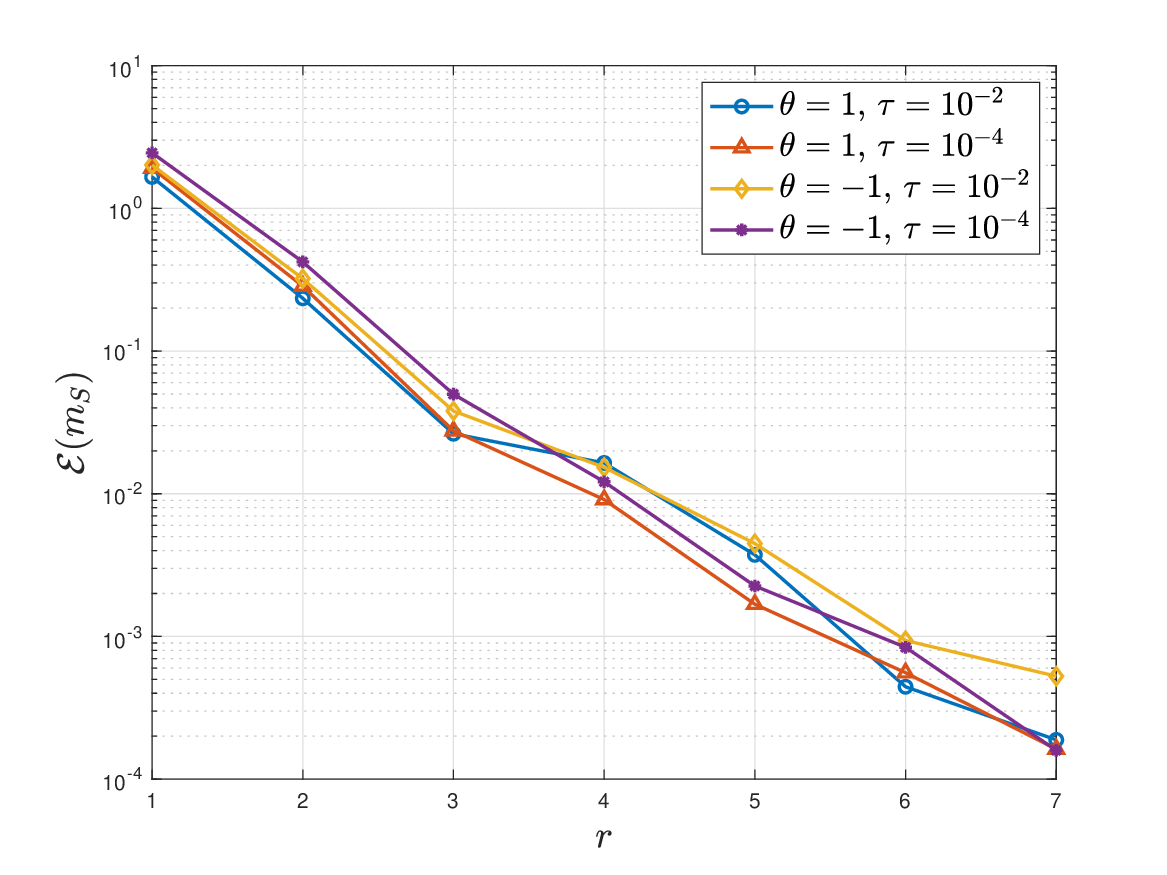}}
    \hspace{5mm}
    \subfloat[Average $L^2$ error of bi-fidelity approximations for $m_E$]{\includegraphics[width=.45\textwidth]{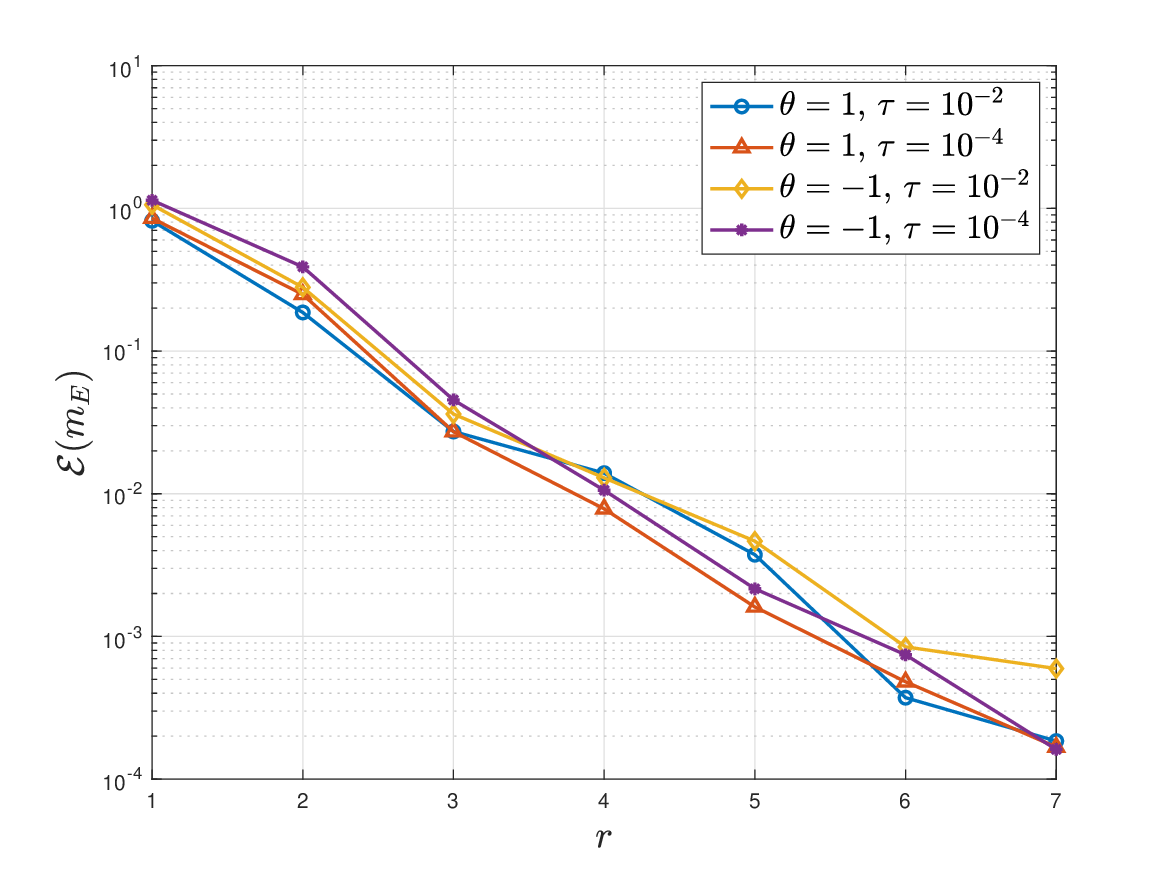}}

    \subfloat[Average $L^2$ error of bi-fidelity approximations for $m_I$]{\includegraphics[width=.45\textwidth]{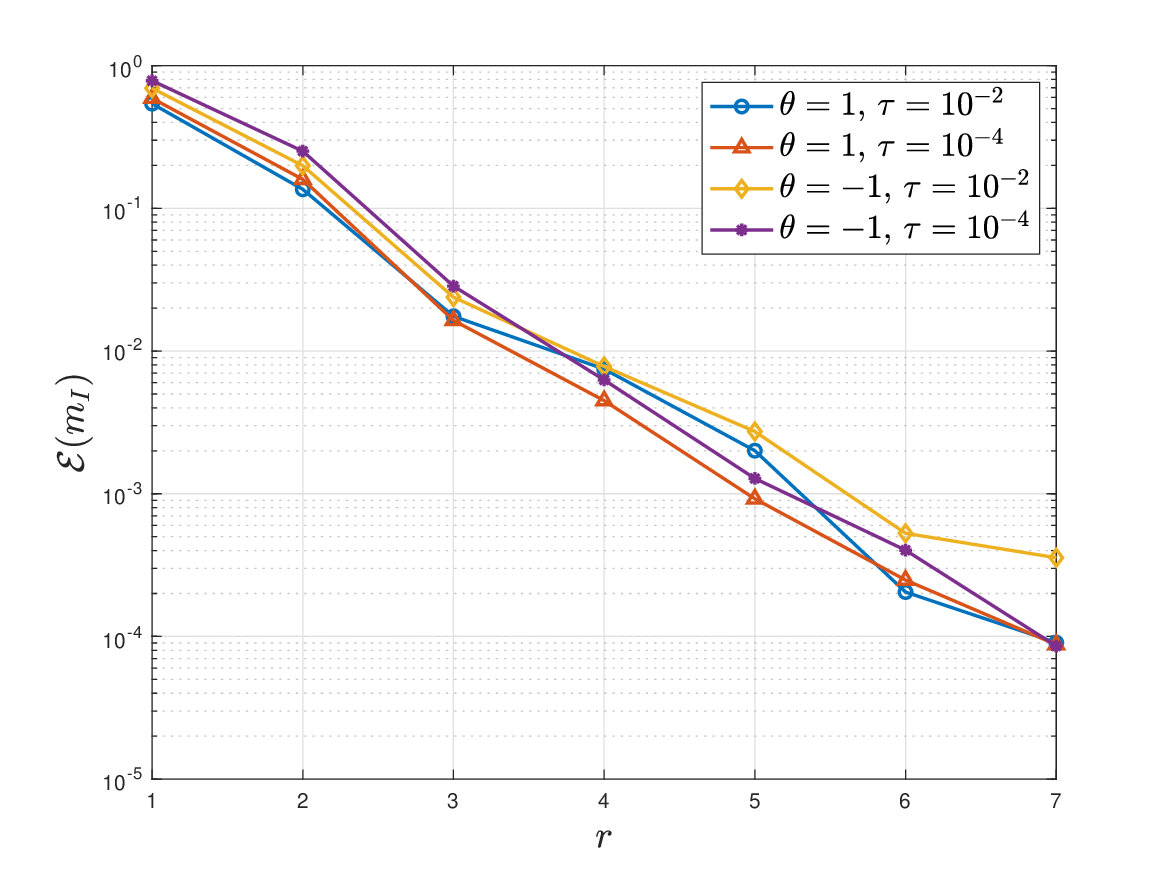}}
    \hspace{5mm}
    \subfloat[Average $L^2$ error of bi-fidelity approximations for $m_R$]{\includegraphics[width=.45\textwidth]{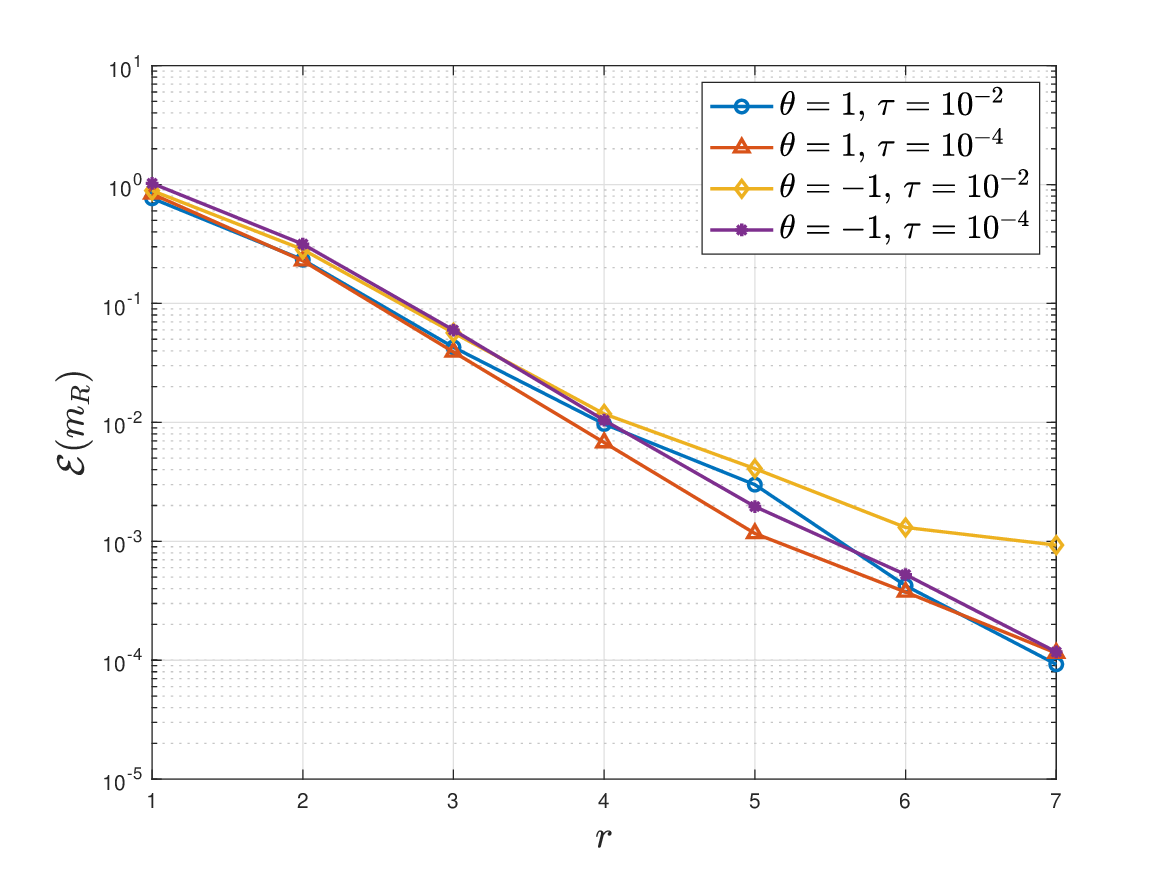}}
    \caption{Test 1: Average $L^2$ error of bi-fidelity approximations for $m_J$ with respect to the number of high-fidelity simulation runs at $\theta=\pm1$ and different $\tau$.}
    \label{Fig.Test1.2.m}
\end{figure*}

Figures~\ref{Fig.Test1.1.m} and Figures~\ref{Fig.Test1.2.m} report the average $L^2$ errors with respect to the number of high-fidelity simulations. In all cases, the errors decrease rapidly when the number of selected high-fidelity samples increases. This shows that the greedy-selected samples capture the main dependence of the solution on the random initial contact moments. The result is consistent with the purpose of the bi-fidelity construction: instead of running the microscopic solver for all random samples, only a small number of carefully selected microscopic simulations are needed to reconstruct the high-fidelity output.

\begin{figure*}
	\centering 
	\subfloat[$\mathbb{E}\big(\rho_S(t,\z)\big)$]{\includegraphics[width=.23\textwidth]{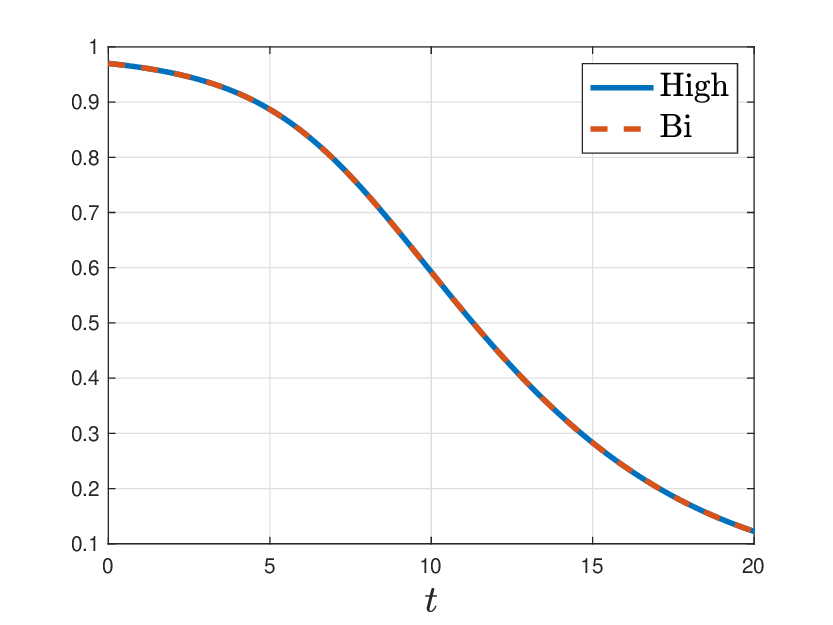}}
    \hspace{1mm}
    \subfloat[$\mathbb{E}\big(\rho_S(t,\z)\big)$]{\includegraphics[width=.23\textwidth]{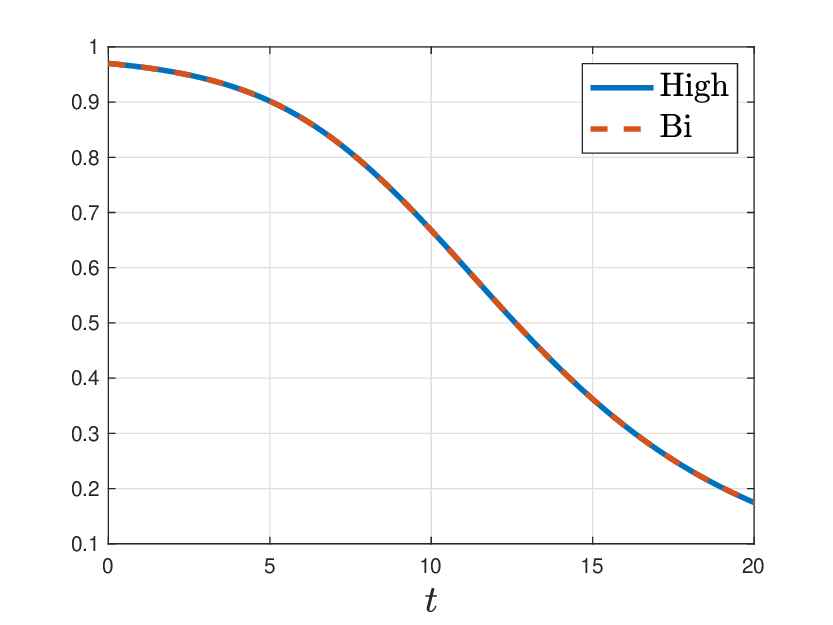}}
    \hspace{1mm}
    \subfloat[$\mathbb{E}\big(\rho_S(t,\z)\big)$]{\includegraphics[width=.23\textwidth]{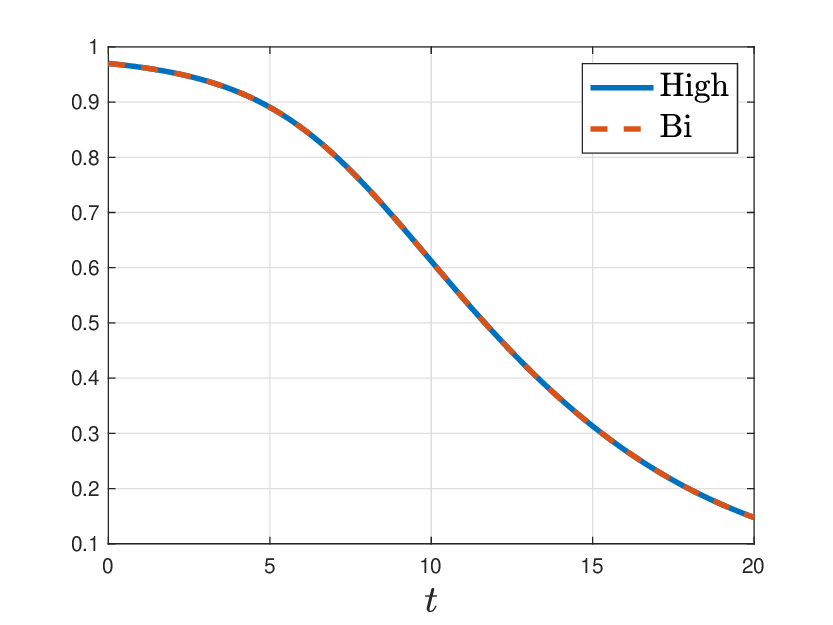}}
    \hspace{1mm}
    \subfloat[$\mathbb{E}\big(\rho_S(t,\z)\big)$]{\includegraphics[width=.23\textwidth]{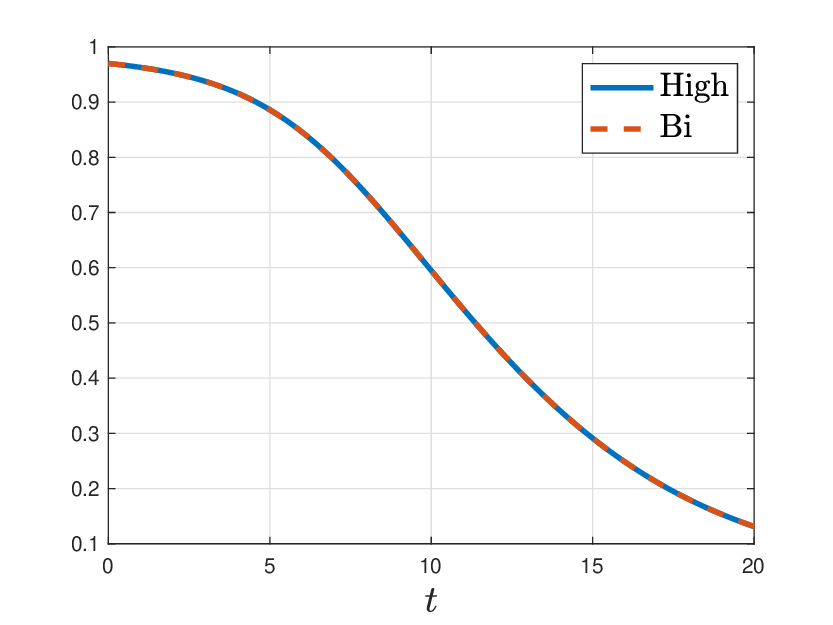}}

    \subfloat[SD$\big(\rho_S(t,\z)\big)$]{\includegraphics[width=.23\textwidth]{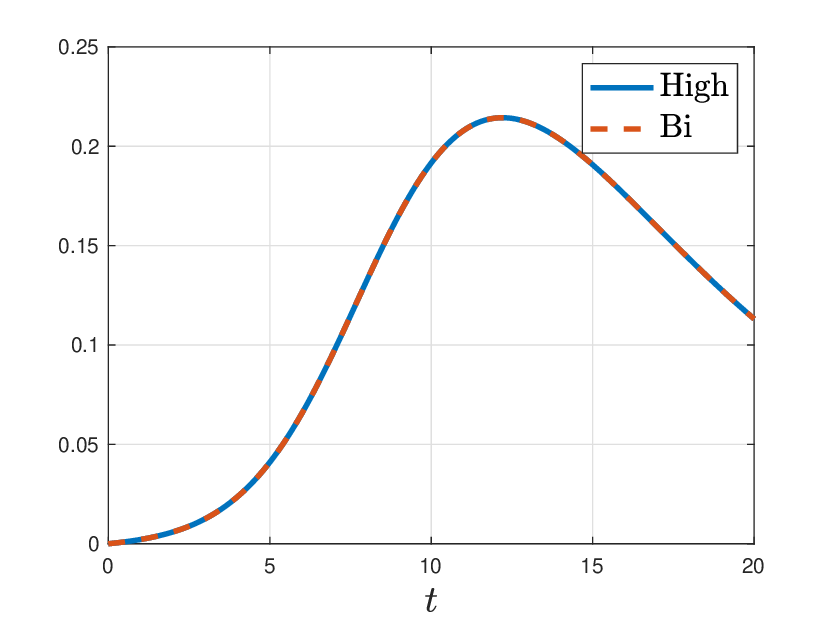}}
    \hspace{1mm}
    \subfloat[SD$\big(\rho_S(t,\z)\big)$]{\includegraphics[width=.23\textwidth]{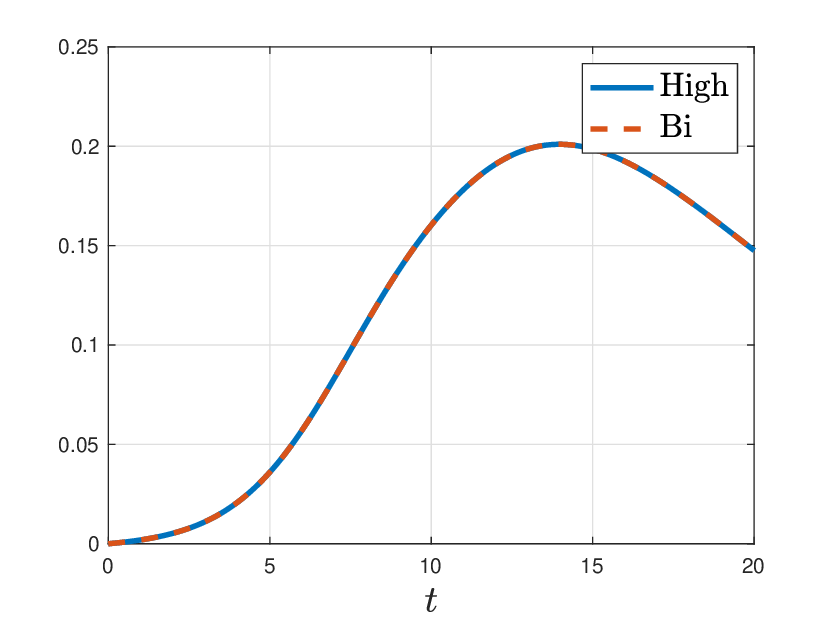}}
    \hspace{1mm}
    \subfloat[SD$\big(\rho_S(t,\z)\big)$]{\includegraphics[width=.23\textwidth]{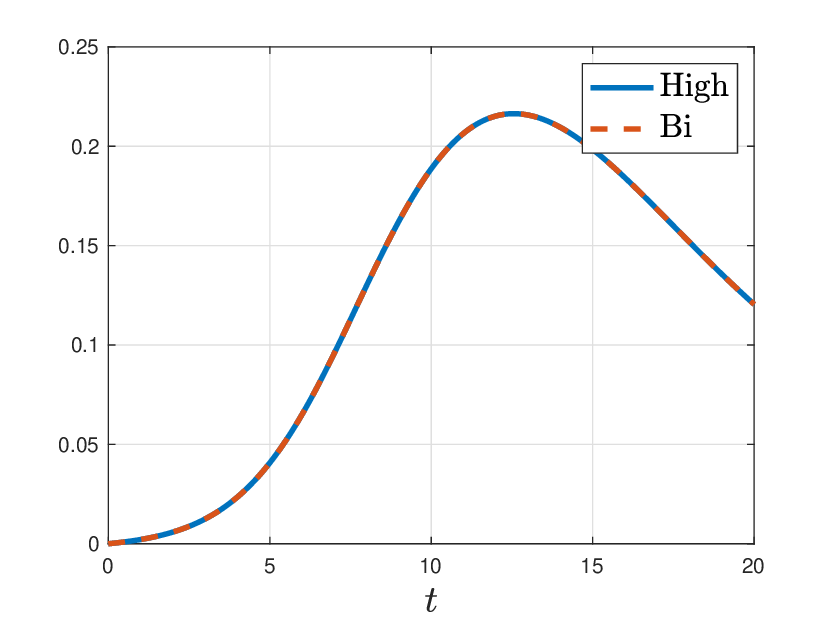}}
    \hspace{1mm}
    \subfloat[SD$\big(\rho_S(t,\z)\big)$]{\includegraphics[width=.23\textwidth]{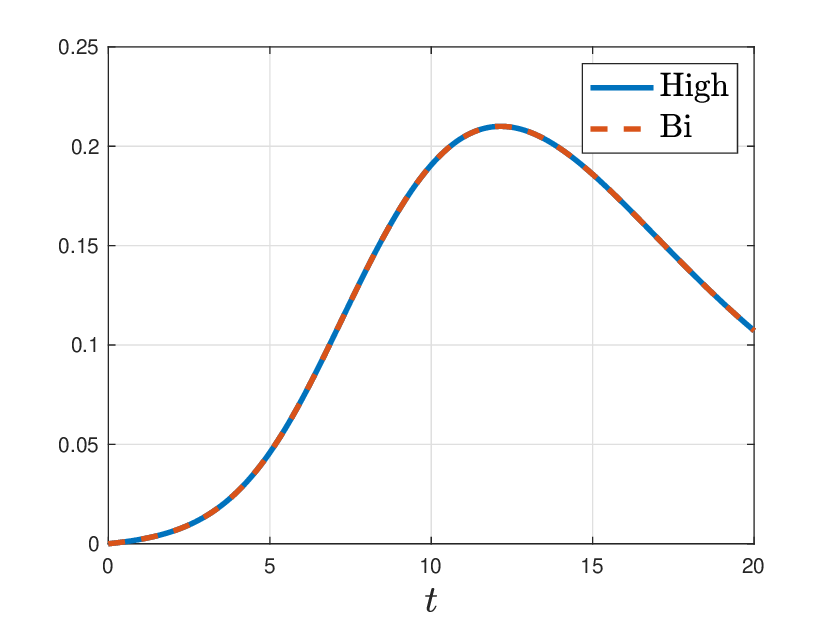}}

    \subfloat[$\mathbb{E}\big(m_S(t,\z)\big)$]{\includegraphics[width=.23\textwidth]{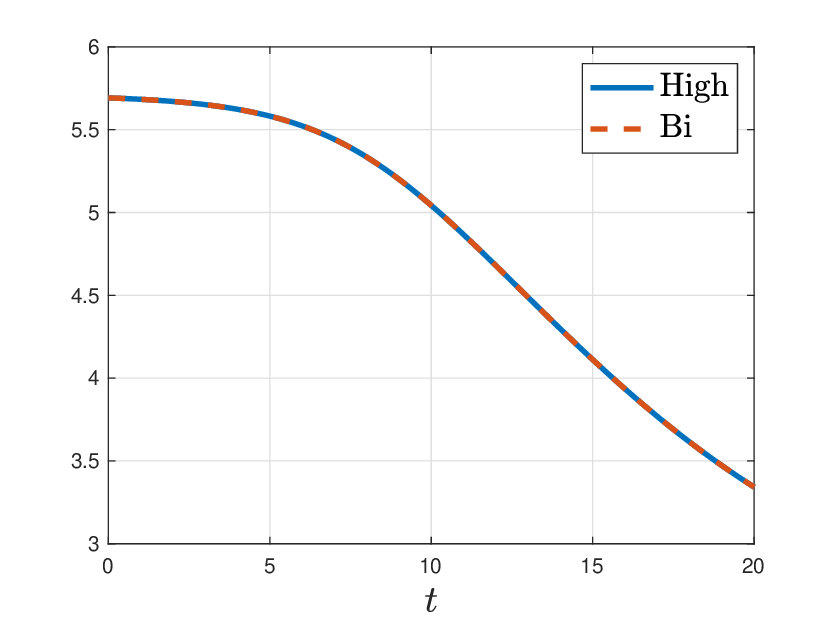}}
    \hspace{1mm}
    \subfloat[$\mathbb{E}\big(m_S(t,\z)\big)$]{\includegraphics[width=.23\textwidth]{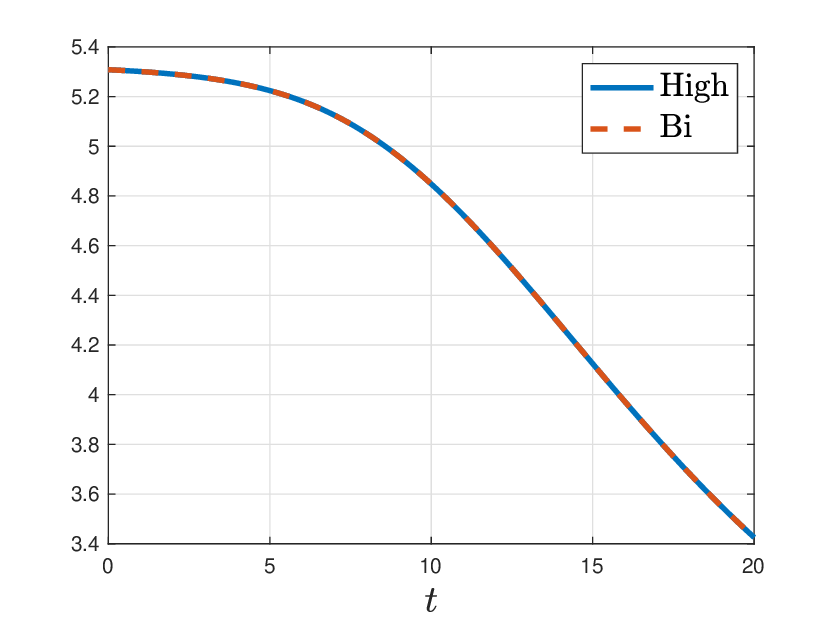}}
    \hspace{1mm}
    \subfloat[$\mathbb{E}\big(m_S(t,\z)\big)$]{\includegraphics[width=.23\textwidth]{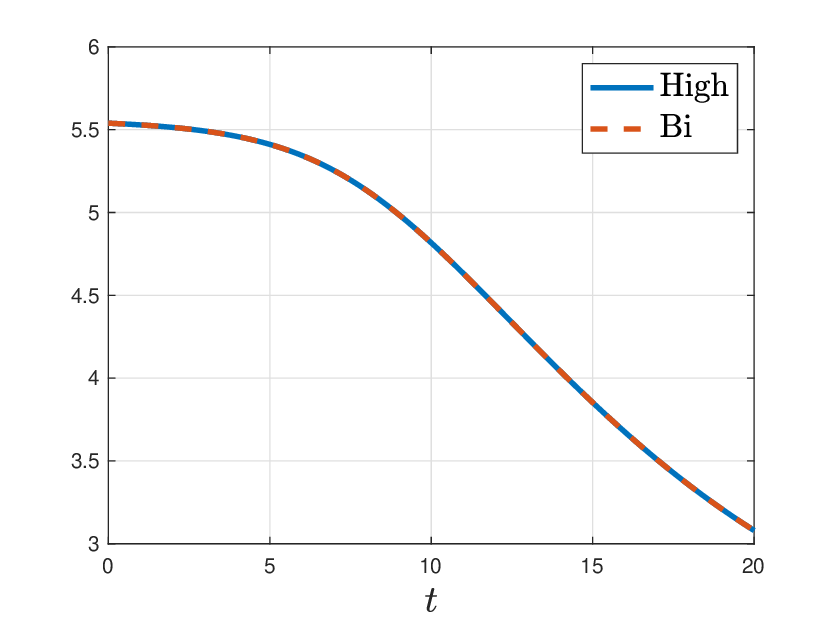}}
    \hspace{1mm}
    \subfloat[$\mathbb{E}\big(m_S(t,\z)\big)$]{\includegraphics[width=.23\textwidth]{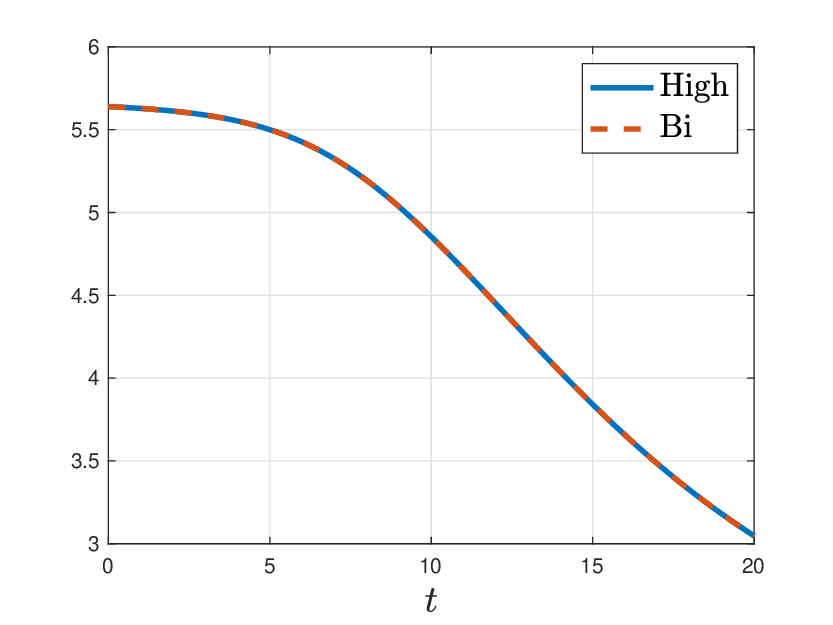}}

    \subfloat[SD$\big(m_S(t,\z)\big)$]{\includegraphics[width=.23\textwidth]{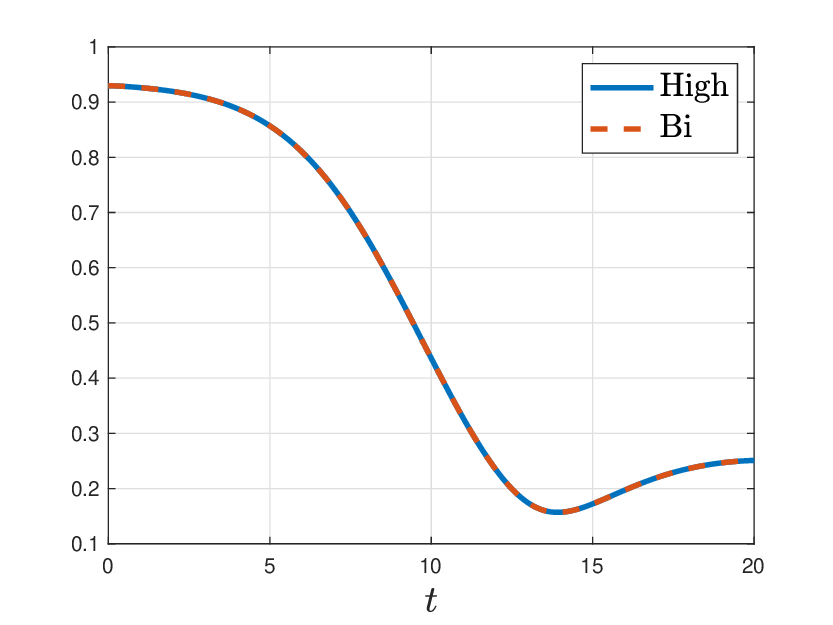}}
    \hspace{1mm}
    \subfloat[SD$\big(m_S(t,\z)\big)$]{\includegraphics[width=.23\textwidth]{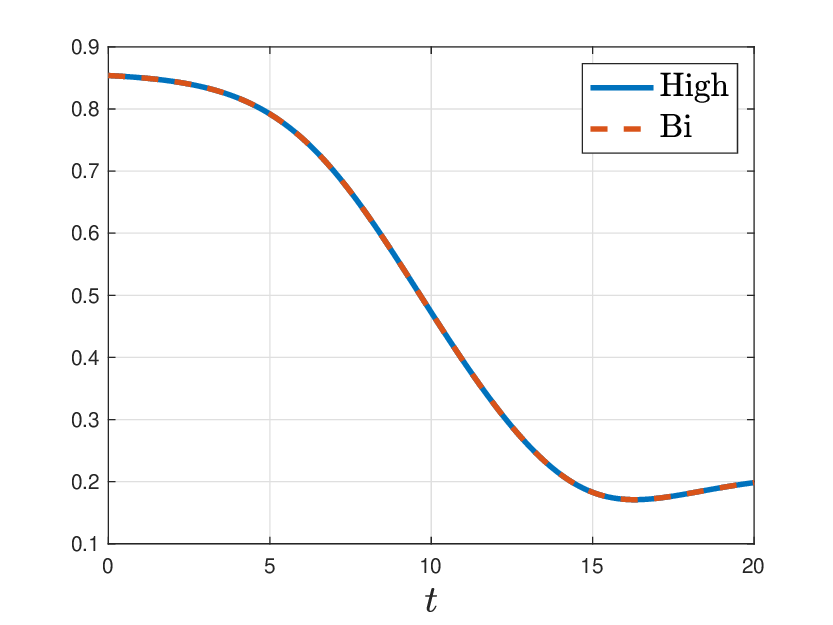}}
    \hspace{1mm}
    \subfloat[SD$\big(m_S(t,\z)\big)$]{\includegraphics[width=.23\textwidth]{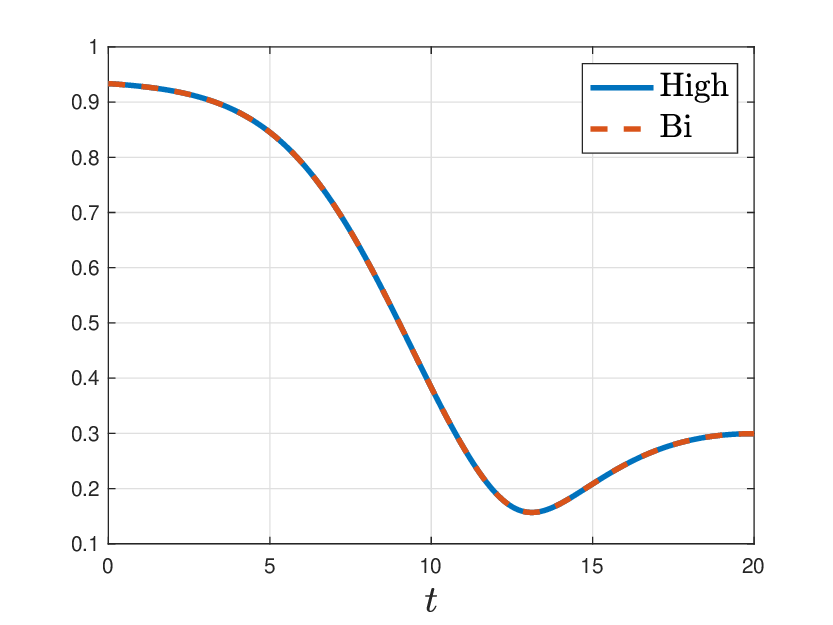}}
    \hspace{1mm}
    \subfloat[SD$\big(m_S(t,\z)\big)$]{\includegraphics[width=.23\textwidth]{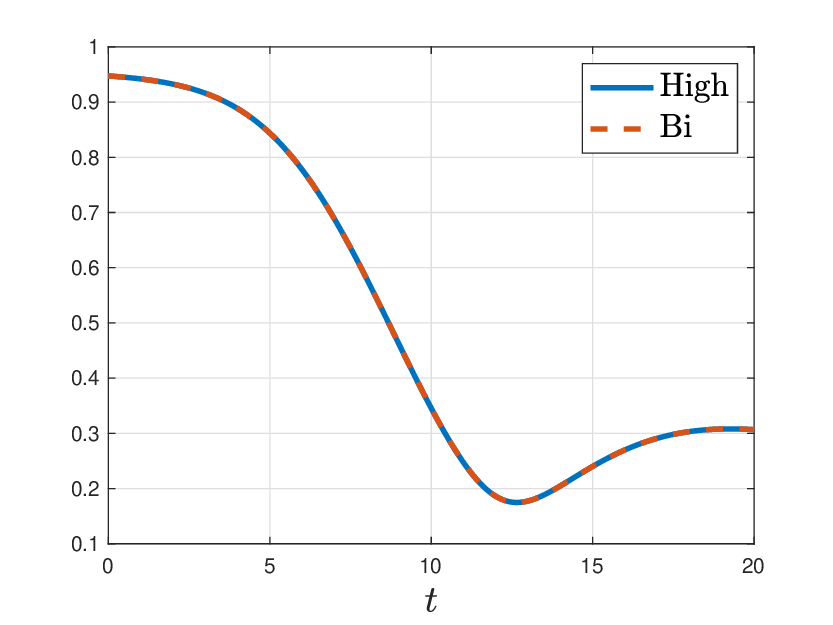}}
	\caption{Test 1: Mean and standard deviation of high- and bi-fidelity solutions of $\rho_S(t,\z)$ and $m_S(t,\z)$ at different $\theta$ and $\tau$. The first column from the left uses $\theta=1$ and $\tau=10^{-2}$. The second column uses $\theta=1$ and $\tau=10^{-4}$. The third column uses $\theta=-1$ and $\tau=10^{-2}$. The fourth column uses $\theta=-1$ and $\tau=10^{-4}$.} 
	\label{Fig.Test1.3.m}
\end{figure*}

Figure~\ref{Fig.Test1.3.m} compares the mean and standard deviation of $\rho_S$ and $m_S$ computed by the high-fidelity and bi-fidelity solvers. The bi-fidelity statistics are obtained using only $7$ high-fidelity simulations. The mean curves are almost indistinguishable at the plotted scale, and the standard deviations also follow the same temporal pattern. This indicates that the proposed method is not only accurate for individual samples, but also reliable for uncertainty propagation. In particular, the method captures both the average epidemic evolution and the spread of the solution caused by random initial contact moments.
The computational saving is significant. When $\tau=10^{-2}$ and $\theta=1$, the high-fidelity microscopic solver costs approximately $3000$ times as much as the low-fidelity macroscopic solver up to the final time $T=20$. In the reported experiment, $500$ runs of the microscopic solver take about $366$ seconds, while the corresponding macroscopic runs take about $0.1$ seconds. This large cost gap is the main reason why the bi-fidelity method is effective: the cheap macroscopic solver can be used many times to explore the random space, while the expensive microscopic solver is only used for a small selected subset.

\begin{figure*}
    \centering
    \subfloat[Solution graph of $\rho_S$ at a certain $\z$]{\includegraphics[width=.45\textwidth]{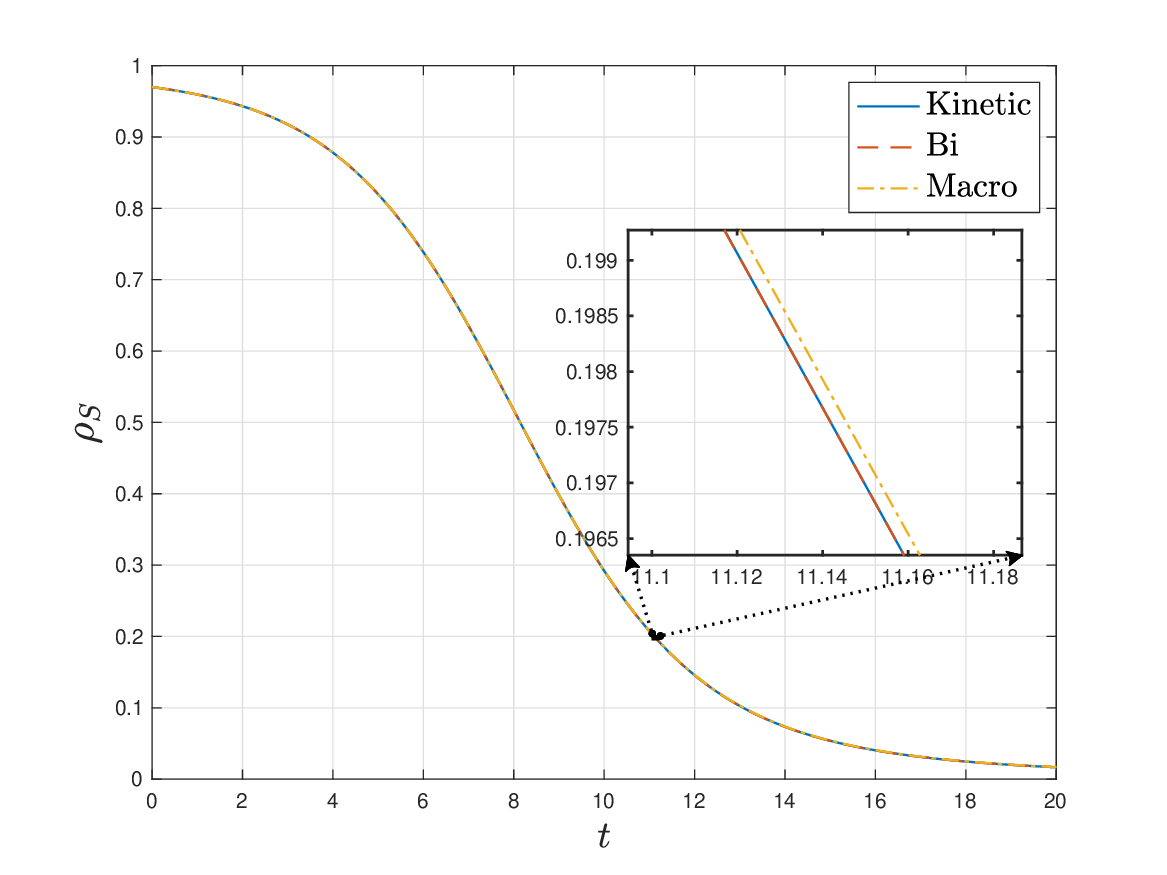}}
    \hspace{5mm}
    \subfloat[Solution graph of $\rho_E$ at a certain $\z$]{\includegraphics[width=.45\textwidth]{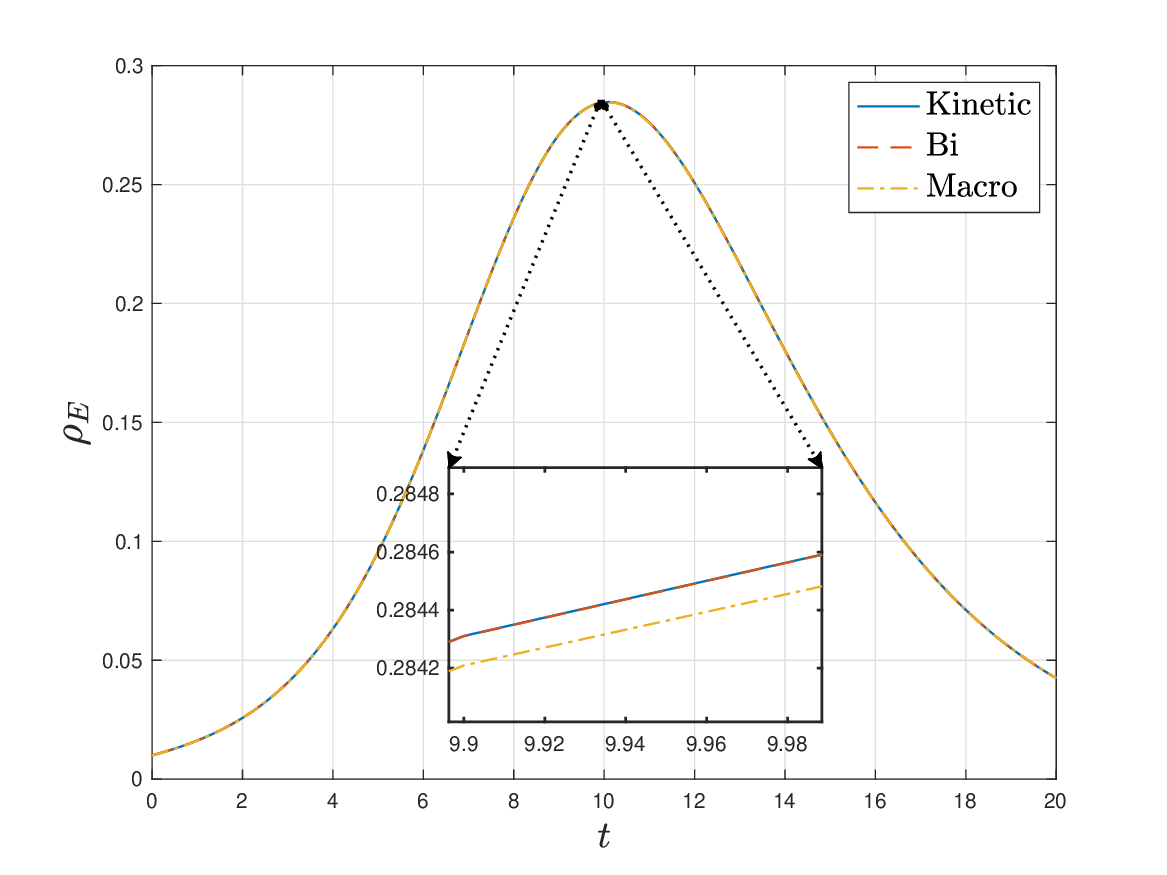}}

    \subfloat[Solution graph of $\rho_I$ at a certain $\z$]{\includegraphics[width=.45\textwidth]{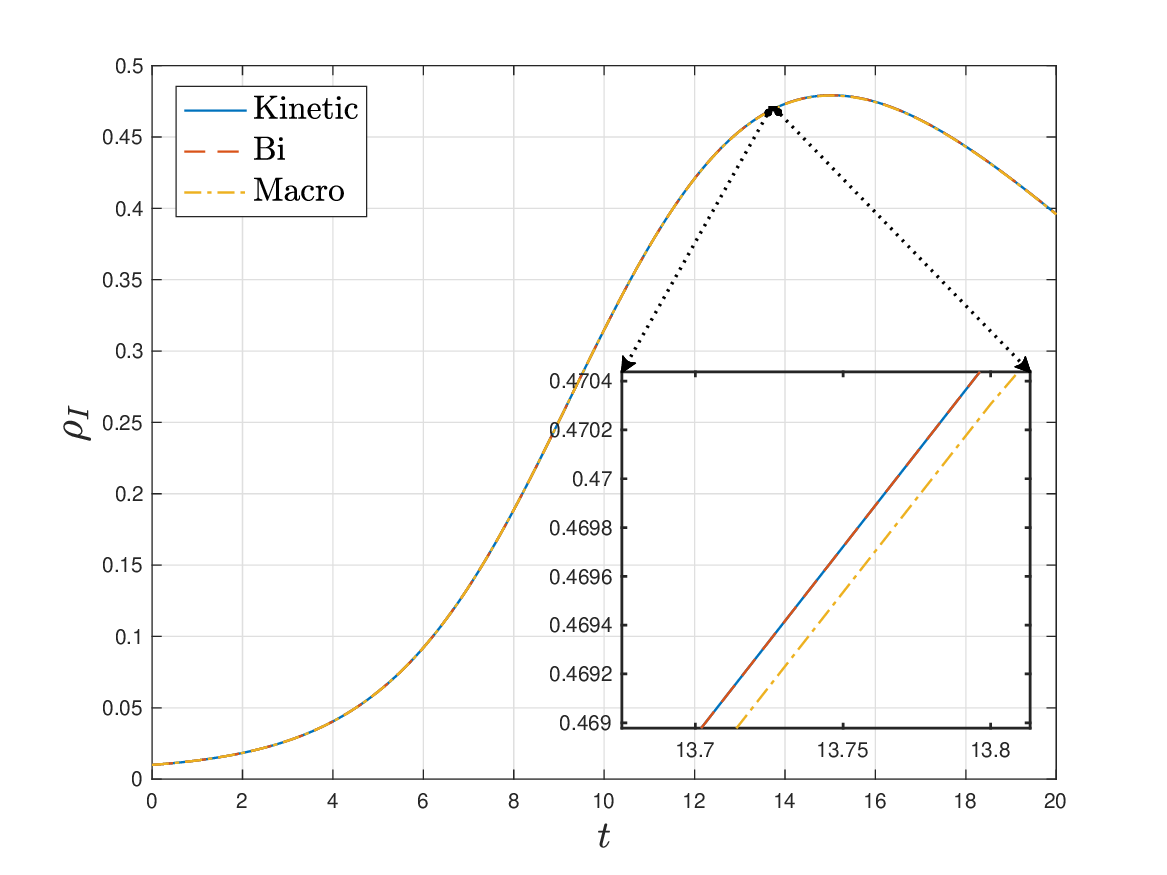}}
    \hspace{5mm}
    \subfloat[Solution graph of $\rho_R$ at a certain $\z$]{\includegraphics[width=.45\textwidth]{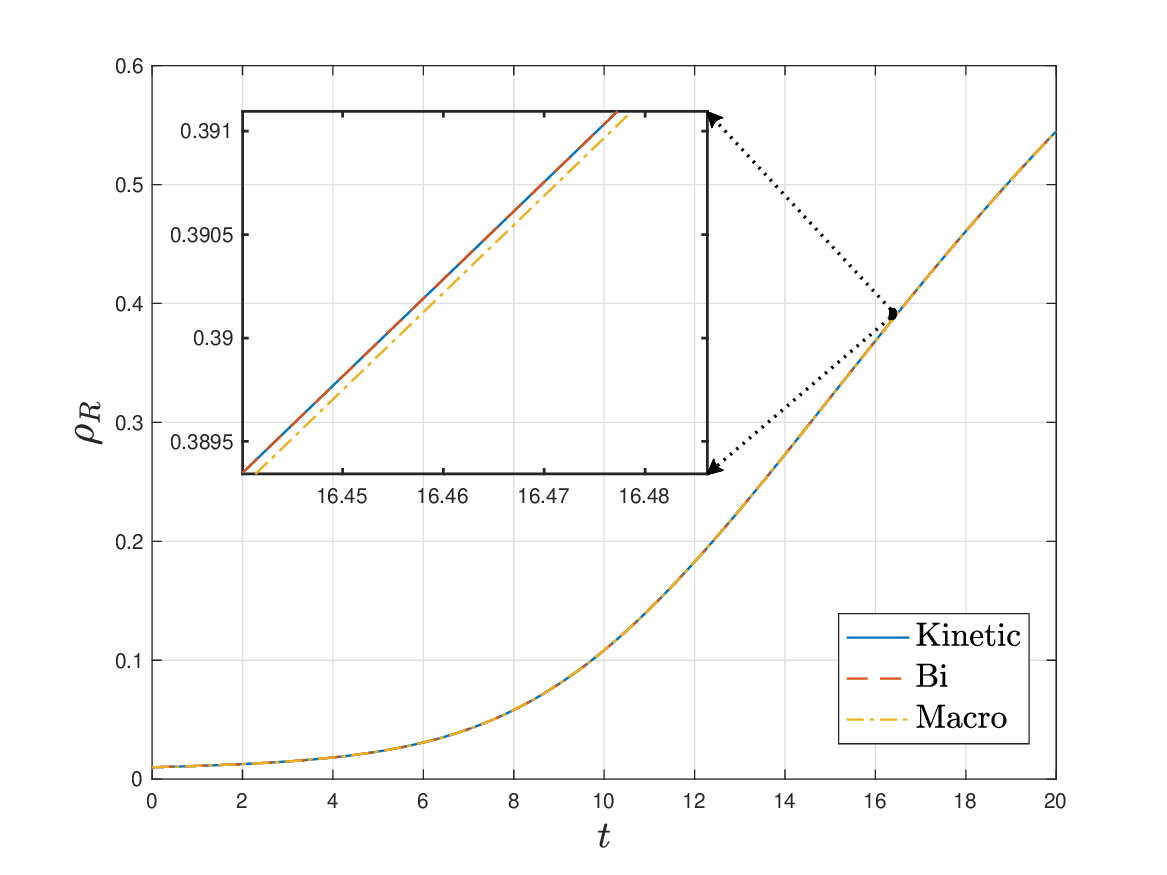}}
    \caption{Test 1: Solution graphs of $\rho_J$ at a certain $\z$ with $\theta=-1$ and $\tau=10^{-4}$.}
    \label{Fig.Test1.4.m}
\end{figure*}

Figures~\ref{Fig.Test1.4.m} shows representative time evolution of $\rho_J$ for a fixed random sample. The low-fidelity solution follows the qualitative trend of the high-fidelity solution, but it is less accurate. In contrast, the bi-fidelity solution almost overlaps with the high-fidelity solution. This confirms the mechanism of the method: the low-fidelity solver provides the correct parametric structure, and the selected high-fidelity snapshots correct the quantitative error. Therefore, even when the macroscopic solver alone is not sufficiently accurate, it can still be useful as a surrogate model inside the bi-fidelity framework.

\subsection{Test 2: Bi-fidelity when $\theta\in[-1,1]$} \label{subsect:test2}
The initial condition of first order moment is specified as
\begin{equation*}
    m_J(t=0,\mathbf{z}) = 
    \begin{cases}\displaystyle
    10\left(1+1.5\sum_{i=1}^{d} \dfrac{z_i\sin z_i}{i}\right),     & \text{if }J=S  \\
    \displaystyle 10\left(1+\sum_{i=1}^{d} \dfrac{z_i\sin z_i}{i}\right),     & \text{if }J=E  \\
    \displaystyle 10\left(1+0.5\sum_{i=1}^{d} \dfrac{z_i\sin z_i}{i}\right),     & \text{if }J=I  \\
    10,     & \text{if }J=R.  
    \end{cases}
\end{equation*}
The initial conditions for the mass fractions and kinetic densities are taken as
\[
\rho_S(0,\z)=0.97,\quad \rho_E(0,\z)=\rho_I(0,\z)=\rho_R(0,\z)=0.01,
\]
while the initial distribution is 
\[
f_J(x,0,\mathbf{z})
=
\frac{\rho_J(0,\z)}{m_J(0,\z)}
\exp\left(-\frac{x}{m_J(0,\z)}\right),
\quad \forall J\in\{S,E,I,R\}.
\]
Compared with the results in Section \ref{subsect:test1} the random perturbation in the initial moments is larger. We set $d=10$ and sample each component independently as $z_i\sim\mathcal U[-1,1]$. The uncertain contact parameter is defined by $\theta(\z)=z_d$, so that $\theta$ may take values in the full interval $[-1,1]$. This setting is more general than Test 1, because the momentum-preserving cases $\theta=\pm1$ are no longer the only possible values. As a result, the macroscopic closure is less direct, and the low-fidelity approximation becomes more challenging. The training set contains $N_z=500$ random samples. The high-fidelity solver uses the following discretization parameters: $50$ time steps on $[0,5]$ and $1000$ grid points on $[0,100]$. 

Recent SEIR-based forecasting studies have emphasized that epidemic parameters may depend on exogenous factors such as environmental conditions, mobility patterns, interventions, and wave-specific characteristics~\cite{colli,viguerie22,ZPPRV2025,VIG_comput}. This motivates us to include uncertainty in several epidemiological parameters in order to represent perturbed transmission and transition rates: $$\beta = 0.025(1 + 0.2 z_1), \qquad \gamma_E = 0.33, \qquad \gamma_I = 0.1. $$

In this regime, the equilibrium distribution \eqref{eq:fqJ} may become numerically ill-conditioned when $\theta$ is close to zero. This difficulty comes from the fact that the expression of the equilibrium density contains terms depending explicitly on $\theta$ in the exponent and in the power of $x$. Therefore, a direct numerical evaluation near $\theta=0$ may be unstable. To avoid this instability, we use the limiting form
\begin{equation}\label{eq:fqJ0}
    f_{J,0}^q(x,t,\z)
    =
    C_{\mu,\sigma^2,m_J}
    x^{-3/2}
    \exp\left\{
    -\frac{\mu}{2\sigma^2}
    \log^2\left(\frac{x}{m_J(t)}\right)
    \right\}.
\end{equation}
In the low-fidelity solver, we introduce a threshold $R_i>0$. When $\theta\in[-R_i,R_i]$, we replace \eqref{eq:fqJ} by the limiting density \eqref{eq:fqJ0}; otherwise, we use the original equilibrium density \eqref{eq:fqJ}. This treatment is only used to stabilize the numerical evaluation near $\theta=0$ and does not change the high-fidelity microscopic solver.

\begin{figure*} 
    \centering
    \subfloat[Average $L^2$ error of bi-fidelity approximations for $\rho_S$]{\includegraphics[width=.45\textwidth]{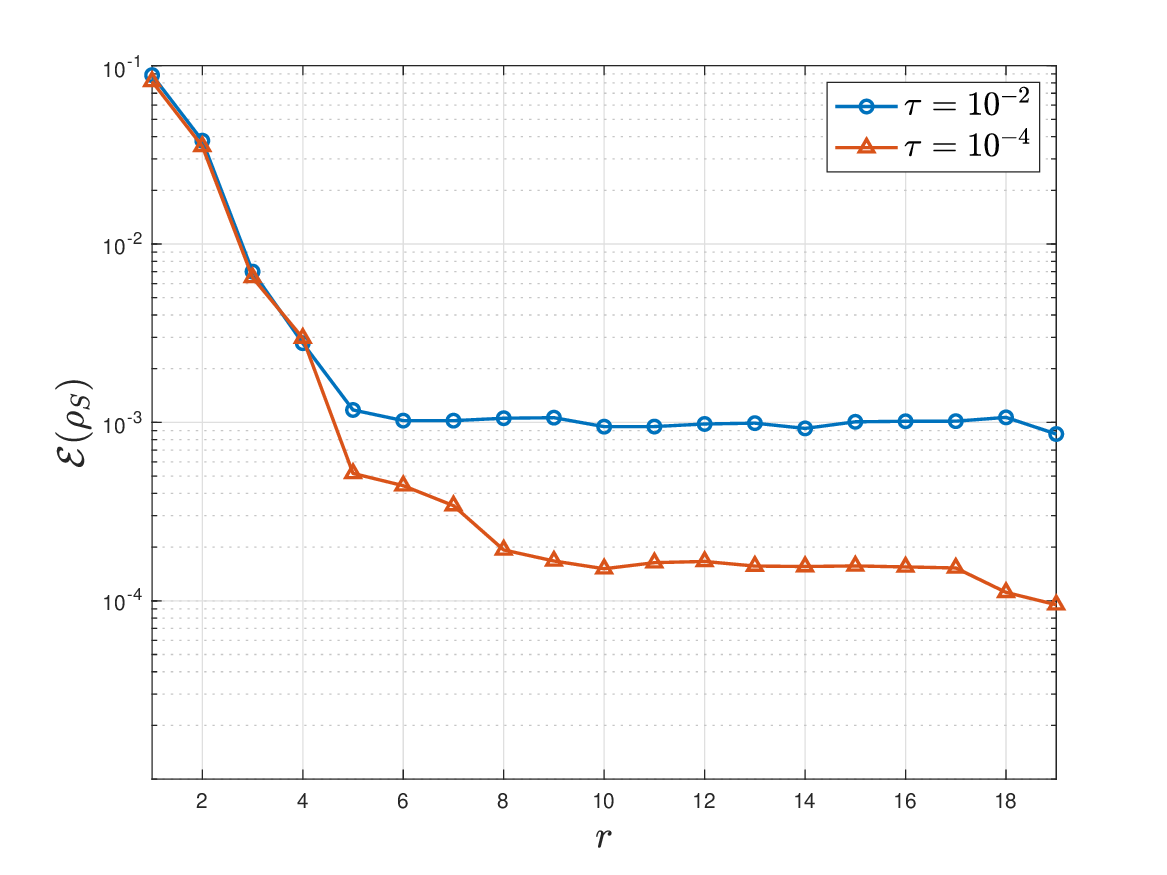}}
    \hspace{5mm}
    \subfloat[Average $L^2$ error of bi-fidelity approximations for $\rho_E$]{\includegraphics[width=.45\textwidth]{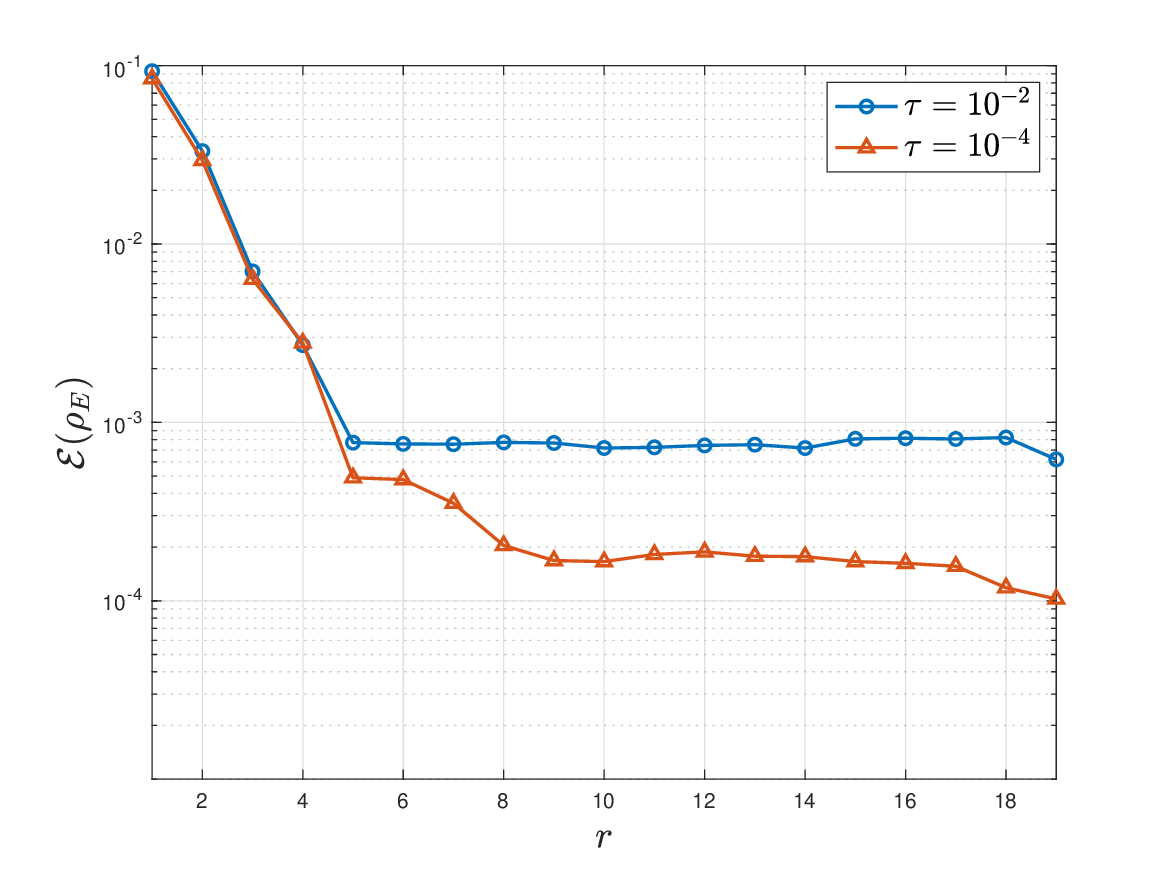}}

    \subfloat[Average $L^2$ error of bi-fidelity approximations for $\rho_I$]{\includegraphics[width=.45\textwidth]{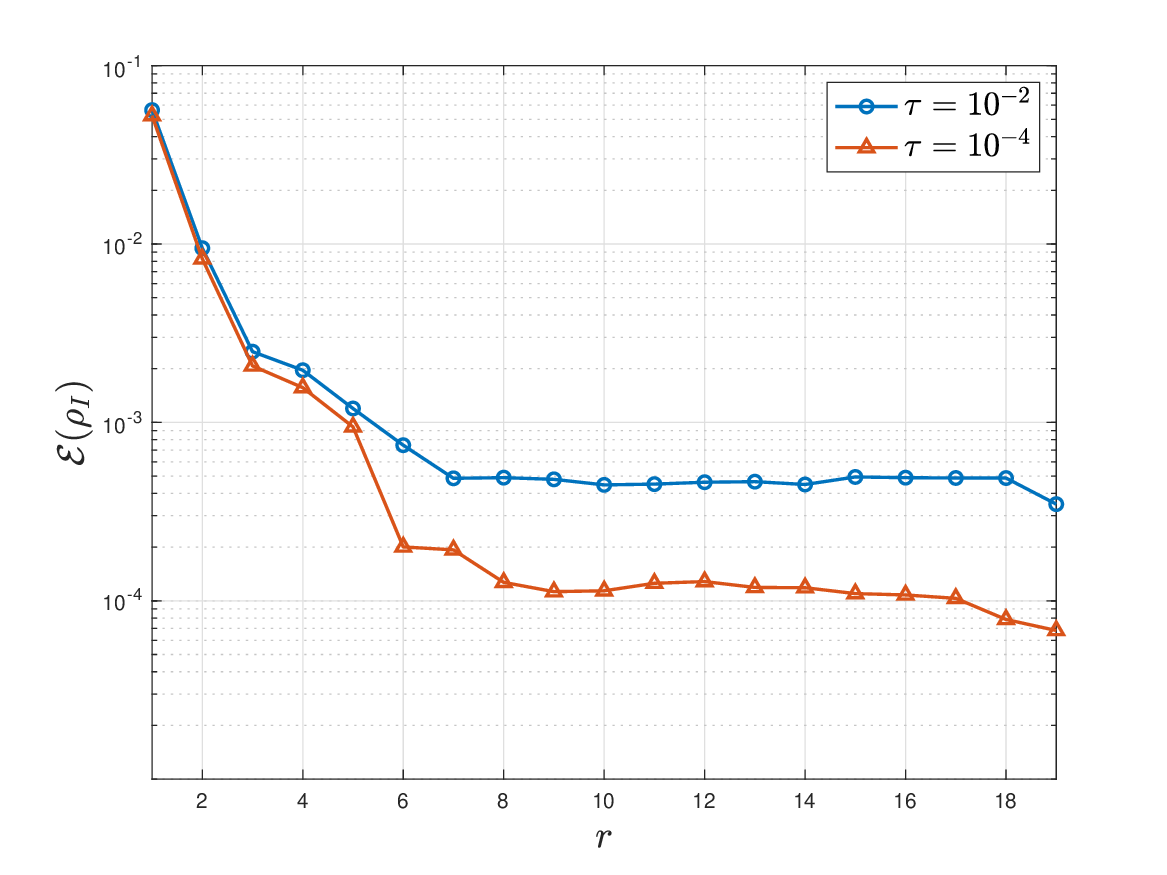}}
    \hspace{5mm}
    \subfloat[Average $L^2$ error of bi-fidelity approximations for $\rho_R$]{\includegraphics[width=.45\textwidth]{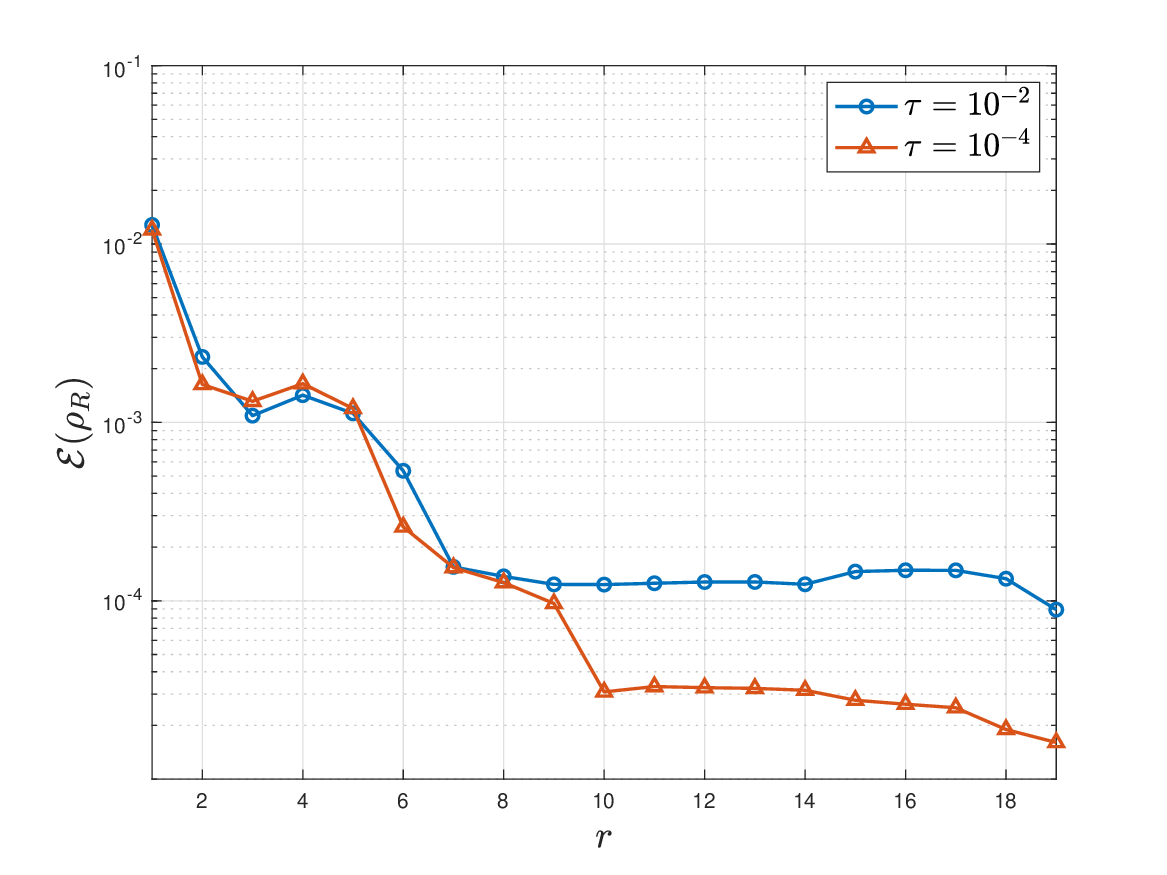}}
    \caption{Test 2: Average $L^2$ error of bi-fidelity approximations for $\rho_J$ with respect to the number of high-fidelity simulation runs at different $\tau$.}
    \label{Fig.Test2.1.m}
\end{figure*}

We test two values of $\tau$ in order to compare different relaxation regimes by using $R_i=5\times10^{-3}$. Figure~\ref{Fig.Test2.1.m} reports the average $L^2$ errors between the bi-fidelity and high-fidelity solutions. Although this test is more difficult than Test 1, the error still decreases as the number of high-fidelity samples increases. This means that the low-fidelity model still provides useful information about the random dependence of the microscopic solution, even when $\theta$ is not restricted to the momentum-preserving values $\pm1$. The error decay also shows that the greedy-selected samples remain representative in the higher-dimensional uncertain setting with both uncertain initial moments and uncertain model parameters.

\begin{figure*} 
	\centering 
	\subfloat[$\mathbb{E}\big(\rho_S(t,\z)\big)$ when $\tau=10^{-2}$]{\includegraphics[width=.4\textwidth]{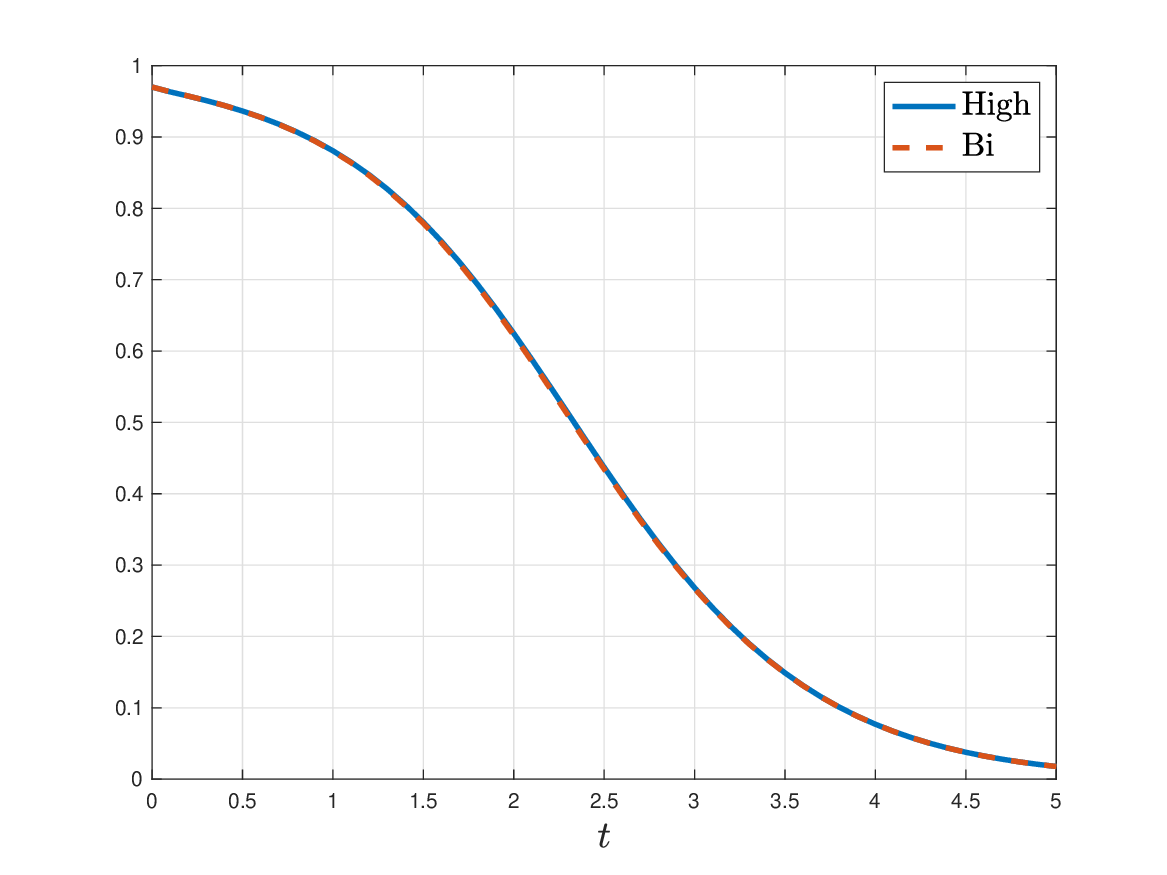}}
    \hspace{5mm}
    \subfloat[SD$\big(\rho_S(t,\z)\big)$ when $\tau=10^{-2}$]{\includegraphics[width=.4\textwidth]{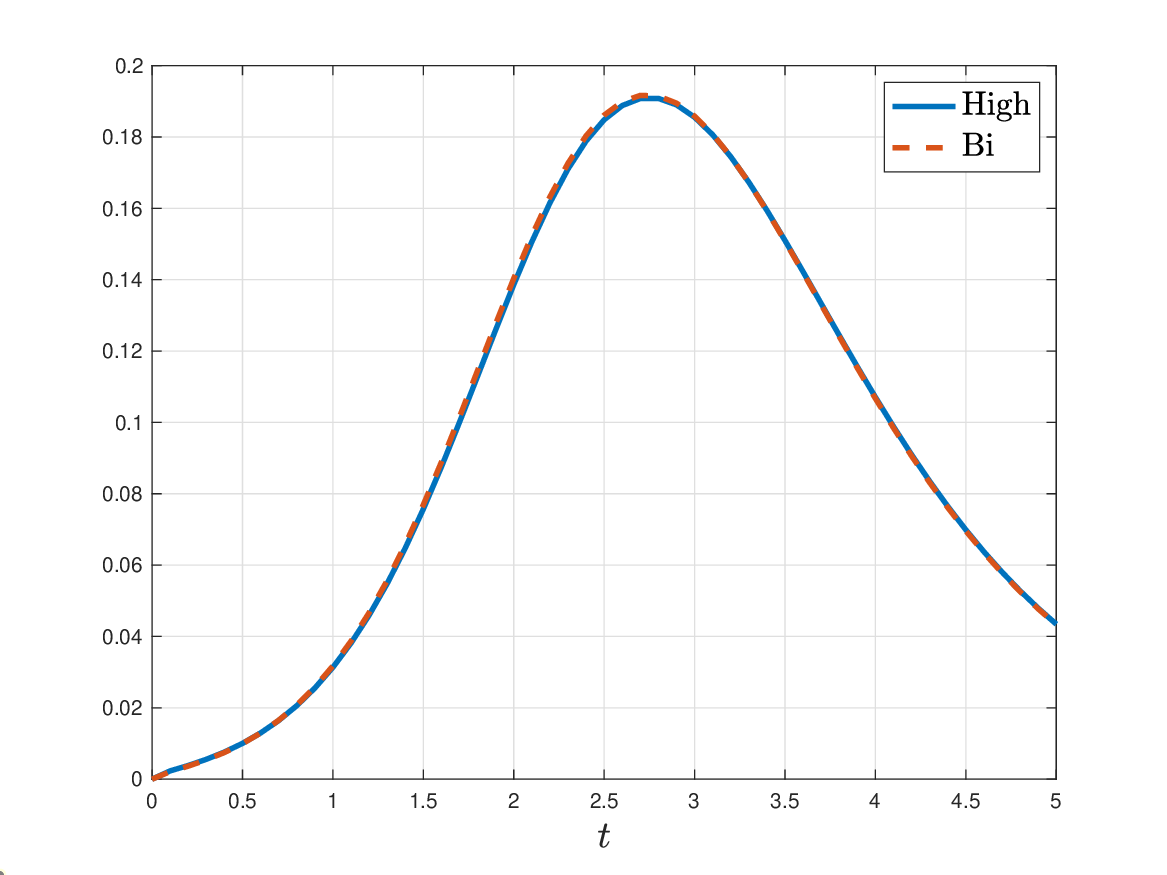}}

    \subfloat[$\mathbb{E}\big(\rho_S(t,\z)\big)$ when $\tau=10^{-4}$]{\includegraphics[width=.4\textwidth]{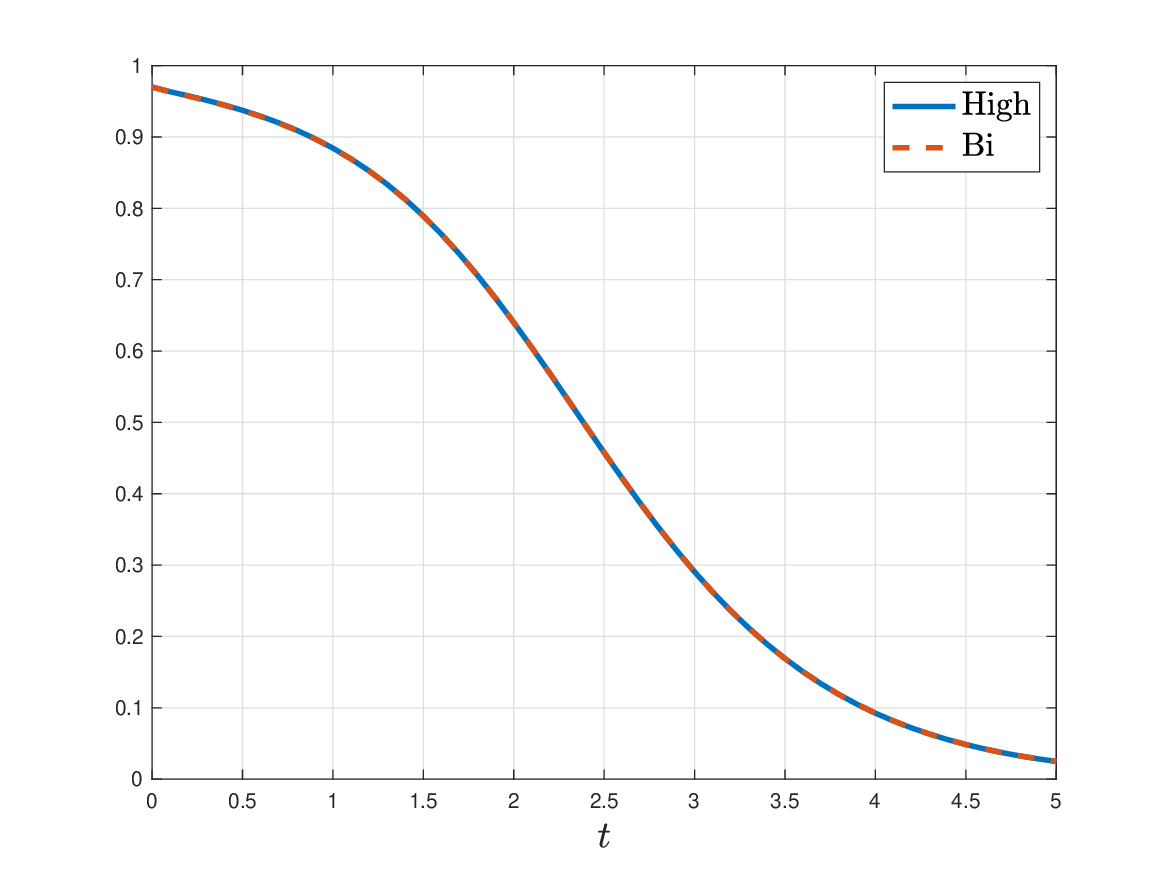}}
    \hspace{5mm}
    \subfloat[SD$\big(\rho_S(t,\z)\big)$ when $\tau=10^{-4}$]{\includegraphics[width=.4\textwidth]{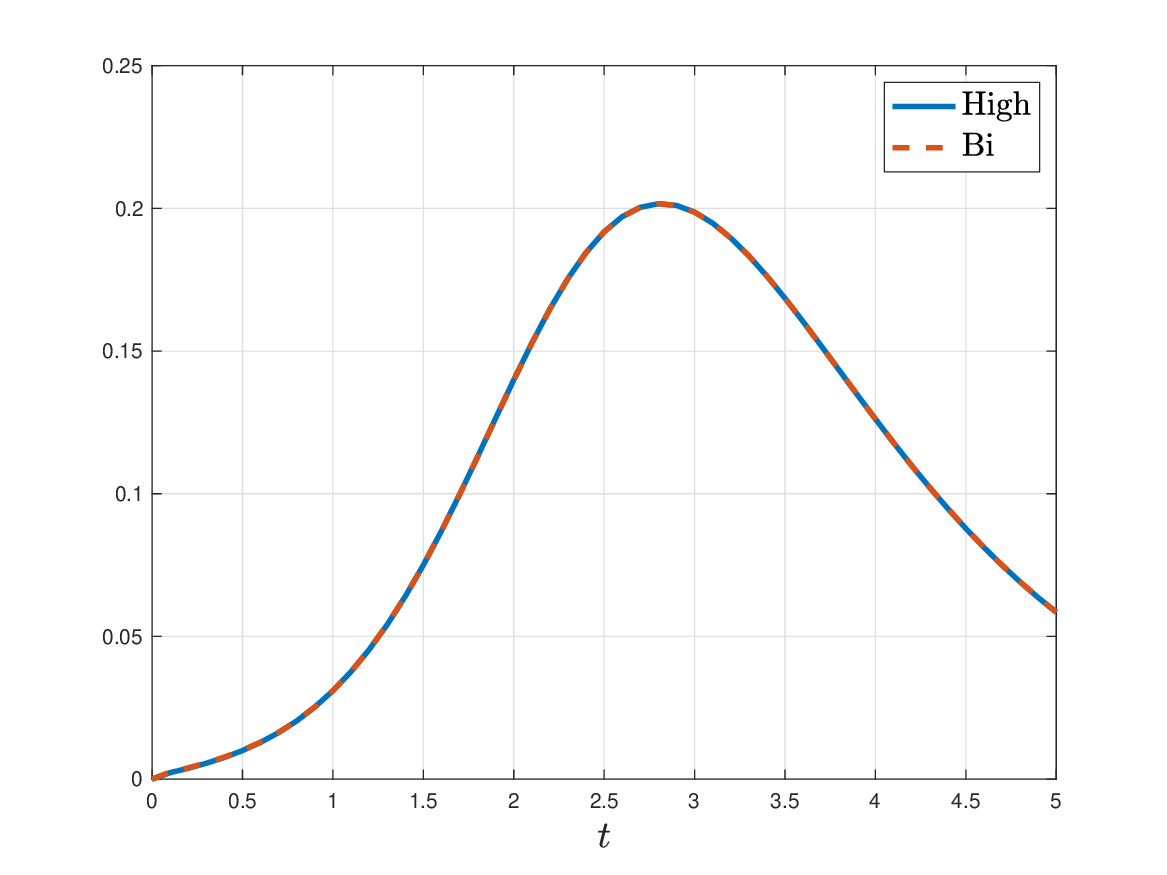}}
    
	\caption{Test 2: Mean and standard deviation of bi-fidelity solutions of $\rho_S(t,\z)$ at different $\tau$. } 
	\label{Fig.Test2.2.m}
\end{figure*}

Figure~\ref{Fig.Test2.2.m} shows the mean and standard deviation of $\rho_S$. Using only $19$ high-fidelity runs, the bi-fidelity approximation remains close to the high-fidelity reference. The agreement of the mean indicates that the method correctly captures the average decay of the susceptible population. The agreement of the standard deviation indicates that the method also reproduces the uncertainty level induced by the random inputs. This point is important because an accurate mean alone would not be sufficient for uncertainty quantification.

\begin{figure*}
    \centering
    \subfloat[Solution graph of $\rho_S$ at a certain $\z$]{\includegraphics[width=.45\textwidth]{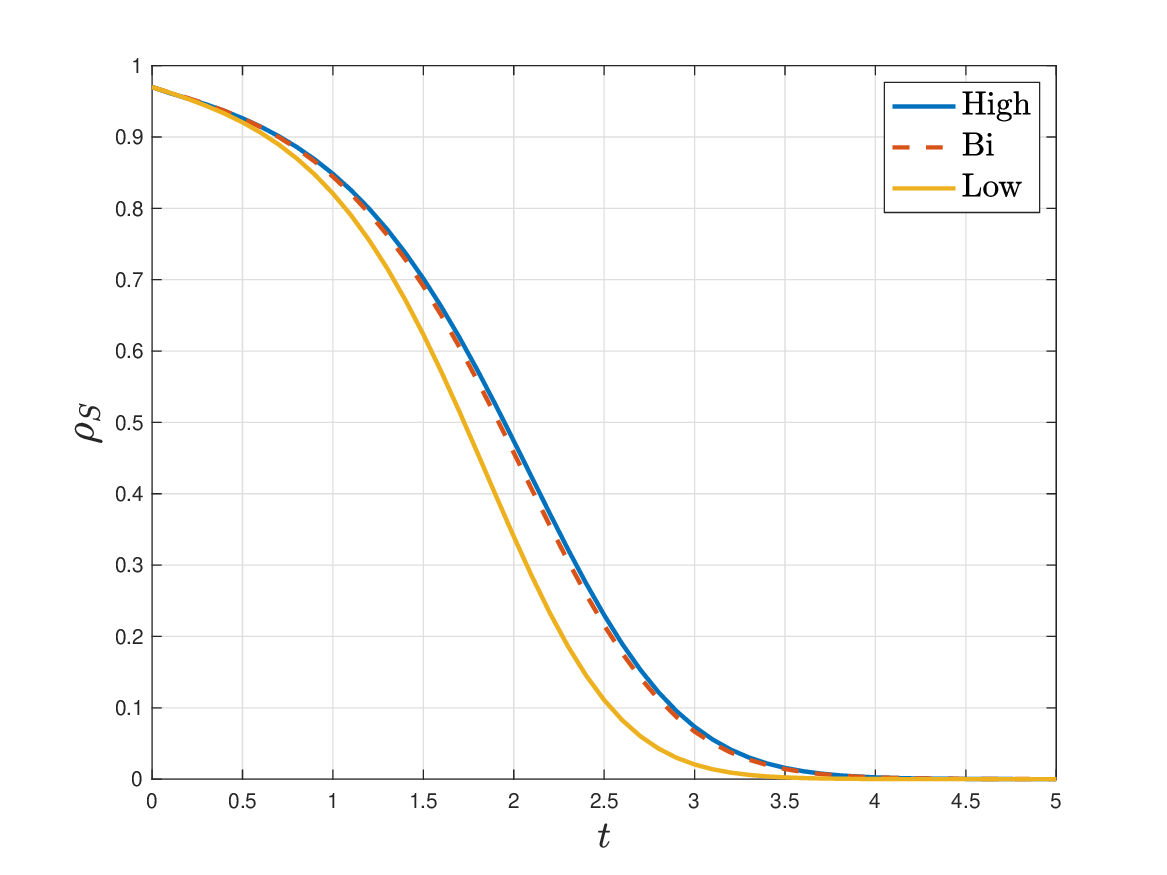}}
    \hspace{5mm}
    \subfloat[Solution graph of $\rho_E$ at a certain $\z$]{\includegraphics[width=.45\textwidth]{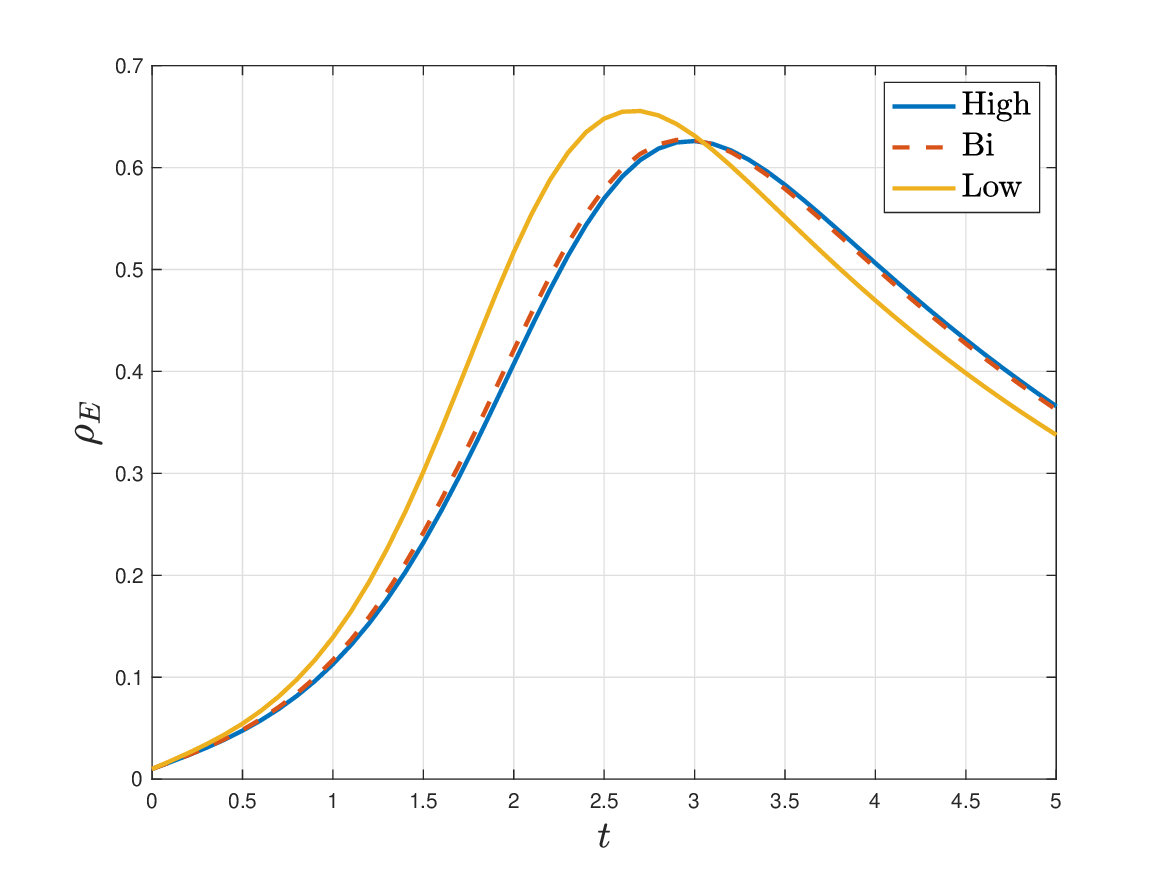}}

    \subfloat[Solution graph of $\rho_I$ at a certain $\z$]{\includegraphics[width=.45\textwidth]{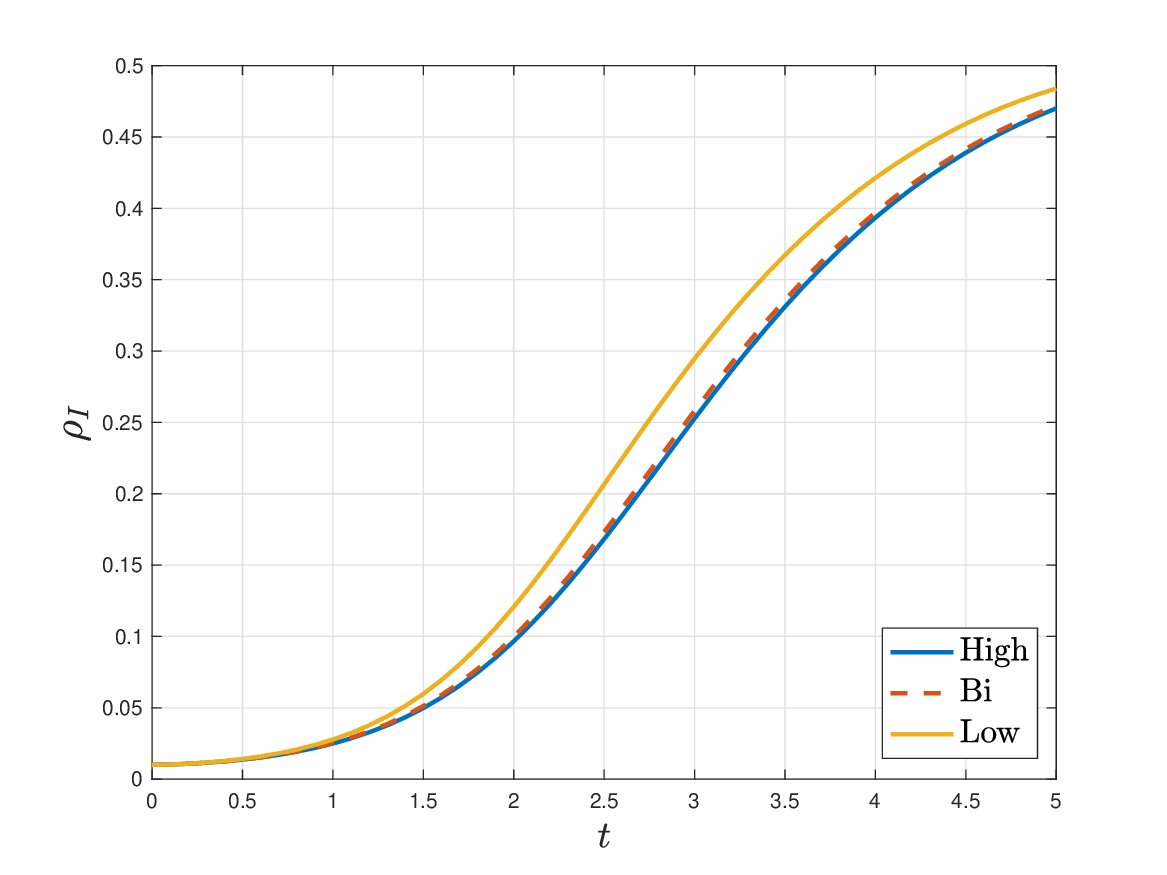}}
    \hspace{5mm}
    \subfloat[Solution graph of $\rho_R$ at a certain $\z$]{\includegraphics[width=.45\textwidth]{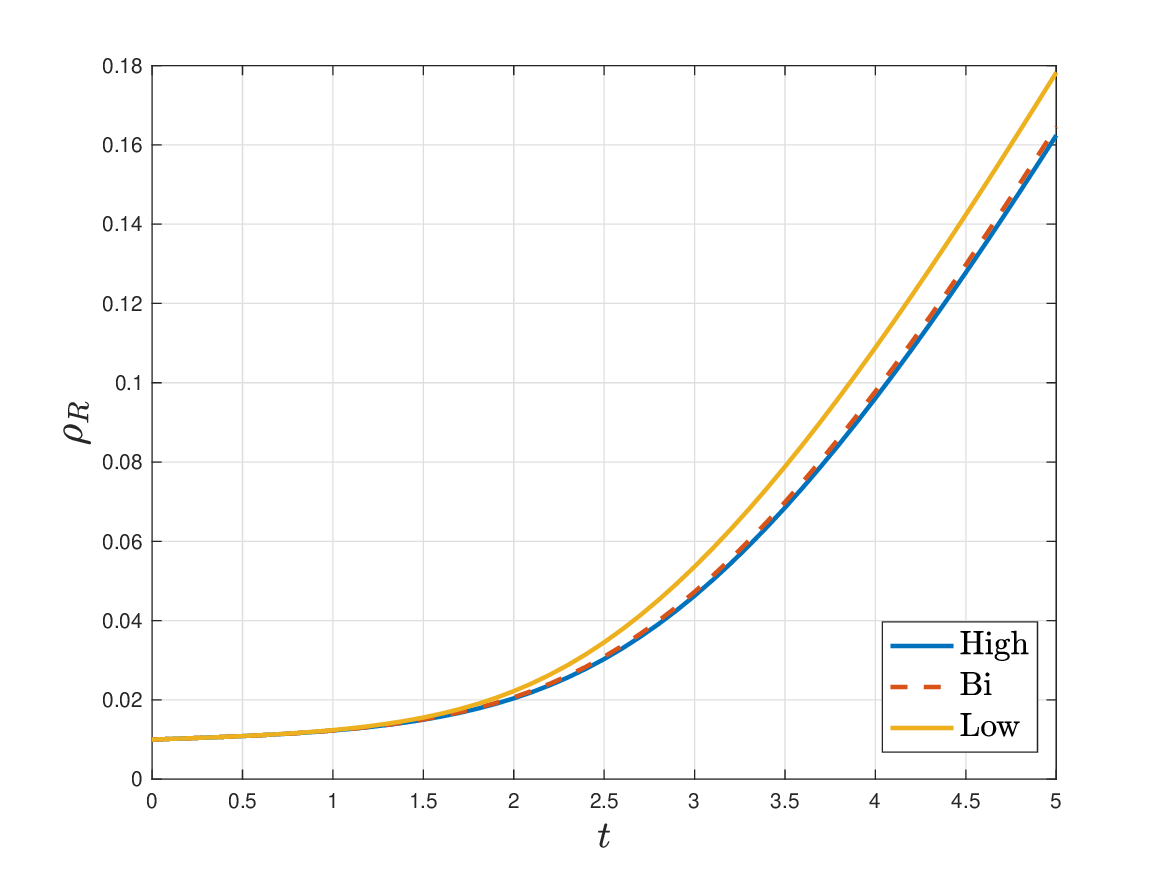}}
    \caption{Test 2: Solution graphs of $\rho_J$ at a certain $\z$ with $\tau=10^{-2}$.}
    \label{Fig.Test2.3.m}
\end{figure*}

When $\tau=10^{-4}$, the high-fidelity microscopic solver costs approximately twice as much as the low-fidelity macroscopic solver up to the final time $T=5$. The cost gap is smaller than in Test 1, but the multi-fidelity method is still useful because it reduces the number of required high-fidelity simulations. Figure~\ref{Fig.Test2.3.m} shows representative trajectories of $\rho_S$ and $\rho_E$ for one random sample. The low-fidelity solver captures the qualitative change of the solution, but it has a visible quantitative discrepancy. The bi-fidelity solution corrects this discrepancy and follows the high-fidelity solution more closely. This confirms that the method does not require the low-fidelity solver to be highly accurate by itself; it only needs the low-fidelity solution space to contain enough information about the parametric structure.

\subsection{Test 3: Tri-fidelity algorithm} \label{subsect:test3}

In Test 3, we apply the tri-fidelity method to the same type of uncertain epidemic problem. The purpose is to compare the bi-fidelity and tri-fidelity constructions when three levels of models are available. The high-fidelity solver is still the microscopic kinetic model \eqref{eq:kinetic}--\eqref{eq:QJ}. The medium-fidelity solver is the macroscopic model \eqref{eq:macro1} together with the moment equations \eqref{eq:macro2}. The low-fidelity solver is a cheaper simplified macroscopic model used mainly for the greedy point selection.

The initial condition of the first-order moment is
\begin{equation*}
    m_J(t=0,\mathbf{z}) = 
    \begin{cases}\displaystyle
    10\left(1+1.5\sum_{i=1}^{d} \dfrac{z_i\sin z_i}{i}\right),     & \text{if }J=S,  \\[0.4em]
    \displaystyle 10\left(1+\sum_{i=1}^{d} \dfrac{z_i\sin z_i}{i}\right),     & \text{if }J=E,  \\[0.4em]
    \displaystyle 10\left(1+0.5\sum_{i=1}^{d} \dfrac{z_i\sin z_i}{i}\right),     & \text{if }J=I,  \\[0.4em]
    10,     & \text{if }J=R.
    \end{cases}
\end{equation*}
The initial mass fractions and kinetic densities are given by
\[
\begin{cases} 
\rho_S(0,\z)=0.97,\quad \rho_E(0,\z)=\rho_I(0,\z)=\rho_R(0,\z)=0.01,\\[0.4em]
\displaystyle
f_J(x,0,\mathbf{z})
=
\frac{1}{m_J(0,\z)}
\exp\left(-\frac{x}{m_J(0,\z)}\right),
\quad \forall J\in\{S,E,I,R\}.
\end{cases}
\]

The epidemic constants are chosen as: $$\beta = 0.025(1+0.2z_1),\quad \gamma_E = 0.33(1+0.2z_2),\quad \gamma_I = 0.1(1+0.2z_3).$$

Other initial data are the same as those used in Test 2. This makes the comparison between the bi-fidelity and tri-fidelity methods clearer, since the improvement or difference comes from the fidelity hierarchy rather than from a different initial condition.

The construction of the three fidelity levels is as follows. The high-fidelity solver is the full kinetic model \eqref{eq:kinetic} with the Fokker--Planck operator \eqref{eq:QJ}. The medium-fidelity solver is the macroscopic model \eqref{eq:macro1} together with the corresponding moment equations \eqref{eq:macro2}. This solver is cheaper than the kinetic solver, but it still keeps more information than the simplified low-fidelity model. The low-fidelity solver is obtained by using the macroscopic equations and fixing the values of
\[
\int_{\mathbb{R}^+}xQ_J(f_J)(t,x,\z)\,\ud x,
\qquad J\in\{E,I,R\},
\]
at their initial values. This approximation is not expected to be highly accurate. Its role is to provide a very cheap way to explore the random space and to select important parameter samples. The more accurate medium-fidelity solver is then used to compute the projection coefficients, which are finally transferred to the high-fidelity space.

\begin{figure*}
    \centering
    \subfloat[Average $L^2$ error of tri-fidelity approximations for $\rho_S$]{\includegraphics[width=.45\textwidth]{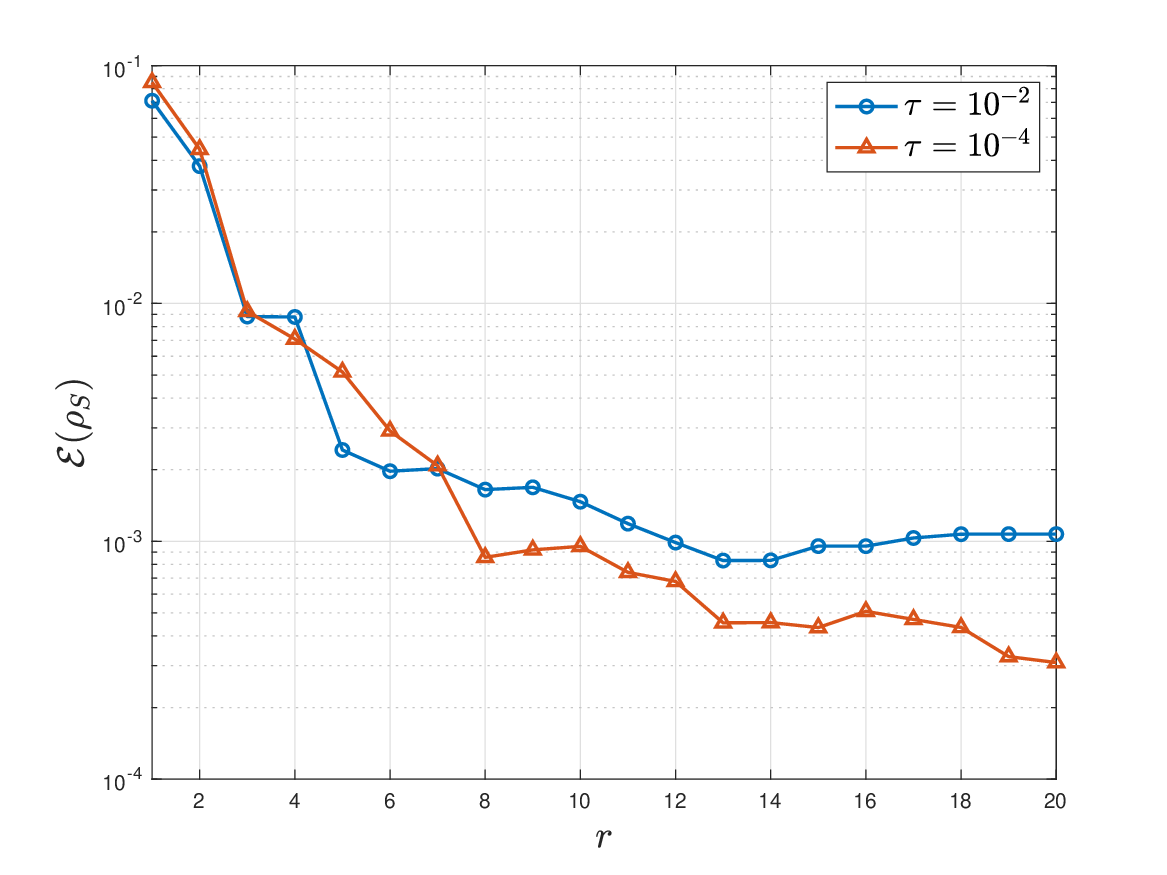}}
    \hspace{5mm}
    \subfloat[Average $L^2$ error of tri-fidelity approximations for $\rho_E$]{\includegraphics[width=.45\textwidth]{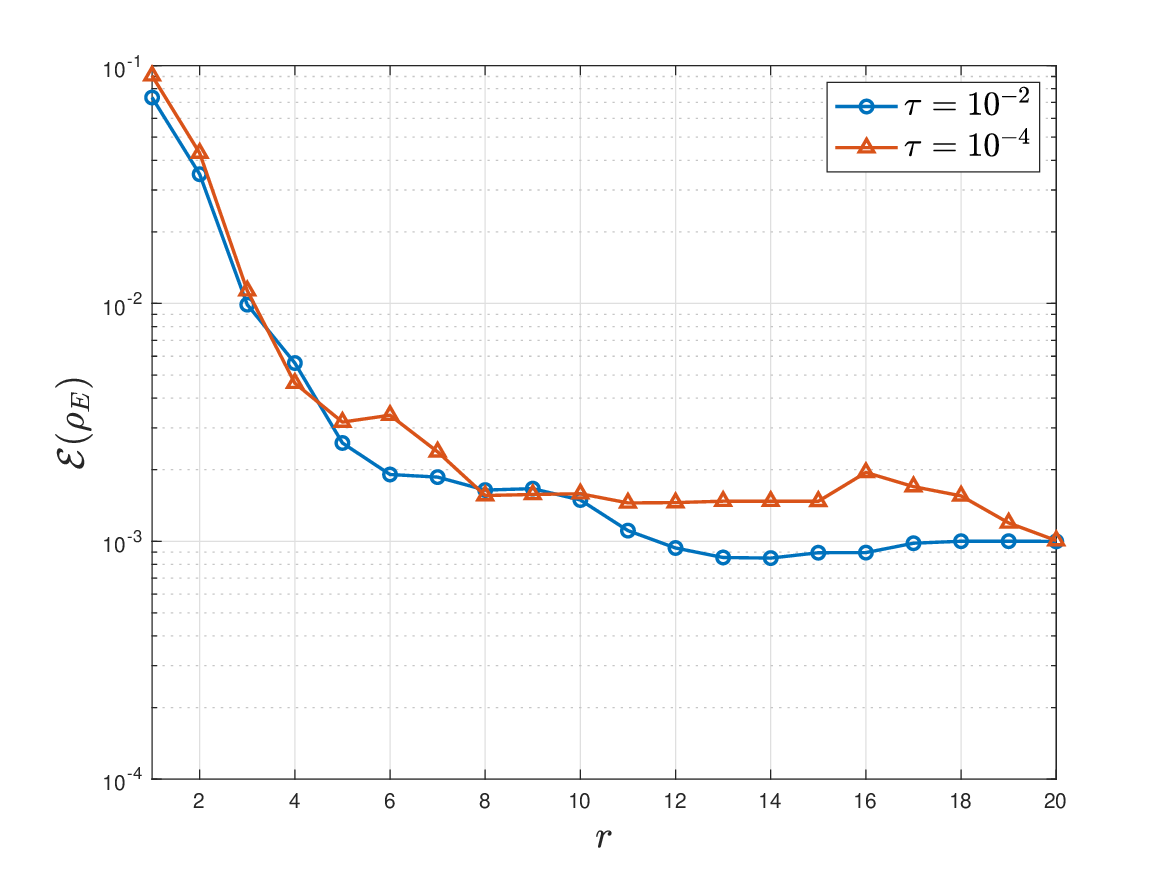}}

    \subfloat[Average $L^2$ error of tri-fidelity approximations for $\rho_I$]{\includegraphics[width=.45\textwidth]{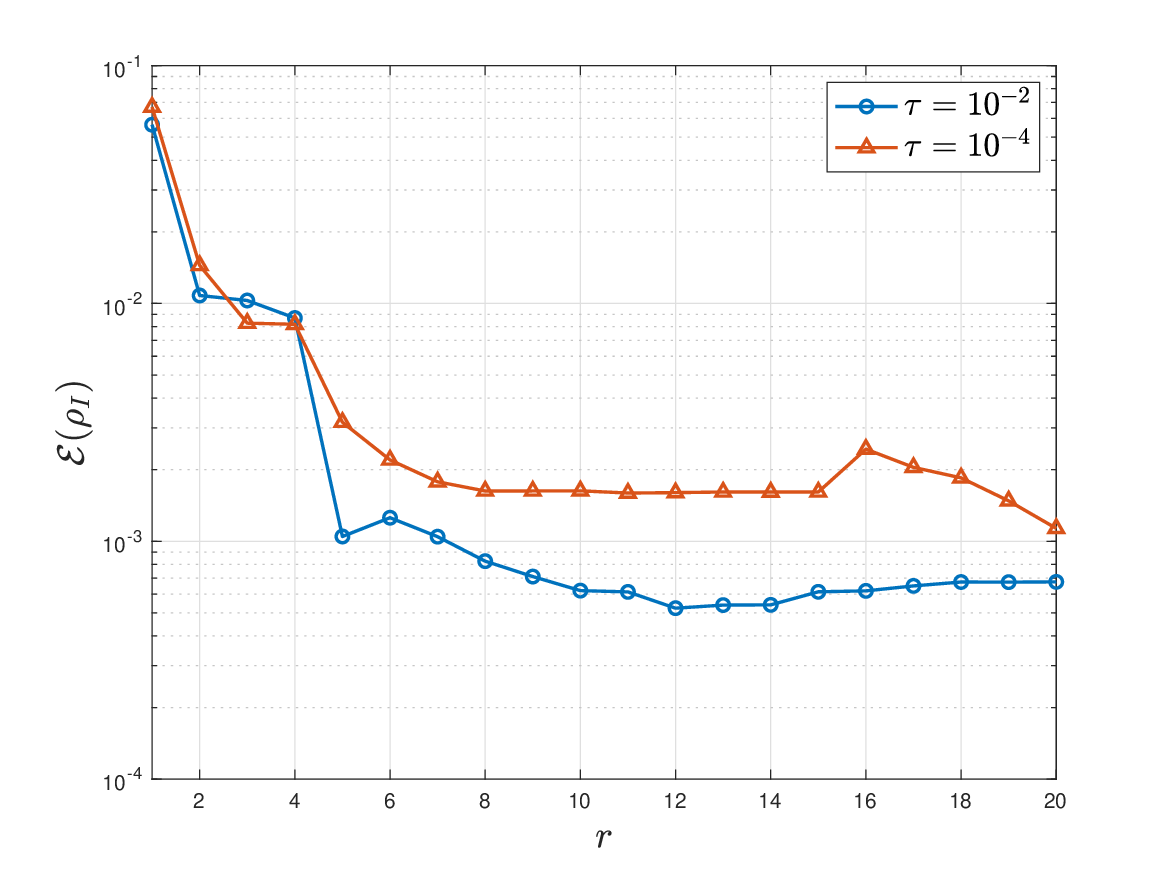}}
    \hspace{5mm}
    \subfloat[Average $L^2$ error of tri-fidelity approximations for $\rho_R$]{\includegraphics[width=.45\textwidth]{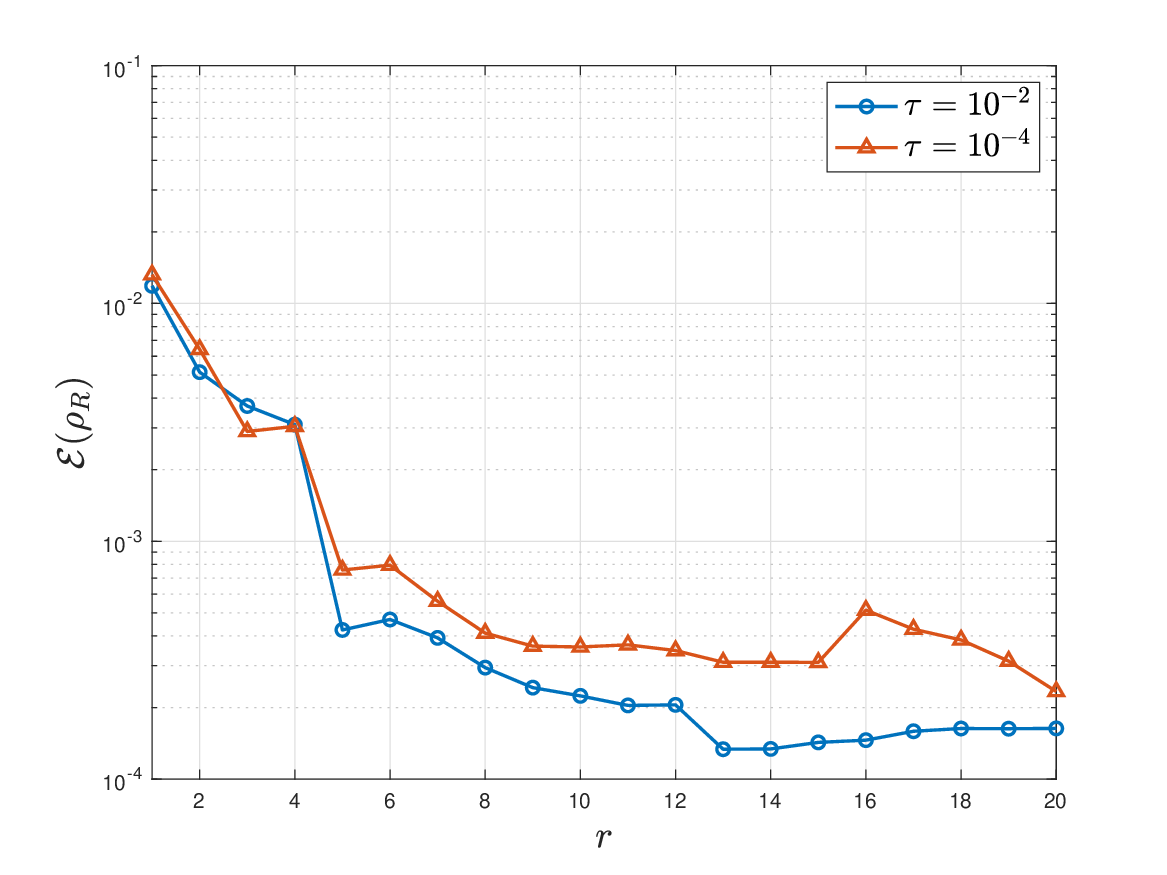}}
    \caption{Test 3: Average $L^2$ error of tri-fidelity approximations for $\rho_J$ with respect to the number of high-fidelity simulation runs at different $\tau$.}
    \label{Fig.Test3.1.m}
\end{figure*}

Figure~\ref{Fig.Test3.1.m} reports the average $L^2$ errors of the tri-fidelity approximation with respect to the number of high-fidelity runs. The errors decrease as more selected high-fidelity samples are used. This again confirms that the selected samples contain useful information about the high-fidelity solution manifold. These results show the performance of the tri-fidelity construction: the low-fidelity model is used only for the cheap selection of important points, while the medium-fidelity model provides a more reliable projection rule. Therefore, the tri-fidelity construction can improve the approximation without requiring the most expensive solver to be evaluated on the whole candidate set. 


\begin{figure*}
	\centering 
	\subfloat[$\mathbb{E}\big(\rho_S(t,\z)\big)$ when $\tau=10^{-2}$]{\includegraphics[width=.4\textwidth]{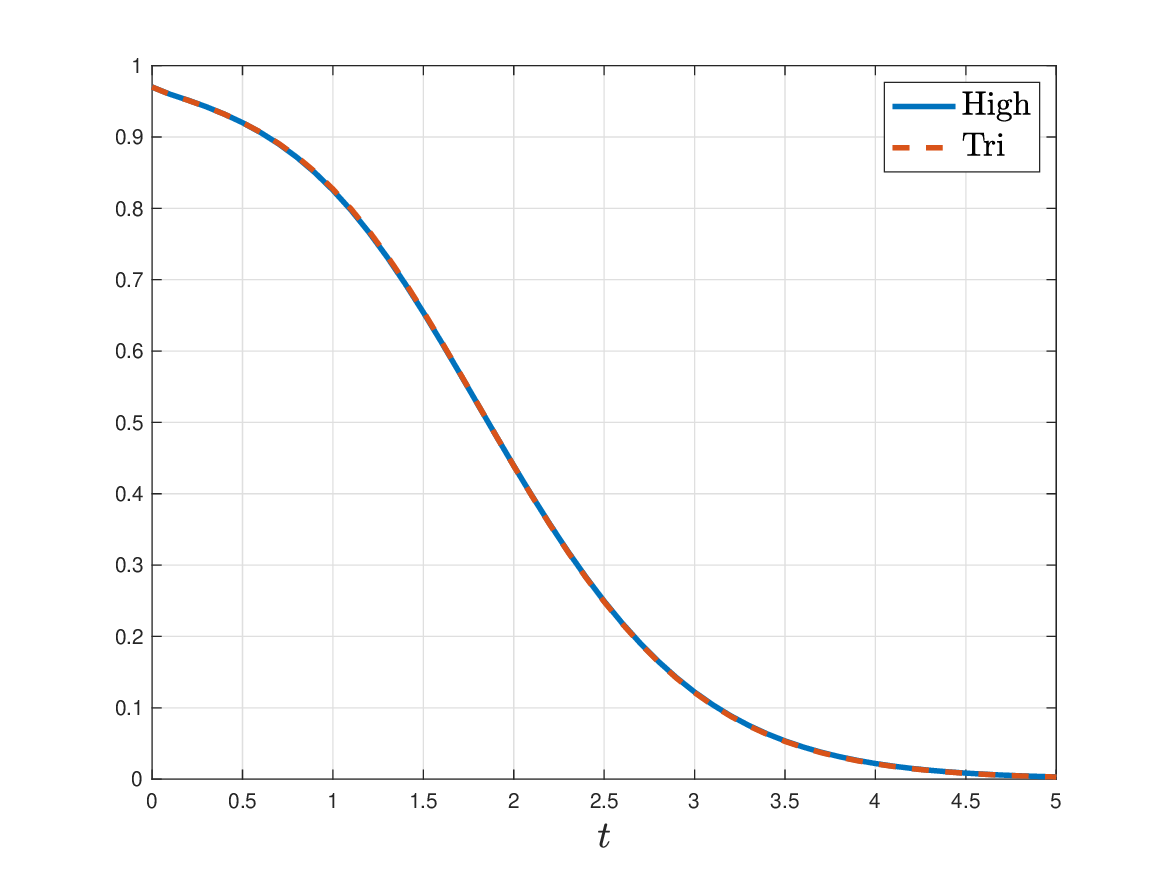}}
    \hspace{5mm}
    \subfloat[SD$\big(\rho_S(t,\z)\big)$ when $\tau=10^{-2}$]{\includegraphics[width=.4\textwidth]{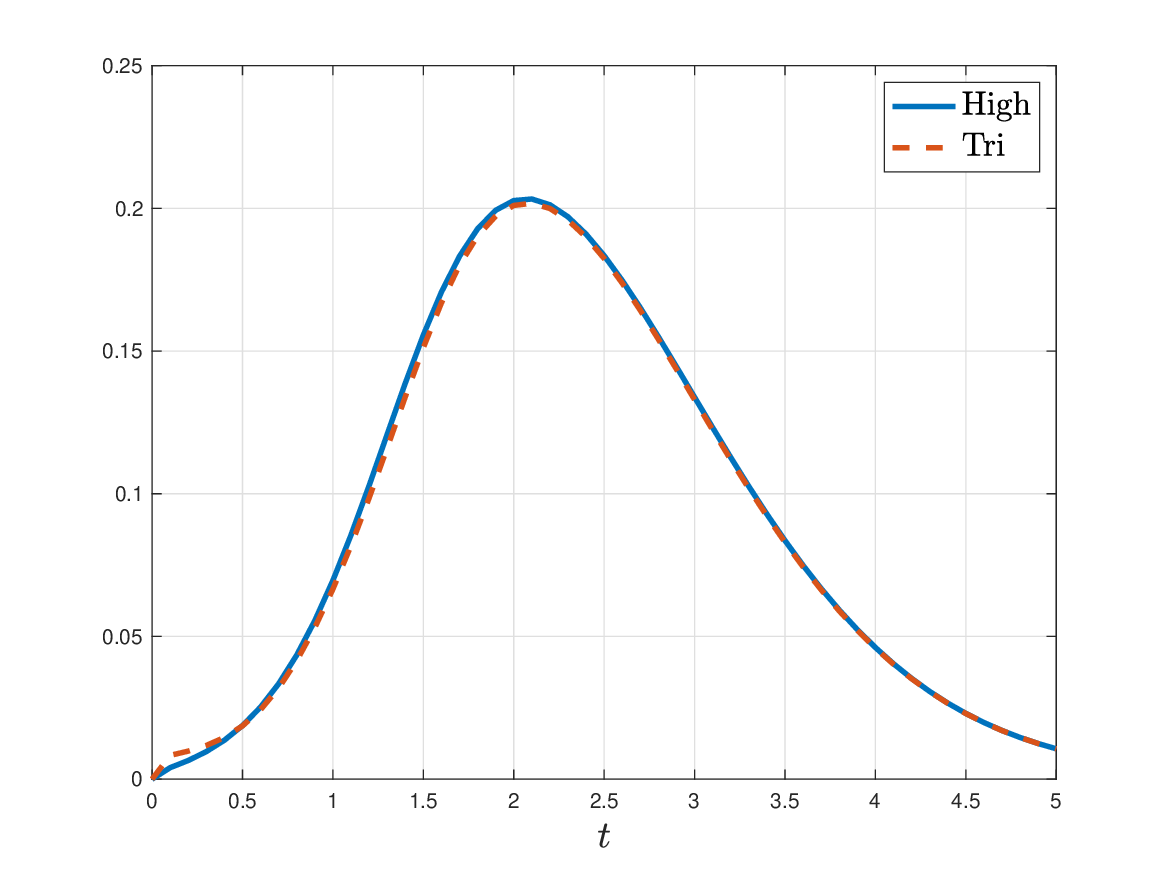}}

    \subfloat[$\mathbb{E}\big(\rho_S(t,\z)\big)$ when $\tau=10^{-4}$]{\includegraphics[width=.4\textwidth]{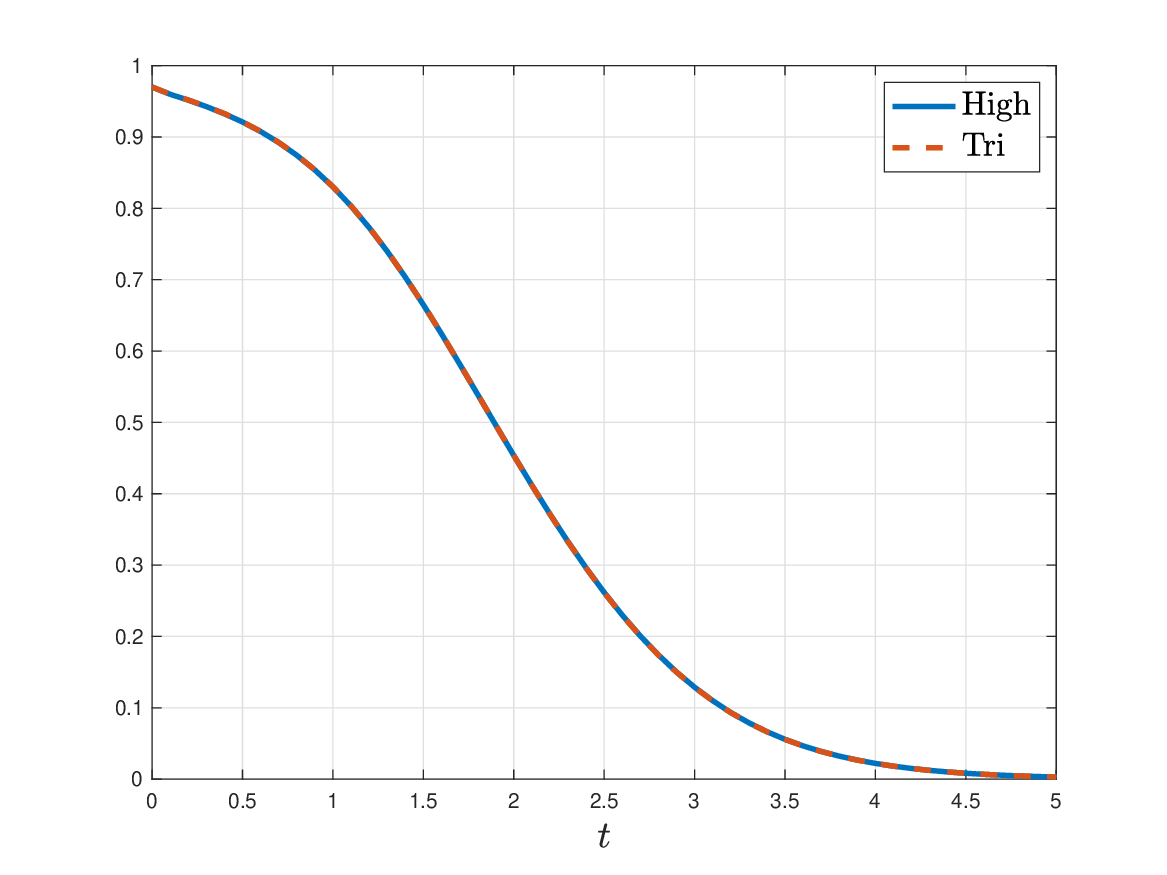}}
    \hspace{5mm}
    \subfloat[SD$\big(\rho_S(t,\z)\big)$ when $\tau=10^{-4}$]{\includegraphics[width=.4\textwidth]{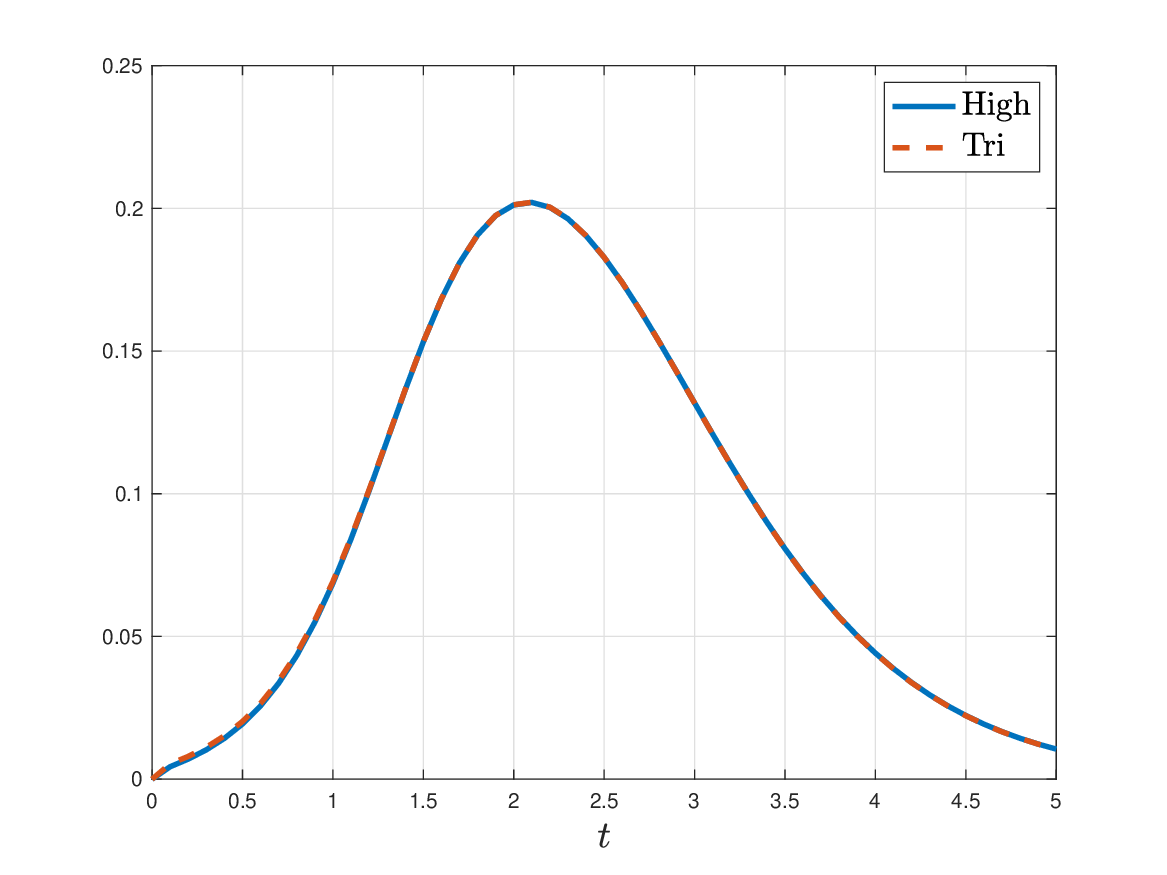}}
    
	\caption{Test 3: Mean and standard deviation of tri-fidelity solutions of $\rho_S(t,\z)$ at different $\tau$. } 
	\label{Fig.Test3.2.m}
\end{figure*}

Figure~\ref{Fig.Test3.2.m} further compares the mean and standard deviation of the tri-fidelity approximation with the high-fidelity reference. The agreement of these statistical quantities shows that the tri-fidelity method can reproduce not only samplewise trajectories but also the main uncertainty information. This is especially relevant in the present test because the uncertainty range crosses $\theta=0$, where the contact equilibrium and the macroscopic approximation are more delicate.

\begin{figure*}
    \centering
    \subfloat[Solution graph of $\rho_S$ at a certain $\z$]{\includegraphics[width=.45\textwidth]{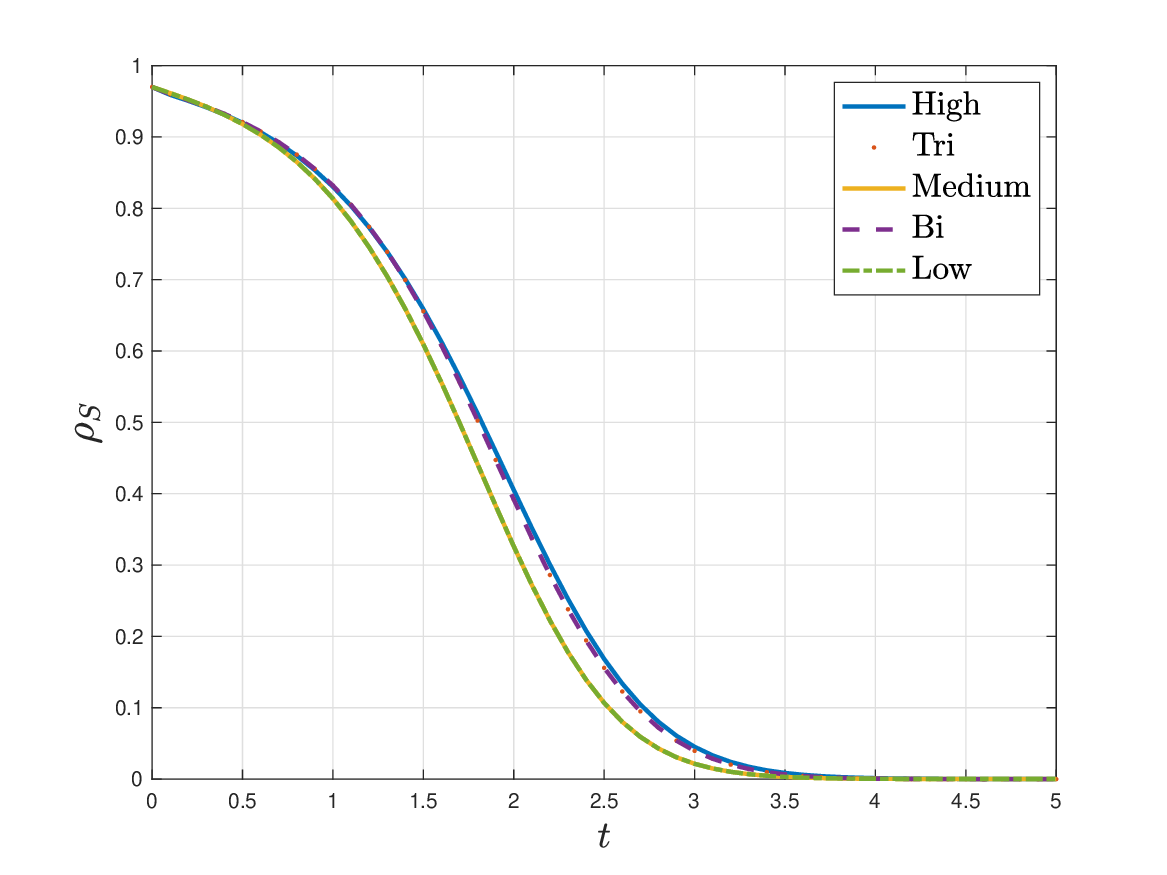}}
    \hspace{5mm}
    \subfloat[Solution graph of $\rho_E$ at a certain $\z$]{\includegraphics[width=.45\textwidth]{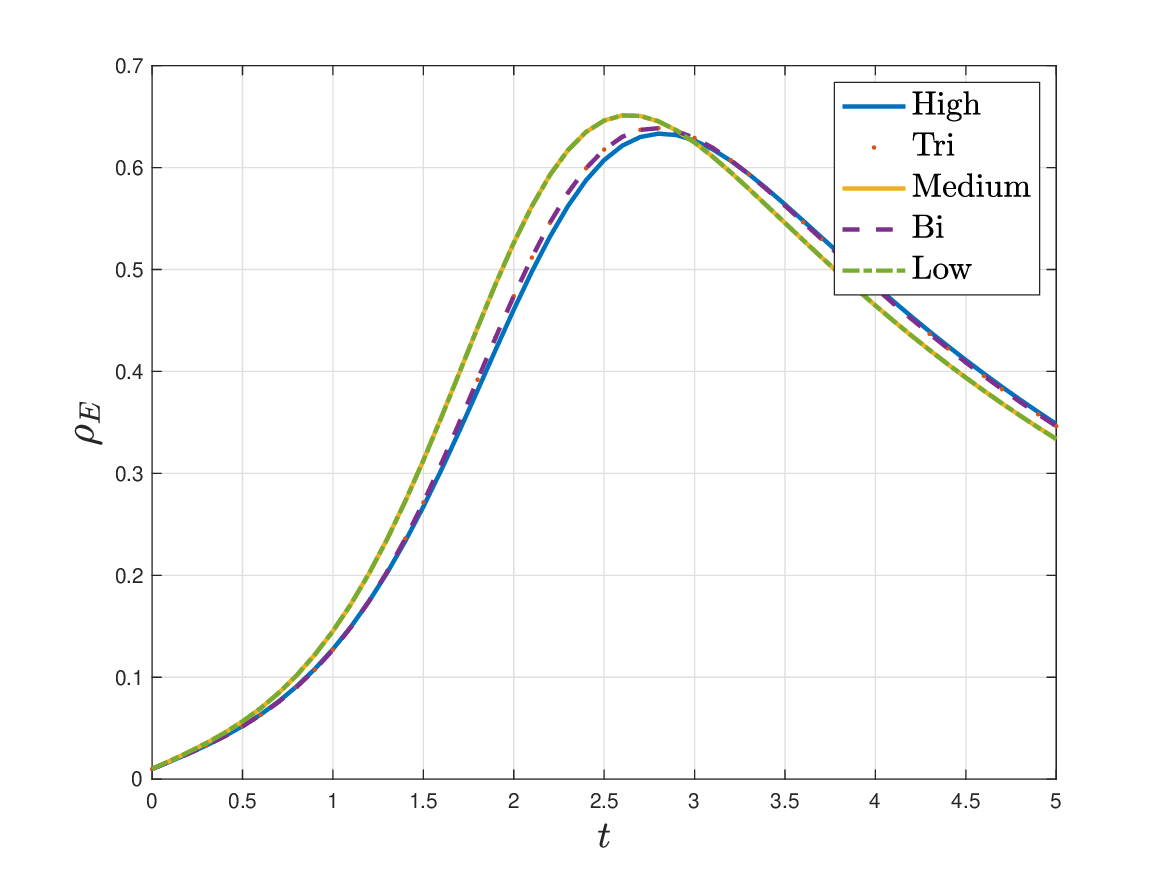}}

    \subfloat[Solution graph of $\rho_I$ at a certain $\z$]{\includegraphics[width=.45\textwidth]{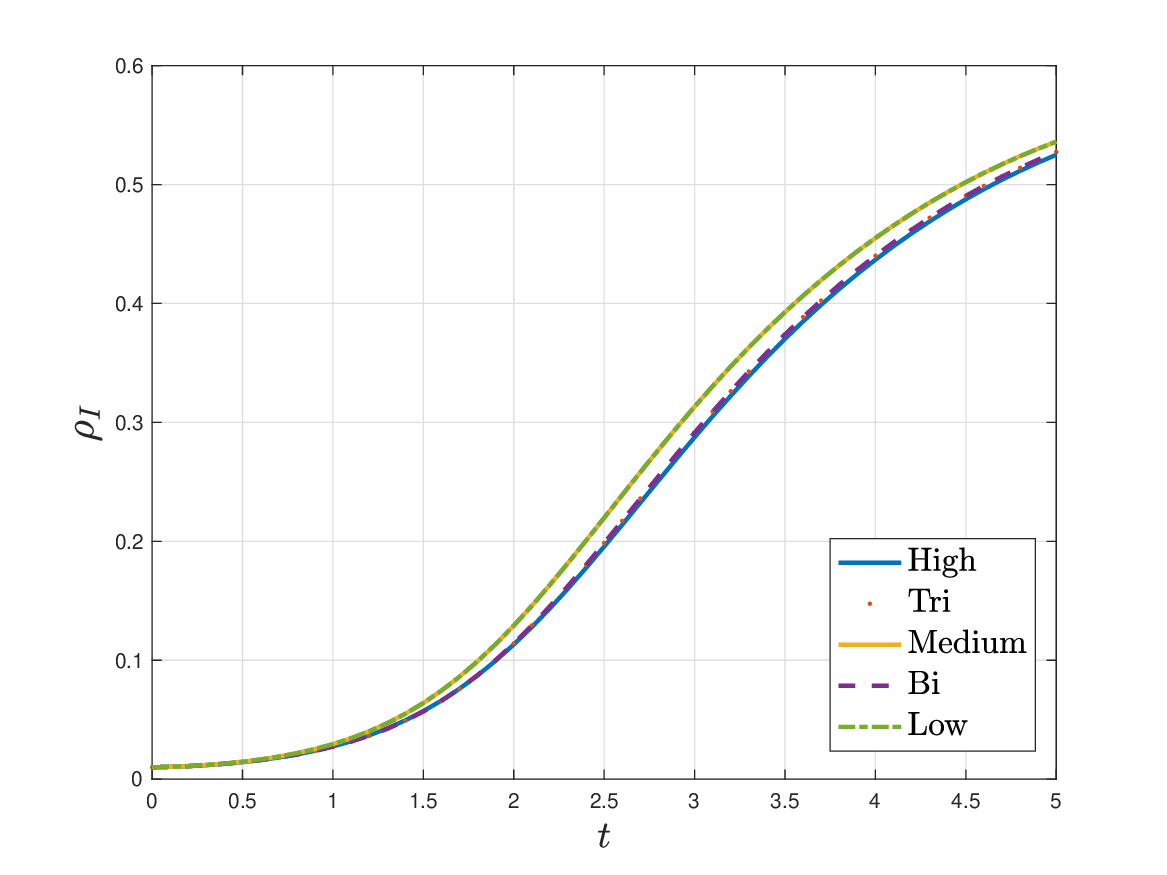}}
    \hspace{5mm}
    \subfloat[Solution graph of $\rho_R$ at a certain $\z$]{\includegraphics[width=.45\textwidth]{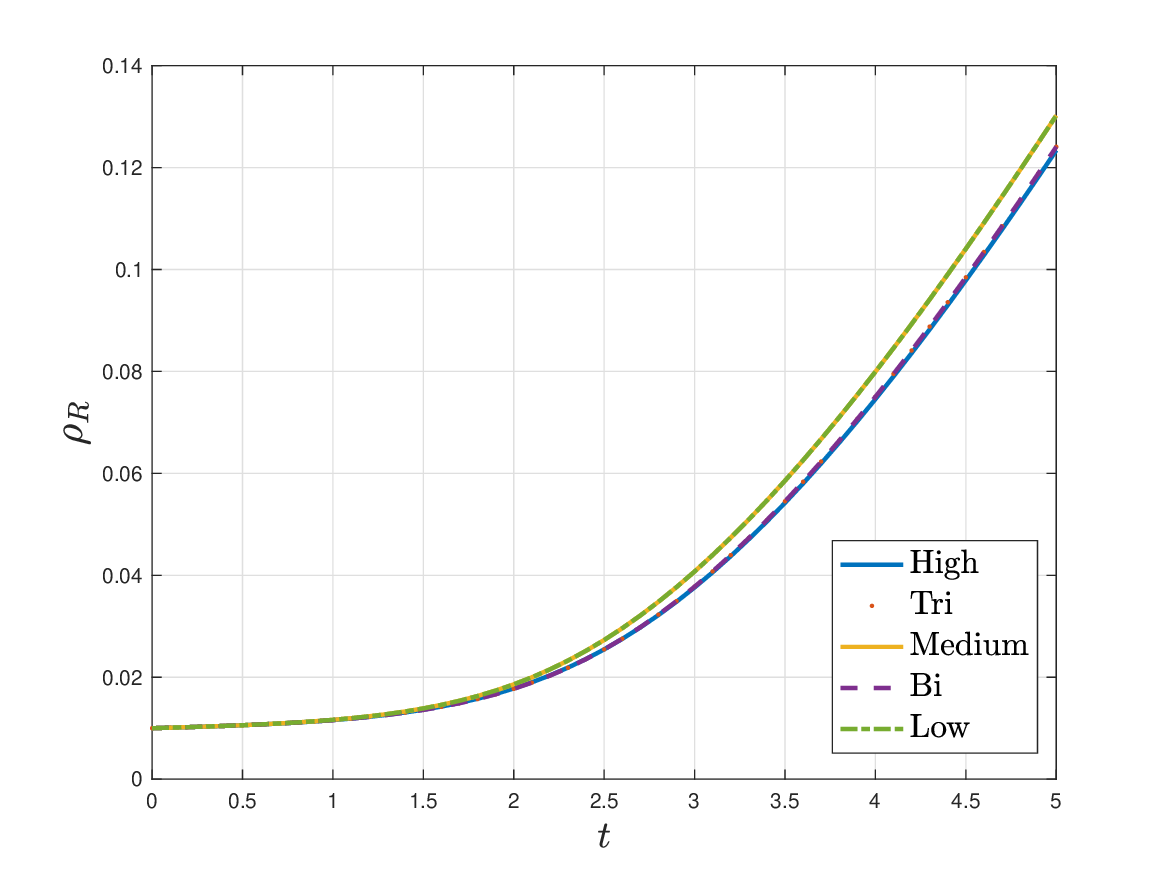}}
    \caption{Test 3: Solution graphs of $\rho_J$ at a certain $\z$ with $\tau=10^{-2}$. In Test 3, the bi-fidelity solution uses the low-fidelity and high-fidelity solvers. }
    \label{Fig.Test3.3.m}
\end{figure*}

Figure~\ref{Fig.Test3.3.m} shows representative time evolutions for $\rho_J$. The low-fidelity and medium-fidelity solvers capture the general temporal behavior, but they still show visible discrepancies from the high-fidelity solution. By contrast, the bi-fidelity and tri-fidelity approximations are much closer to the high-fidelity curves. This observation supports the main idea of the multi-fidelity approach: even when cheaper solvers are not accurate enough to replace the microscopic model directly, they can still guide the construction of an accurate surrogate when combined with a small number of high-fidelity simulations.

Overall, Test 3 shows that the tri-fidelity method is a useful extension when an intermediate model is available. The cheapest model reduces the cost of selecting representative samples, the medium-fidelity model improves the projection coefficients, and the high-fidelity model supplies the accurate snapshots. This division of roles is particularly suitable for kinetic epidemic models, where fully resolved simulations are expensive but several reduced models are available.

\section{Conclusion and Future Work}

In this paper, we proposed a multi-fidelity method for kinetic epidemic models with uncertain contact dynamics. The method uses the hierarchy between kinetic and macroscopic models. The high-fidelity kinetic solver gives accurate results, but it is expensive. The low- and medium-fidelity solvers are cheaper and are used to select important samples or compute projection coefficients. The numerical results show that the proposed methods can approximate the high-fidelity solution well with only a small number of high-fidelity runs. In the bi-fidelity tests, the method gives accurate results for both the special cases $\theta=\pm 1$ and the more general case $\theta \in [-1,1]$. It also captures the mean and standard deviation of the main epidemic quantities. In the tri-fidelity test, the use of a medium-fidelity model improves the projection step and gives reliable approximations. 

Beyond the specific model considered here, the proposed framework is general enough to accommodate a broad class of kinetic descriptions, and its structure naturally lends itself to several promising extensions. These include the incorporation of network-based contact models, the integration with real epidemiological data for data-driven uncertainty quantification, and the application to optimal control problems aimed at assessing the effects of intervention strategies on epidemic dynamics. 

The proposed methods reduce the computational cost of uncertainty quantification while maintaining good accuracy for important epidemic observables. Taken together, these directions point toward a broader role for multi-fidelity kinetic methods in epidemic forecasting and public health decision-making under uncertainty.

\section*{Acknowledgements}
A.M. and M.Z. wish to acknowledge GNFM (National Group of Mathematical Physics)  of INdAM (National Institute of High Mathematics). M.Z. acknowledges partial support from the PRIN2022PNRR project No.P2022Z7ZAJ, European Union - NextGenerationEU and by ICSC - Centro Nazionale di Ricerca in High Performance
Computing, Big Data and Quantum Computing, funded by European Union - NextGenerationEU. A.M. acknowledges the support by Fondo Italiano per la Scienza (FIS2023-01334) advanced
grant ”ADvanced numerical Approaches for MUltiscale Systems with uncertainties” -
ADAMUS.

\appendix
\section{Pivoted Cholesky sample selection}\label{app:chol-selection}

This appendix details the point-selection step used in the multi-fidelity construction. Given a candidate set $\Gamma_N=\{z_1,\ldots,z_N\}$, the low-fidelity snapshots are vectorized epidemic observables, for example time histories of $\rho_S$, $\rho_E$, or selected moments. In Algorithm~\ref{alg:chol-greedy-epidemic}, we use the notation $f^L_\ell:=\mathbf U^L(\z_\ell)$, $\ell\in\mathcal N=\{1,\ldots,N\}$, to keep the presentation close to the standard pivoted Cholesky selection procedure. The inner product $\langle\cdot,\cdot\rangle^L$ is the discrete inner product induced by the chosen observable vector.

\begin{algorithm}
    \caption{Pivoted Cholesky selection of representative epidemic samples}
    \label{alg:chol-greedy-epidemic}
    \begin{algorithmic}[1]
        \STATE \textbf{Input:} ensemble $\mathcal F^L=\{f^L_1,\ldots,f^L_N\}$, where $f^L_\ell=\mathbf U^L(\z_\ell)$ for $\ell\in\mathcal N$, maximum number $M$ of samples to be selected, and tolerance $\delta>0$.
        \STATE Initialize the index vector $\gamma=[1,\ldots,M]$ and the inner-product matrix $L=0_{N\times M}$.
        \STATE Initialize the diagonal values $w$ of the Gramian matrix of $\mathcal F^L$:
        \[
            w_\ell=\|f^L_\ell\|_{L}^{2}
            =
            \langle f^L_\ell,f^L_\ell\rangle^L,
            \qquad \ell=1,\ldots,N .
        \]
        \FOR{$n=1,\ldots,M$}
            \IF{$\max_{\ell\in\{n,\ldots,N\}} w_\ell<\delta$}
                \STATE $n\leftarrow n-1$ and stop the loop.
            \ENDIF
            \STATE Select the next pivot:
            \[
                \gamma(n)=\arg\max_{\ell\in\{n,\ldots,N\}} w_\ell .
            \]
            \STATE Exchange row $n$ and row $\gamma(n)$ in $L$, and exchange entry $n$ and entry $\gamma(n)$ in $\mathcal F^L$, $\Gamma_N$, and $w$.
            \STATE Compute, for $\ell=n+1,\ldots,N$,
            \[
                r_\ell
                =
                \langle f^L_\ell,f^L_n\rangle^L
                -
                \sum_{j=1}^{n-1}L_{\ell,j}L_{n,j}.
            \]
            \STATE Set
            \[
                L_{n,n}=\sqrt{w_n},
                \qquad
                L_{\ell,n}=\frac{r_\ell}{L_{n,n}},
                \quad \ell=n+1,\ldots,N .
            \]
            \STATE Update the residual diagonal:
            \[
                w_\ell\leftarrow w_\ell-L_{\ell,n}^2,
                \qquad \ell=n+1,\ldots,N .
            \]
        \ENDFOR
        \STATE Truncate rows $M+1,\ldots,N$ of $L$.
        \STATE \textbf{Output:} selected index set $\gamma$, corresponding sample set $\Gamma_\gamma=\{\z_k:k\in\gamma\}$, and Cholesky factor $L$ for $\{f^L_k:k\in\gamma\}$.
    \end{algorithmic}
\end{algorithm}


 


\bibliographystyle{plain}
\bibliography{Ref.bib}

@article{GF,
author = {Giambiagi-Ferrari, C. and Pinasco, J.P. and Saintier, N.},
title =  {Coupling epidemiological models with social
dynamics},
journal = {Bull. Math. Biol.},
volume = {83},
number = {74},
year = {2021}
}

@article{VIG_comput,
    author = {Viguerie, A. and Veneziani, A. and Lorenzo, G. and Baroli, D. and Aretz-Nellesen, N. and Patton, A. and Yankeelov, T. E. and Reali, A. and Hughes, T. J. R. and Auricchio, F. },
    title = {Diffusion–reaction compartmental models formulated in a continuum mechanics framework: application to {COVID}-19, mathematical analysis, and numerical study} ,
    journal = {Comput. Mech.},
    year = {2020},
    volume ={66},
pages = {1131--1152}
}

@article{colli,
title = { Chemotaxis-inspired PDE model for airborne infectious disease transmission: analysis and simulations}, 
author = {Colli, P. and Marinoschi, G. and Rocca, E. and Viguerie, A.},
journal={J. Nonlinear. Sci.},
volume ={35},
pages = {28},
year ={2025}
}

@article{viguerie22,
title = {Coupled and uncoupled dynamic mode decomposition in multi-compartmental systems with applications to epidemiological and additive manufacturing problems},
journal = {Comput. Methods Appl. Mech. Eng.},
volume = {391},
pages = {114600},
year = {2022},
issn = {0045-7825},
doi = {https://doi.org/10.1016/j.cma.2022.114600},
url = {https://www.sciencedirect.com/science/article/pii/S0045782522000202},
author = {Viguerie, A. and  Barros, G. F. and Grave, M. and Reali, A. and  Coutinho, A.L.G.A.}
}

@InProceedings{proc_MZ,
author="Medaglia, A.
and Zanella, M.",
editor="Barbante, Paolo
and Belgiorno, Francesco D.
and Lorenzani, Silvia
and Valdettaro, Lorenzo",
title="Kinetic and Macroscopic Epidemic Models in Presence of Multiple Heterogeneous Populations",
booktitle="From Kinetic Theory to Turbulence Modeling",
year="2023",
publisher="Springer Nature Singapore",
address="Singapore",
pages="191--201",
isbn="978-981-19-6462-6"
}

@article {FMZ,
    AUTHOR = {Franceschi, J. and Medaglia, A. and Zanella, M.},
     TITLE = {On the optimal control of kinetic epidemic models with
              uncertain social features},
   JOURNAL = {Optimal Control Appl. Methods},
  FJOURNAL = {Optimal Control Applications \& Methods},
    VOLUME = {45},
      YEAR = {2024},
    NUMBER = {2},
     PAGES = {494--522},
}

@article {DPTZ,
    AUTHOR = {Dimarco, G. and Perthame, B. and Toscani, G. and Zanella, M.},
     TITLE = {Kinetic models for epidemic dynamics with social
              heterogeneity},
   JOURNAL = {J. Math. Biol.},
  FJOURNAL = {Journal of Mathematical Biology},
    VOLUME = {83},
      YEAR = {2021},
    NUMBER = {1},
     PAGES = {Paper No. 4, 32},
}

@article {DTZ, 
    AUTHOR = {Dimarco, G. and Toscani, G. and Zanella, M.},
     TITLE = {Optimal control of epidemic spreading in the presence of
              social heterogeneity},
   JOURNAL = {Philos. Trans. Roy. Soc. A},
  FJOURNAL = {Philosophical Transactions of the Royal Society A.
              Mathematical, Physical and Engineering Sciences},
    VOLUME = {380},
      YEAR = {2022},
    NUMBER = {2224},
     PAGES = {Paper No. 20210160, 16},
}

@article {DeVore,
    AUTHOR = {DeVore, R. and Petrova, G. and Wojtaszczyk, P.},
     TITLE = {Greedy algorithms for reduced bases in {B}anach spaces},
   JOURNAL = {Constr. Approx.},
  FJOURNAL = {Constructive Approximation. An International Journal for
              Approximations and Expansions},
    VOLUME = {37},
      YEAR = {2013},
    NUMBER = {3},
     PAGES = {455--466},
}

@article {ZNX14,
    AUTHOR = {Zhu, X. and Narayan, A. and Xiu, D.},
     TITLE = {Computational aspects of stochastic collocation with
              multifidelity models},
   JOURNAL = {SIAM/ASA J. Uncertain. Quantif.},
  FJOURNAL = {SIAM/ASA Journal on Uncertainty Quantification},
    VOLUME = {2},
      YEAR = {2014},
    NUMBER = {1},
     PAGES = {444--463},
}

@article {LL25,
    AUTHOR = {Lin, Y. and Liu, L.},
     TITLE = {On a class of multi-fidelity methods for the semiclassical {S}chr\"odinger equation with uncertainties},
   JOURNAL = {SIAM J. Sci. Comput.},
  FJOURNAL = {SIAM Journal on Scientific Computing},
    VOLUME = {47},
    NUMBER = {5},
      YEAR = {2025},
}

@article{Liu2020,
title = {A bi-fidelity method for the multiscale {B}oltzmann equation with random parameters},
Journal = {J. Comput. Phys.},
FJournal = {Journal of Computational Physics},
volume = {402},
pages = {108914},
year = {2020},
issn = {0021-9991},
author = {Liu, L. and Zhu, X.},
}

@article{Dimarco2020,
	author = {Dimarco, G. and Pareschi, L. and Toscani, G. and Zanella, M.},
	journal = {Phys. Rev. E},
	number = {022303},
	title = {Weatlh distribution under the spread of infectious diseases},
	volume = {102},
	year = {2020}
}

@article{Albi2021,
author = {Albi, G. and Pareschi, L. and Zanella, M.},
	journal = {J. Math. Biol.},
	number = {63},
	title = {Control with uncertain data of socially structured compartmental epidemic models},
	volume = {82},
	year = {2021}
}

@article{Albi2021.b,
	author = {Albi, G. and Pareschi, L. and Zanella, M.},
	journal = {Math. Biosci. Eng.},
	number = {6},
	pages = {7161--7190},
	title = {Modelling lockdown measures in epidemic outbreaks using selective socio-economic containment with uncertainty},
	volume = {18},
	year = {2021}
}

@article{Gatto2020,
	author = {Gatto, M. and Bertuzzo, E. and Mari, L. and Miccoli, S. and Carraro, L. and Casagrandi, R. and Rinaldo, A.},
	doi = {10.1073/pnas.2004978117},
	issn = {0027-8424},
	journal = {Proc. Natl. Acad. Sci. U.S.A.},
	month = {May},
	number = {19},
	pages = {10484--10491},
	title = {Spread and dynamics of the {COVID}-19 epidemic in {I}taly: Effects of emergency containment measures},
	url = {https://europepmc.org/articles/PMC7229754},
	volume = {117},
	year = {2020},
}

@article{
Sun21,
author = {Sun, K.  and Wang, W.  and Gao, L.  and  Wang, Y.  and Luo, K.  and Ren, L.  and Zhan, Z.  and  Chen, X.  and Zhao, S.  and Huang, Y.  and Sun, Q.  and Liu, Z.  and Litvinova, M.  and Vespignani, A.  and Ajelli, M.  and Viboud, C.  and Yu, H. },
title = {Transmission heterogeneities, kinetics, and controllability of {SARS}-{C}o{V}-2},
journal = {Science},
volume = {371},
number = {6526},
pages = {eabe2424},
year = {2021},
doi = {10.1126/science.abe2424}
}

@article{martalo26,
    author = {Martalò, G. and Toscani, G. and Zanella, M.},
    title = {Individual-based foundation of SIR-type epidemic models: mean-field limit and large-time behaviour},
    journal = {Proc. R. Soc. A},
    volume = {482},
    number = {2331},
    pages = {20250633},
    year = {2026},
    month = {02}
}

@article{mossong08,
    doi = {10.1371/journal.pmed.0050074},
    author = {Mossong, J. AND Hens, N. AND Jit, M. AND Beutels, P. AND Auranen, K. AND Mikolajczyk, R. AND Massari, M. AND Salmaso, S. AND Tomba, G. S.  AND Wallinga, J. AND Heijne, J. AND Sadkowska-Todys, M. AND Rosinska, M. AND Edmunds, W. J.},
    journal = {PLOS Medicine},
    publisher = {Public Library of Science},
    title = {Social Contacts and Mixing Patterns Relevant to the Spread of Infectious Diseases},
    year = {2008},
    month = {03},
    volume = {5},
    url = {https://doi.org/10.1371/journal.pmed.0050074},
    pages = {1-1},
    number = {3},
}

@article{visi226,
title = {Age Group Sensitivity Analysis in age stratified epidemic models: Investigating the impact of contact matrix structure},
journal = {Epidemics},
volume = {55},
pages = {100915},
year = {2026},
issn = {1755-4365},
doi = {https://doi.org/10.1016/j.epidem.2026.100915},
url = {https://www.sciencedirect.com/science/article/pii/S1755436526000319},
author = { Vizi, Z. and Korir, E. K.  and Bogya, N. and Rosztóczy, C. and Kökény, Z. and Makay, G. and  Boldog, P.},
}

@article{fumanelli,
    doi = {10.1371/journal.pcbi.1002673},
    author = {Fumanelli, L. AND Ajelli, M. AND Manfredi, P. AND Vespignani, A. AND Merler, S.},
    journal = {PLOS Computational Biology},
    publisher = {Public Library of Science},
    title = {Inferring the Structure of Social Contacts from Demographic Data in the Analysis of Infectious Diseases Spread},
    year = {2012},
    month = {09},
    volume = {8},
    url = {https://doi.org/10.1371/journal.pcbi.1002673},
    pages = {1-10},
    number = {9}
}

@article{
britton20,
author = {Britton, T.  and Ball, F.  and Trapman, P.},
title = {A mathematical model reveals the influence of population heterogeneity on herd immunity to SARS-CoV-2},
journal = {Science},
volume = {369},
number = {6505},
pages = {846-849},
year = {2020},
doi = {10.1126/science.abc6810},
URL = {https://www.science.org/doi/abs/10.1126/science.abc6810},
eprint = {https://www.science.org/doi/pdf/10.1126/science.abc6810}
}

@article{beraud15,
    doi = {10.1371/journal.pone.0133203},
    author = {Béraud, G. AND Kazmercziak, S. AND Beutels, P. AND Levy-Bruhl, D. AND Lenne, X. AND Mielcarek, N. AND Yazdanpanah, Y. AND Boëlle, P.-Y. AND Hens, N. AND Dervaux, B.},
    journal = {PLOS ONE},
    publisher = {Public Library of Science},
    title = {The French Connection: The First Large Population-Based Contact Survey in France Relevant for the Spread of Infectious Diseases},
    year = {2015},
    month = {07},
    volume = {10},
    url = {https://doi.org/10.1371/journal.pone.0133203},
    pages = {1-22},
    number = {7},
}

@article{BARTHELEMY2005275,
title = {Dynamical patterns of epidemic outbreaks in complex heterogeneous networks},
journal = {J. Theoret. Biol.},
volume = {235},
number = {2},
pages = {275-288},
year = {2005},
issn = {0022-5193},
doi = {https://doi.org/10.1016/j.jtbi.2005.01.011},
url = {https://www.sciencedirect.com/science/article/pii/S0022519305000251},
author = { Barthélemy, M. and  Barrat, A. and  Pastor-Satorras, R. and Vespignani, A.},
}

@article{block,
	author = {Block, P. and Hoffman, M. and Raabe, I. J. and Dowd, J. B. and Rahal, C. and Kashyap, R. and Mills, M. C.},
	date = {2020/06/01},
	doi = {10.1038/s41562-020-0898-6},
	id = {Block2020},
	isbn = {2397-3374},
	journal = {Nature Human Behaviour},
	number = {6},
	pages = {588--596},
	title = {Social network-based distancing strategies to flatten the COVID-19 curve in a post-lockdown world},
	url = {https://doi.org/10.1038/s41562-020-0898-6},
	volume = {4},
	year = {2020}}

@incollection{zanella_rev,
    author = {Zanella, M.},
    title = {Derivation of macroscopic epidemic models from multi-agent systems},
    booktitle = {Modeling, Analysis, and Control of Multi-Agent Systems Across Scales},
    publisher = {EMS Series of Congress Reports},
    year = {2026}
}

@Inbook{Albi2022,
author="Albi, G.
and Bertaglia, G.
and Boscheri, W.
and Dimarco, G.
and Pareschi, L.
and Toscani, G.
and Zanella, M.",
editor="Bellomo, Nicola
and Chaplain, Mark A. J.",
title="Kinetic Modelling of Epidemic Dynamics: Social Contacts, Control with Uncertain Data, and Multiscale Spatial Dynamics",
bookTitle="Predicting Pandemics in a Globally Connected World, Volume 1: Toward a Multiscale, Multidisciplinary Framework through Modeling and Simulation",
year="2022",
publisher="Springer International Publishing",
address="Cham",
pages="43--108",
}

@article{zanella_medaglia,
author = {Zanella, M. and Medaglia, A.},
title = {Control of overpopulated tails in kinetic epidemic models},
journal = {J. Hyperbolic Differ. Equ.},
volume = {23},
number = {01},
pages = {151-177},
year = {2026},
doi = {10.1142/S0219891626400072}
}

@article{LPZ2022,
title = {A bi-fidelity stochastic collocation method for transport equations with diffusive scaling and multi-dimensional random inputs},
journal = {J. Comput. Phys.},
volume = {462},
pages = {111252},
year = {2022},
issn = {0021-9991},
doi = {https://doi.org/10.1016/j.jcp.2022.111252},
author = { Liu, L. and Pareschi, L. and Zhu, X.}
}

@article{BLPZ2022,
title = {Bi-fidelity stochastic collocation methods for epidemic transport models with uncertainties},
journal = {Networks and Heterogeneous Media},
volume = {17},
number = {3},
pages = {401-425},
year = {2022},
issn = {1556-1801},
doi = {10.3934/nhm.2022013},
author = { Bertaglia, G. and Liu, L. and Pareschi, L. and Zhu, X.}
}

@article{DLPZ2021,
  title = {MULTI-FIDELITY METHODS FOR UNCERTAINTY PROPAGATION IN KINETIC EQUATIONS},
  author = {Dimarco, G. and Liu, L. and Pareschi, L. and Zhu, X.},
  JOURNAL = {{Panoramas et Synth{\`e}ses}},
  PUBLISHER = {{SMF}},
  YEAR = {2021}
}

@article{JLZC2026,
author = { Jin, X. and Liu, L. and Zhong, X. and Chung, E. T.},
title = {Efficient Numerical Method for the {S}chr\"{o}dinger Equation with High-Contrast Potentials},
year = {2026},
issue_date = {Dec 2025},
publisher = {Society for Industrial and Applied Mathematics},
address = {USA},
volume = {23},
number = {4},
issn = {1540-3459},
url = {https://doi.org/10.1137/25M173288X},
doi = {10.1137/25M173288X},
journal = {Multiscale Model. Simul.},
month = dec,
pages = {1581–1606},
numpages = {26}
}

@article{ZPPRV2025,
title = {A model learning framework for inferring the dynamics of transmission rate depending on exogenous variables for epidemic forecasts},
journal = {Comput. Meth. Appl. Mech. Eng.},
volume = {437},
pages = {117796},
year = {2025},
issn = {0045-7825},
doi = {https://doi.org/10.1016/j.cma.2025.117796},
author = {Ziarelli, G. and Pagani, S. and Parolini, N. and Regazzoni, F. and Verani, M.}
}

\end{document}